\documentclass[11pt]{article}

\usepackage{latexsym,amsmath,amsfonts,amssymb,amstext,theorem,euscript,array,enumerate,mathrsfs}
\usepackage{color}
\usepackage{comment}

\newtheorem{Theorem}{Theorem}[section]
\newtheorem{Definition}{Definition}[section]
\newtheorem{Proposition}{Proposition}[section]

\newtheorem{Lemma}{Lemma}[section]
\newtheorem{Corollary}{Corollary}[section]
\newtheorem{Remark}{Remark}[section]

\numberwithin{equation}{section}

\def\esssup_#1{\underset{#1}{\mathrm{ess\,sup\, }}}
\def\essinf_#1{\underset{#1}{\mathrm{ess\,inf\, }}}

\def \trans{^{\scriptscriptstyle{\intercal}}}

\def\sqr#1#2{{\vcenter{\vbox{\hrule height .#2pt \hbox{\vrule
 width .#2pt height#1pt \kern#1pt \vrule
width .#2pt} \hrule height .#2pt}}}}

\def\ds{\begin{displaystyle}}
\def\eds{\end{displaystyle}}
\def\dis{\displaystyle }
\def\<{\langle }
\def\>{\rangle }

\def \N{\mathbb{N}}
\def \R{\mathbb{R}}

\def \E{\mathbb{E}}
\def \F{\mathbb{F}}

\def \P{\mathbb{P}}
\def \Q{\mathbb{Q}}

\def \Ac{{\cal A}}
\def \Bc{{\cal B}}

\def \Ec{{\cal E}}
\def \Fc{{\cal F}}
\def \Gc{{\cal G}}

\def \Ic{{\cal I}}
\def \Kc{{\cal K}}

\def \Pc{{\cal P}}

\def \Nc{{\cal N}}
\def \Rc{{\cal R}}
\def \Sc{{\cal S}}

\def \Vc{{\cal V}}
\def \Wc{{\cal W}}

\def \Vc{{\cal V}}

\def\calf{{\cal F}}

\def \eps{\varepsilon}

\def \ep{\hbox{ }\hfill$\Box$}

\def\reff#1{{\rm(\ref{#1})}}

\def\beqs{\begin{eqnarray*}}
\def\enqs{\end{eqnarray*}}
\def\beq{\begin{eqnarray}}
\def\enq{\end{eqnarray}}

\allowdisplaybreaks

\begin{document}

\title{Randomized dynamic programming principle and Feynman-Kac representation for optimal control of McKean-Vlasov dynamics}

\author{Erhan BAYRAKTAR\thanks{Department of Mathematics, University of Michigan; e-mail: \texttt{erhan@umich.edu}. E. Bayraktar is supported in part by the National Science Foundation under grant  DMS-1613170 and the Susan M. Smith Professorship.}
\and
Andrea COSSO\thanks{Politecnico di Milano, Dipartimento di Matematica, via Bonardi 9, 20133 Milano, Italy; e-mail: \texttt{andrea.cosso@polimi.it}}
\and
Huy\^{e}n PHAM\thanks{Laboratoire de Probabilit\'es et Mod\`eles Al\'eatoires, CNRS, UMR 7599, Universit\'e Paris Diderot, and CREST-ENSAE; e-mail: \texttt{pham@math.univ-paris-diderot.fr}. 
H. Pham  is supported in  part by  the ANR project CAESARS (ANR-15-CE05-0024).
}
}

\maketitle

\begin{abstract}
We analyze a stochastic optimal control problem, where the state process follows a McKean-Vlasov dynamics and the diffusion coefficient can be degenerate.  We prove that its value function $V$ admits a nonlinear Feynman-Kac representation in terms of a class of forward-backward stochastic differential equations, with an autonomous forward process. 
We exploit this probabilistic representation to rigorously prove the dynamic programming principle (DPP) for $V$. The Feynman-Kac representation we obtain has an important role beyond its intermediary role in obtaining our main result: in fact it would be useful in developing probabilistic numerical schemes for $V$. The DPP is important in obtaining a characterization of the value function as a solution of a non-linear partial differential equation (the so-called Hamilton-Jacobi-Belman equation), in this case on the Wasserstein space of measures. We should note that the usual way of solving these equations is through the Pontryagin maximum principle, which requires some convexity assumptions. There were attempts in using the dynamic programming approach before, but these works assumed a priori that the controls were of Markovian feedback type, which helps write the problem only in terms of the distribution of the state process (and the control problem becomes a deterministic problem). In this paper, we will consider open-loop controls and derive the dynamic programming principle in this most general case. In order to obtain the Feynman-Kac representation and the randomized dynamic programming principle, we implement the so-called randomization method, which consists in formulating a new McKean-Vlasov control problem, expressed in weak form taking the supremum over a family of equivalent probability measures. One of the main results of the paper is the proof that this latter control problem has the same value function $V$ of the original control problem. \end{abstract}

\vspace{5mm}


\vspace{5mm}


\maketitle

\date{}

\newpage

\section{Introduction}
\label{S:Introduction}

In the present paper we study a stochastic optimal control problem of McKean-Vlasov type. More precisely, let $T>0$ be a finite time horizon, $(\Omega,\Fc,\P)$ a complete probability space, $B=(B_t)_{t\geq0}$ a $d$-dimensional Brownian motion defined on $(\Omega,\Fc,\P)$, $\F^B=(\Fc_t^B)_{t\geq 0}$ the $\P$-completion of the filtration generated by $B$, and $\Gc$ a sub-$\sigma$-algebra of $\Fc$ independent of $B$. 
Let also $\mathscr P_{\text{\tiny$2$}}(\R^n)$ denote the set of all probability measures on $(\R^n,\Bc(\R^n))$ with a finite second-order moment. We endow $\mathscr P_{\text{\tiny$2$}}(\R^n)$ with the 2-Wasserstein metric $\Wc_{\text{\tiny$2$}}$, and  assume that $\Gc$ is rich enough in the sense that $\mathscr P_{\text{\tiny$2$}}(\R^n)=\{\P_{\text{\tiny$\xi$}}\colon\xi\in L^2(\Omega,\Gc,\P;\R^n)\}$, where $\P_{\text{\tiny$\xi$}}$ denotes the law of $\xi$ under $\P$. Then, the controlled state equations are given by
\begin{align}
X_s^{t,\xi,\alpha} \ &= \ \xi + \int_t^s b\big(r,X_r^{t,\xi,\alpha},\P_{\text{\tiny$X_r^{t,\xi,\alpha}$}},\alpha_r\big)\,dr + \int_t^s \sigma\big(r,X_r^{t,\xi,\alpha},\P_{\text{\tiny$X_r^{t,\xi,\alpha}$}},\alpha_r\big)\,dB_r, \label{StateEq1_Intro} \\
X_s^{t,x,\xi,\alpha} \ &= \ x + \int_t^s b\big(r,X_r^{t,x,\xi,\alpha},\P_{\text{\tiny$X_r^{t,\xi,\alpha}$}},\alpha_r\big)\,dr + \int_t^s \sigma\big(r,X_r^{t,x,\xi,\alpha},\P_{\text{\tiny$X_r^{t,\xi,\alpha}$}},\alpha_r\big)\,dB_r, \label{StateEq2_Intro}
\end{align}
for all $s\in[t,T]$, where $(t,x,\xi)\in[0,T]\times\R^n\times L^2(\Omega,\Gc,\P;\R^n)$, and $\alpha$ is an admissible control process, namely an $\F^B$-progressive process $\alpha\colon\Omega\times[0,T]\rightarrow A$, with $A$ Polish space. We denote by $\Ac$ the set of admissible control processes. On the coefficients $b\colon[0,T]\times\R^n\times\mathscr P_{\text{\tiny$2$}}(\R^n)\times A\rightarrow\R^n$ and $\sigma\colon[0,T]\times\R^n\times\mathscr P_{\text{\tiny$2$}}(\R^n)\times A\rightarrow\R^{n\times d}$ we impose standard Lipschitz and linear growth conditions, which guarantee existence and uniqueness of a pair $(X_s^{t,\xi,\alpha},X_s^{t,x,\xi,\alpha})_{s\in[t,T]}$ of continuous $(\Fc_s^B\vee\Gc)_s$-adapted processes solution to equations \eqref{StateEq1_Intro}-\eqref{StateEq2_Intro}. Notice that $X^{t,x,\xi,\alpha}$ depends on $\xi$ only through its law $\pi:=\P_{\text{\tiny$\xi$}}$. Therefore, we define $X^{t,x,\pi,\alpha}:=X^{t,x,\xi,\alpha}$.

The control problem consists in maximizing over all admissible control processes $\alpha\in\Ac$ the following functional
\[
J(t,x,\pi,\alpha) \ = \ \E\bigg[\int_t^T f\big(s,X_s^{t,x,\pi,\alpha},\P_{\text{\tiny$X_s^{t,\xi,\alpha}$}},\alpha_s\big)\,ds + g\big(X_T^{t,x,\pi,\alpha},\P_{\text{\tiny$X_T^{t,\xi,\alpha}$}}\big)\bigg],
\]
for any $(t,x,\pi)\in[0,T]\times\R^n\times\mathscr P_{\text{\tiny$2$}}(\R^n)$, where $f\colon[0,T]\times\R^n\times\mathscr P_{\text{\tiny$2$}}(\R^n)\times A\rightarrow\R$ and $g\colon\R^n\times\mathscr P_{\text{\tiny$2$}}(\R^n)\rightarrow\R$ satisfy suitable continuity and growth conditions, see Assumptions {\bf (A1)} and {\bf (A2)}. We define the value function
\begin{equation}\label{Value_Intro}
V(t,x,\pi) \ = \ \sup_{\alpha\in\Ac} J(t,x,\pi,\alpha),
\end{equation}
for all $(t,x,\pi)\in[0,T]\times\R^n\times\mathscr P_{\text{\tiny$2$}}(\R^n)$. We will show in Proposition \ref{P:V_MKV=V} that the mapping $V$ is the disintegration of the value function 
\begin{equation}\label{V_MKV_App_Intro}
V_{\textup{\tiny MKV}}(t,\xi) \ = \ \sup_{\alpha\in\Ac_{\text{\tiny$\xi$}}} \E\bigg[\int_t^T f\big(s,X_s^{t,\xi,\alpha},\P_{\text{\tiny$X_s^{t,\xi,\alpha}$}}^{\text{\tiny$\xi$}},\alpha_s\big)\,ds + g\big(X_T^{t,\xi,\alpha},\P_{\text{\tiny$X_T^{t,\xi,\alpha}$}}^{\text{\tiny$\xi$}}\big)\bigg],
\end{equation}
for every $(t,\xi)\in[0,T]\times L^2(\Omega,\Gc,\P;\R^n)$, where $\Ac_{\text{\tiny$\xi$}}$ denotes the set of $A$-valued $(\Fc_s^B\vee\sigma(\xi))$-progressive processes, and 
$\P_{\text{\tiny$X_s^{t,\xi,\alpha}$}}^{\text{\tiny$\xi$}}$ denotes the regular conditional distribution of the random variable $X_s^{t,\xi,\alpha}\colon\Omega\rightarrow\R^n$ with respect to $\sigma(\xi)$. 
That is,
\beq \label{desin}
V_{\textup{\tiny MKV}}(t,\xi) \; = \; \int V(t,x,\pi) \pi(dx).  
\enq
Notice that at time $t$ $=$ $0$, when $\xi$ $=$ $x_0$ is a constant, then $V_{\textup{\tiny MKV}}(0,x_0)$ is the natural formulation of the McKean-Vlasov control problem as in \cite{CarmDel15}.

Optimal control of McKean-Vlasov dynamics  is a new type of stochastic control problem related to, but different from, what is well-known as mean field games (MFG), and which has attracted a surge of interest in the stochastic control community since the lectures by P.L. Lions at Coll\`ege de France, see \cite{lio12} and \cite{car12}, and the recent books \cite{benbook} and \cite{carbook}.
 Both of these problems describe equilibriums states of large population of weakly interacting symmetric players and we refer to \cite{CarmDelLach} for a discussion pointing out the differences between the two frameworks: In a nutshell MFGs describe Nash equilibrium in large populations and the optimal control of McKean-Vlasov dynamics describes the Pareto optimality, as heuristically shown in \cite{CarmDelLach}, and 
recently proved in \cite{Lacker}. As an example we mention the model of systemic risk due to \cite{CFS}, where, using our notation, $X^{t,\xi,\alpha}$ (as well as the auxiliary process $X^{t,x,\xi,\alpha}$) represents the $\log$-reserve of the representative bank, and $\alpha$ is the rate of borrowing/lending to a central bank.

In the literature McKean-Vlasov control problem  is tackled  by two different approaches: On the one hand,  the  stochastic Pontryagin maximum principle allows one to characterize solutions to the controlled McKean-Vlasov systems in terms of an adjoint backward stochastic differential equation (BSDE) coupled with a forward SDE:  see \cite{anddje10}, \cite{bucetal11} in which the state dynamics depend upon moments of the distribution, and \cite{CarmDel15} for a deep  investigation in a more general setting.     
On the other hand,  the dynamic programming (DP) method (also called Bellman principle), which is known to be a powerful tool for standard Markovian stochastic control problem 
and  does not require any convexity assumption usually imposed in Pontryagin principle, was first used in \cite{LaurPironn14} and \cite{BensFY} for a specific  McKean-Vlasov SDE and cost functional, depending only upon statistics  
like the mean of the distribution of the state variable. These papers assume a priori that the state variables marginals at all times have a density. 
Recently, \cite{PhamWei} managed to drop the density assumption, but
still restricted the admissible controls to be of closed-loop (a.k.a. feedback) type, i.e., deterministic and Lipschitz functions of the current value of the state, which is somewhat restrictive. This feedback form  on the class of controls allows one to reformulate the McKean-Vlasov control problem \reff{V_MKV_App_Intro} as a deterministic control problem in an infinite dimensional space with the marginal distribution as the state variable.  In this paper we will consider the most general case and allow the controls to be open-loop. In this case reformulation mentioned above is no more possible. We will instead work with a proper disintegration of the value function, which we described in \eqref{V_MKV_App_Intro}. The disintegration formula 
\reff{desin} was pointed out heuristically  in \cite{CarmDel14}, see their formulae (40) and (41), but the value function $V$ was not identified. The idea of formulating the McKean-Vlasov control problem as in \eqref{Value_Intro} (rather than as in \eqref{V_MKV_App_Intro}) is inspired by \cite{BLPR14}, where the uncontrolled case is addressed. We will then generalize the randomization approach developed by \cite{KP15} to the McKean-Vlasov control problem corresponding to $V$.


 The DPP that we will prove is the so-called \emph{randomized dynamic programming principle}  (see \cite{BCFP15}), which is the dynamic programming principle for an intensity control problem for a Poisson random measure whose marks leave in
  a subclass of control processes which is dense with respect to the Krylov metric (see Definition 3.2.3 in \cite{Krylov80}).
See \eqref{ValueRand} for the definition of the randomized control problem, Theorem~\ref{Thm:Equivalence} for the equivalence to $V$ (in itself is one of the main technical contributions), and Theorem~\ref{thm:rdpp}, which is our main result, for the statement of the randomized dynamic programming principle. Although, the approach of replacing the original control problem with a randomized version is also taken in  \cite{BCFP15} and \cite{FuhrmanPham15}, our contribution here is in identifying the correct randomization that corresponds to the McKean-Vlasov problem. The
McKean-Vlasov nature of the control problem makes this task rather difficult and as a result the marks of the Poisson random measure live in an abstract space of processes.
We should also emphasize that another relevant issue resolved in this paper concerns the flow properties for the solutions to equations \eqref{StateEq1_Intro} and \eqref{StateEq2_Intro}, see Section \ref{SubS:Flow}. The importance of the flow properties is to prove an identification formula (Lemma \ref{Identification}) between $V$ and the solution to the BSDE, which in turn allows to derive the randomized dynamic programming principle for $V$. 
Our aim is then to use the randomized dynamic programming principle to characterize $V$ through a Hamilton-Jacobi-Bellman equation on the Wasserstein space $\mathscr P_{\text{\tiny$2$}}(\R^n)$, using the recent notion of Lions' differentiability.

Although it is an intermediary step in deriving the randomized DPP, we see Theorem~\ref{Thm:Feynman-Kac} as the second main result of our paper. Here we derive the nonlinear Feynman-Kac representation of the value function $V$ in terms of a class of forward-backward stochastic differential equations with constrained jumps, where the forward process is autonomous.  This representation has been derived in \cite{KP15} for the case of classical stochastic optimal control problems and here we are generalizing it to McKean-Vlasov control problems. The importance of this representation, beyond its intermediary role, 
is  that it would be useful in developing probabilistic numerical schemes for $V$ (see \cite{KLP15} for the case treated in \cite{KP15}). 

The rest of the paper is organized as follows. Section \ref{S:Formulation} is devoted to the formulation of the McKean-Vlasov control problem, and its continuity properties. In Section \ref{S:RandomizedProblem} we introduce the randomized McKean-Vlasov control problem and we prove the fundamental equivalence result between $V$ and $V^\Rc$ (Theorem \ref{Thm:Equivalence}). In Section \ref{S:Feynman-Kac} we prove the nonlinear Feynman-Kac representation for $V$ in terms of the so-called randomized equation, namely BSDE \eqref{BSDE}. In Section \ref{S:DPP} we derive the randomized dynamic programming principle, proving the flow properties (Lemma \ref{L:Flow}) and the identification between $V$ and the solution to the BSDE (Lemma \ref{Identification}). Finally, in the Appendix we prove some convergence results with respect to the 2-Wasserstein metric $\Wc_{\text{\tiny$2$}}$ (Appendix \ref{App:Wass}), we report the proofs of the measurability Lemmata \ref{L:hatP} and \ref{L:P}  (Appendix \ref{App:Proofs}), we state and prove a stability result with respect to the Krylov metric $\tilde\rho$ (Appendix \ref{App:Stability}), we consider an alternative randomization McKean-Vlasov control problem, more similar to the randomized problems studied for instance in \cite{BCFP15,CC,FuhrmanPham15,KP15} (Appendix \ref{App:Randomization}).

\section{Formulation of the McKean-Vlasov control problem}
\label{S:Formulation}

\subsection{Notations}
\label{SubS:Notation}

Consider a complete probability space $(\Omega,\calf,\P)$ and a $d$-dimensional Brownian motion $B=(B_t)_{t\geq0}$ defined on it. Let $\F^B=(\Fc_t^B)_{t\geq 0}$ denote the $\P$-completion of the filtration generated by $B$. Fix a finite time horizon $T>0$ and a Polish space $A$, endowed with a metric $\rho$. We suppose, without loss of generality, that $\rho<1$ (if this is not the case, we replace $\rho$ with the equivalent metric $\rho/(1+\rho)$). We indicate by $\Bc(A)$ the family of Borel subsets of $A$.

Let $\mathscr P_{\text{\tiny$2$}}(\R^n)$ denote the set of all probability measures on $(\R^n,\Bc(\R^n))$ with a finite second-order moment. We endow $\mathscr P_{\text{\tiny$2$}}(\R^n)$ with the 2-Wasserstein metric $\Wc_{\text{\tiny$2$}}$ defined as follows:
\[
\Wc_{\text{\tiny$2$}}(\pi,\pi') = \inf\bigg\{\bigg(\int_{\R^n\times\R^n} |x - x'|^2\,\boldsymbol\pi(dx,dx')\bigg)^{1/2}\colon\boldsymbol\pi\in\mathscr P_{\text{\tiny$2$}}(\R^n\times\R^n)\text{ with marginals $\pi$ and $\pi'$}\bigg\},
\]
for all $\pi,\pi'\in\mathscr P_{\text{\tiny$2$}}(\R^n)$. We recall from Theorem 6.18 in \cite{Villani} that $(\mathscr P_{\text{\tiny$2$}}(\R^n),\Wc_{\text{\tiny$2$}})$ is a complete separable metric space. Notice that
\begin{equation}\label{WassMetricEstimate}
\Wc_{\text{\tiny$2$}}(\P_{\text{\tiny$\xi$}},\P_{\text{\tiny$\xi'$}}) \ \leq \ (\E[|\xi-\xi'|^2])^{1/2}, \qquad \text{for every pair }\xi,\xi'\in L^2(\Omega,\Fc,\P;\R^n),
\end{equation}
where $\P_{\text{\tiny$\xi$}}$ denotes the law under $\P$ of the random variable $\xi\colon\Omega\rightarrow\R^n$. We also denote by $\|\pi\|_{\text{\tiny$2$}}$ the square root of the second-order moment of $\pi\in\mathscr P_{\text{\tiny$2$}}(\R^n)$:
\begin{equation}\label{||pi||_2^2}
\Wc_{\text{\tiny$2$}}(\pi,\delta_0) \ = \ \|\pi\|_{\text{\tiny$2$}} \ = \ \bigg(\int_{\R^n} |x|^2\,\pi(dx)\bigg)^{\frac{1}{2}}, \qquad \text{for all }\pi\in\mathscr P_{\text{\tiny$2$}}(\R^n),
\end{equation}
where $\delta_0$ is the Dirac measure on $\R^n$ concentrated at the origin. We denote $\Bc(\mathscr P_{\text{\tiny$2$}}(\R^n))$ the Borel $\sigma$-algebra on $\mathscr P_{\text{\tiny$2$}}(\R^n)$ induced by the 2-Wasserstein metric $\Wc_{\text{\tiny$2$}}$.

We assume that there exists a sub-$\sigma$-algebra $\Gc\subset\Fc$ such that $B$ is independent of $\Gc$ and $\mathscr P_{\text{\tiny$2$}}(\R^n)=\{\P_{\text{\tiny$\xi$}}\colon\xi\in L^2(\Omega,\Gc,\P;\R^n)\}$.

Finally, we denote $C_{\text{\tiny$2$}}(\R^n)$ the set of real-valued continuous functions with at most quadratic growth, and $\mathscr B_{\text{\tiny$2$}}(\R^n)$ the set of real-valued Borel measurable functions with at most quadratic growth.

\begin{Remark}\label{R:MeasWass}
{\rm
For every $\varphi\in C_{\text{\tiny$2$}}(\R^n)$, let $\Lambda_{_\varphi}\colon\mathscr P_{\text{\tiny$2$}}(\R^n)\rightarrow\R$ be given by
\[
\Lambda_{_\varphi}(\pi) \ = \ \int_{\R^n}\varphi(x)\,\pi(dx), \qquad \text{for every }\pi\in\mathscr P_{\text{\tiny$2$}}(\R^n).
\]
We notice that (as remarked on pages 6-7 in \cite{GangboKimPacini}) $\Bc(\mathscr P_{\text{\tiny$2$}}(\R^n))$ coincides with the $\sigma$-algebra generated by the family of maps $\Lambda_{_\varphi}$, $\varphi\in C_{\text{\tiny$2$}}(\R^n)$. As a consequence, we observe that, given a measurable space $(E,\Ec)$ and a map $F\colon E\rightarrow\mathscr P_{\text{\tiny$2$}}(\R^n)$, then $F$ is measurable if and only if $\Lambda_{_\varphi}\circ F$ is measurable, for every $\varphi\in C_{\text{\tiny$2$}}(\R^n)$. Finally, we notice that if $\varphi\in\mathscr B_{\text{\tiny$2$}}(\R^n)$ then the map $\Lambda_{_\varphi}$ is $\Bc(\mathscr P_{\text{\tiny$2$}}(\R^n))$-measurable. This latter property can be proved using a monotone class argument, noting that $\Lambda_{_\varphi}$ is $\Bc(\mathscr P_{\text{\tiny$2$}}(\R^n))$-measurable whenever $\varphi\in C_{\text{\tiny$2$}}(\R^n)$.
\ep}
\end{Remark}

\subsection{Optimal control of McKean-Vlasov dynamics}
\label{Primal}

Let $\Ac$ denote the set of admissible control processes, which are $\F^B$-progressive processes $\alpha\colon\Omega\times[0,T]\rightarrow A$. Given $(t,x,\xi)\in[0,T]\times\R^n\times L^2(\Omega,\Gc,\P;\R^n)$ and $\alpha\in\Ac$, the controlled state equations are given by:
\begin{align}
dX_s^{t,\xi,\alpha} \ &= \ b\big(s,X_s^{t,\xi,\alpha},\P_{\text{\tiny$X_s^{t,\xi,\alpha}$}},\alpha_s\big)\,ds + \sigma\big(s,X_s^{t,\xi,\alpha},\P_{\text{\tiny$X_s^{t,\xi,\alpha}$}},\alpha_s\big)\,dB_s, \quad &X_t^{t,\xi,\alpha} \ = \ \xi, \label{StateEq1} \\
dX_s^{t,x,\xi,\alpha} \ &= \ b\big(s,X_s^{t,x,\xi,\alpha},\P_{\text{\tiny$X_s^{t,\xi,\alpha}$}},\alpha_s\big)\,ds + \sigma\big(s,X_s^{t,x,\xi,\alpha},\P_{\text{\tiny$X_s^{t,\xi,\alpha}$}},\alpha_s\big)\,dB_s, &X_t^{t,x,\xi,\alpha} \ = \ x, \label{StateEq2}
\end{align}
for all $s\in[t,T]$. The coefficients $b\colon[0,T]\times\R^n\times\mathscr P_{\text{\tiny$2$}}(\R^n)\times A\rightarrow\R^n$ and $\sigma\colon[0,T]\times\R^n\times\mathscr P_{\text{\tiny$2$}}(\R^n)\times A\rightarrow\R^{n\times d}$ are assumed to be Borel measurable. Recall that $\P_{\text{\tiny$X_s^{t,\xi,\alpha}$}}$ denotes the law under $\P$ of the random variable $X_s^{t,\xi,\alpha}\colon\Omega\rightarrow\R^n$. Notice that $(\P_{\text{\tiny$X_s^{t,\xi,\alpha}$}})_{s\in[t,T]}$ depends on $\xi$ only through its law $\pi=\P_{\text{\tiny$\xi$}}$, and $\pi$ is an element of $\mathscr P_{\text{\tiny$2$}}(\R^n)$. As a consequence, $X^{t,x,\xi,\alpha}=(X_s^{t,x,\xi,\alpha})_{s\in[t,T]}$ depends on $\xi$ only through $\pi$. Therefore, we denote $X^{t,x,\xi,\alpha}$ simply by $X^{t,x,\pi,\alpha}$, whenever $\pi=\P_{\text{\tiny$\xi$}}$.
By misuse of notations, we keep the same letter $X$ for the solution to \reff{StateEq1} and \reff{StateEq2}, but we emphasize that 
in \reff{StateEq2}, the coefficients depend on the law of the first component and the SDE for \reff{StateEq2} 
should be viewed as a standard SDE with initial date $(t,x)$ given a control $\alpha$. 

Our aim is to maximize, over all $\alpha\in\Ac$, the following functional
\begin{equation}\label{FunctionalJ}
J(t,x,\pi,\alpha) \ = \ \E\bigg[\int_t^T f\big(s,X_s^{t,x,\pi,\alpha},\P_{\text{\tiny$X_s^{t,\xi,\alpha}$}},\alpha_s\big)\,ds + g\big(X_T^{t,x,\pi,\alpha},\P_{\text{\tiny$X_T^{t,\xi,\alpha}$}}\big)\bigg],
\end{equation}
where $f\colon[0,T]\times\R^n\times\mathscr P_{\text{\tiny$2$}}(\R^n)\times A\rightarrow\R$ and $g\colon\R^n\times\mathscr P_{\text{\tiny$2$}}(\R^n)\rightarrow\R$ are Borel measurable. We impose the following assumptions.

\vspace{3mm}

\noindent {\bf (A1)}
\begin{itemize}
\item [(i)] For every $t$ $\in$ $[0,T]$,  $b(t,\cdot)$, $\sigma(t,\cdot)$ and $f(t,\cdot)$ are continuous on $\R^n\times\mathscr P_{\text{\tiny$2$}}(\R^n)\times A$, and $g$ is continuous on $\R^n\times\mathscr P_{\text{\tiny$2$}}(\R^n)$.
\item [(ii)] For every $(t,x,x',\pi,\pi',a)\in[0,T]\times\R^n\times\R^n\times\mathscr P_{\text{\tiny$2$}}(\R^n)\times\mathscr P_{\text{\tiny$2$}}(\R^n)\times A$,
\begin{align*}
|b(t,x,\pi,a) - b(t,x',\pi',a)| + |\sigma(t,x,\pi,a)-\sigma(t,x',\pi',a)| \ &\leq \ L\big(|x-x'| + \Wc_{\text{\tiny$2$}}(\pi,\pi')\big), \\
|b(t,0,\delta_0,a)| + |\sigma(t,0,\delta_0,a)| \ &\leq \ L, \\
|f(t,x,\pi,a)| + |g(x,\pi)| \ &\leq \ h(\|\pi\|_{\text{\tiny$2$}}) \big(1 + |x|^p\big),
\end{align*}
for some positive constants $L$ and $p$, and some continuous function $h\colon\R_+\rightarrow\R_+$.
\end{itemize}

\vspace{2mm}

Under Assumption {\bf (A1)}, and recalling property \eqref{WassMetricEstimate}, it can be proved by standard arguments that there exists a unique (up to indistinguishability) pair $(X_s^{t,\xi,\alpha},X_s^{t,x,\pi,\alpha})_{s\in[t,T]}$ of continuous $(\Fc_s^B\vee\Gc)_s$-adapted processes solution to equations \eqref{StateEq1}-\eqref{StateEq2}, satisfying
\begin{equation}\label{Estimate}
\sup_{\alpha\in\Ac}\E\Big[\sup_{s\in[t,T]}\big(\big|X_s^{t,\xi,\alpha}\big|^2 + \big|X_s^{t,x,\pi,\alpha}\big|^q\big)\Big] \ < \ \infty,
\end{equation}
for all $q\geq1$. The estimate $\sup_{\alpha\in\Ac}\E[\sup_{s\in[t,T]}|X_s^{t,\xi,\alpha}|^q]<\infty$ holds whenever $|\xi|^q$ is integrable. Notice that $(X_s^{t,x,\pi,\alpha})_{s\in[t,T]}$ is $\F^B$-adapted.

Recalling $\mathscr P_{\text{\tiny$2$}}(\R^n)=\{\P_{\text{\tiny$\xi$}}\colon\xi\in L^2(\Omega,\Gc,\P;\R^n)\}$, we see that $J(t,x,\pi,\alpha)$ is defined for every quadruple $(t,x,\pi,\alpha)\in[0,T]\times\R^n\times\mathscr P_{\text{\tiny$2$}}(\R^n)\times\Ac$. The value function of our stochastic control problem is the function $V$ on $[0,T]\times\R^n\times\mathscr P_{\text{\tiny$2$}}(\R^n)$ defined as
\begin{equation}\label{Value}
V(t,x,\pi) \ = \ \sup_{\alpha\in\Ac} J(t,x,\pi,\alpha),
\end{equation}
for all $(t,x,\pi)\in[0,T]\times\R^n\times\mathscr P_{\text{\tiny$2$}}(\R^n)$.

From estimate \eqref{Estimate}, we see that $\|\P_{X_s^{t,\xi,\alpha}}\|_{\text{\tiny$2$}}\leq M$, for some positive constant $M$ independent of $\alpha\in\Ac$ and $s\in[t,T]$. It follows from the continuity of $h$ that the quantity $h(\|\P_{X_s^{t,\xi,\alpha}}\|_{\text{\tiny$2$}})$ is bounded uniformly with respect to $\alpha$ and $s$. Therefore, by the polynomial growth condition on $f$ and $g$ in Assumption {\bf (A1)}(ii), we deduce that the value function $V$ in \eqref{Value} is always a finite real number on its domain $[0,T]\times\R^n\times\mathscr P_{\text{\tiny$2$}}(\R^n)$, namely $V\colon[0,T]\times\R^n\times\mathscr P_{\text{\tiny$2$}}(\R^n)\rightarrow\R$. In particular, it is easy to see that, under Assumption {\bf (A1)}, $V$ satisfies the following growth condition:
\begin{equation}\label{GrowthV}
|V(t,x,\pi)| \ \leq \ \psi(\|\pi\|_{\text{\tiny$2$}}) \big(1 + |x|^p\big),
\end{equation}
for some continuous function $\psi\colon\R_+\rightarrow\R_+$.

We now study the continuity of $V$. Firstly, we impose the following additional assumption.

\vspace{3mm}

\noindent {\bf (A2)} \qquad For every $t$ $\in$ $[0,T]$ and $R>0$, the map $(x,\pi)\mapsto f(t,\cdot,\cdot,a)$ is uniformly continuous and bounded on $\{(x,\pi)\in\R^n\times\mathscr P_2(\R^n)\colon|x|,\|\pi\|_{\text{\tiny$2$}}\leq R\}$, uniformly with respect to $a\in A$. For every $R>0$, the map $g$ is uniformly continuous and bounded on $\{(x,\pi)\in\R^n\times\mathscr P_2(\R^n)\colon|x|,\|\pi\|_{\text{\tiny$2$}}\leq R\}$.

\begin{Proposition}\label{P:VCont}
Under Assumptions {\bf (A1)} and {\bf (A2)}, for every $t\in[0,T]$ the map $(x,\pi)\mapsto V(t,x,\pi)$ is continuous on $\R^n\times\mathscr P_{\text{\tiny$2$}}(\R^n)$.
\end{Proposition}
\textbf{Proof.}
We begin noting that, as a consequence of Assumption {\bf (A2)}, for every $t\in[0,T]$ and $R>0$, there exists a modulus of continuity $\delta_t^R\colon[0,\infty)\rightarrow[0,\infty)$ such that, for $t\in[0,T)$,
\[
\big|f(t,x,\pi,a) - f(t,x',\pi',a)\big| \ \leq \ \delta_t^R\big(|x - x'| + \Wc_{\text{\tiny$2$}}(\pi,\pi')\big),
\]
and, for $t=T$,
\[
\big|f(T,x,\pi,a) - f(T,x',\pi',a)\big| + \big|g(x,\pi) - g(x',\pi)\big| \ \leq \ \delta_T^R\big(|x - x'| + \Wc_{\text{\tiny$2$}}(\pi,\pi')\big),
\]
for all $(x,\pi),(x',\pi')\in\R^n\times\mathscr P_{\text{\tiny$2$}}(\R^n)$, $a\in A$, with $|x|,|x'|,\|\pi\|_{\text{\tiny$2$}},\|\pi'\|_{\text{\tiny$2$}}\leq R$. Recall that, by definition (see for instance \cite{ArPan56}, page 406), the modulus of continuity $\delta_s^R$ is nondecreasing and $\lim_{\eps\rightarrow0^+}\delta_s^R(\eps)=0$. Moreover, by Assumption {\bf (A2)}, we see that $\delta_s^R$ can be taken bounded. In particular, $\limsup_{\eps\rightarrow+\infty}\delta_s^R(\eps)/\eps=0$. Therefore, without loss of generality, we can suppose that $\delta_s^R$ is concave (see for instance Theorem 1, page 406, in \cite{ArPan56}; we refer, in particular, to the concave modulus of continuity constructed in the proof of Theorem 1 and given by formula (1.6) at page 407). Then, we notice that $\delta_s^R$ is also subadditive. 

Now, fix $t\in[0,T]$ and $(x,\pi),(x_m,\pi_m)\in\R^n\times\mathscr P_{\text{\tiny$2$}}(\R^n)$, with $|x_m-x|\rightarrow0$ and $\Wc_{\text{\tiny$2$}}(\pi_m,\pi)\rightarrow0$ as $m$ goes to infinity. Our aim is to prove that
\begin{equation}\label{ContinuityV}
V(t,x_m,\pi_m) \ \overset{m\rightarrow\infty}{\longrightarrow} \ V(t,x,\pi).
\end{equation}
By Lemma \ref{L:W2Skorohod} we know that there exist random variables $\xi,\xi_m\in L^2(\Omega,\Gc,\P;\R^n)$ such that $\pi=\P_{\text{\tiny$\xi$}}$ and $\pi_m=\P_{\text{\tiny$\xi_m$}}$ under $\P$, moreover $\xi_m$ converges to $\xi$ pointwise $\P$-a.s. and in $L^2(\Omega,\Gc,\P;\R^n)$. In particular, $\sup_m\E[|\xi_m|^2]<\infty$. Then, by standard arguments, we have
\[
\max\bigg\{\sup_{s\in[t,T],\,\alpha\in\Ac}\big\|\P_{X_s^{t,\xi,\alpha}}\big\|_{\text{\tiny$2$}},\sup_m\sup_{s\in[t,T],\,\alpha\in\Ac}\big\|\P_{X_s^{t,\xi_m,\alpha}}\big\|_{\text{\tiny$2$}}\bigg\} \ =: \ \bar R,
\]
for some constant $\bar R\geq 0$. For every $R>\bar R$ and $\alpha\in\Ac$, define the set $E_\alpha\in\Fc$ as
\[
E_\alpha \ := \ \Big\{\omega\in\Omega\colon\sup_{s\in[t,T]}|X_s^{t,x,\pi,\alpha}(\omega)|,\,\sup_m\sup_{s\in[t,T]}|X_s^{t,x_m,\pi_m,\alpha}(\omega)|\leq R\Big\}.
\]
Then, we have
\begin{align}
&|V(t,x,\pi) - V(t,x_m,\pi_m)| \notag \\
&\leq \ \sup_{\alpha\in\Ac} \E\bigg[1_{E_\alpha} \int_t^T \delta_s^R\big(\big|X_s^{t,x,\pi,\alpha} - X_s^{t,x_m,\pi_m,\alpha}\big|\big)\,ds + 1_{E_\alpha}\,\delta_T^R\big(\big|X_T^{t,x,\pi,\alpha} - X_T^{t,x_m,\pi_m,\alpha}\big|\big)\bigg] \notag \\
&\quad \ + \sup_{\alpha\in\Ac} \E\bigg[ 1_{E_\alpha} \int_t^T \delta_s^R\big(\Wc_{\text{\tiny$2$}}\big(\P_{X_s^{t,\xi,\alpha}},\P_{X_s^{t,\xi_m,\alpha}}\big)\big)\,ds + 1_{E_\alpha}\,\delta_T^R\big(\Wc_{\text{\tiny$2$}}\big(\P_{X_T^{t,\xi,\alpha}},\P_{X_T^{t,\xi_m,\alpha}}\big)\big)\bigg] \notag \\
&\quad \ + \sup_{\alpha\in\Ac} \E\bigg[1_{E_\alpha^c}\big|g\big(X_T^{t,x,\pi,\alpha},\P_{X_T^{t,\xi,\alpha}}\big) - g\big(X_T^{t,x_m,\pi_m,\alpha},\P_{X_T^{t,\xi_m,\alpha}}\big)\big| \notag \\
&\quad \ + 1_{E_\alpha^c}\int_t^T \big|f\big(s,X_s^{t,x,\pi,\alpha},\P_{X_s^{t,\xi,\alpha}},\alpha_s\big) - f\big(s,X_s^{t,x_m,\pi_m,\alpha},\P_{X_s^{t,\xi_m,\alpha}},\alpha_s\big)\big|\,ds\bigg] \notag \\
&\leq \ \sup_{\alpha\in\Ac} \E\bigg[\int_t^T \delta_s^R\big(\big|X_s^{t,x,\pi,\alpha} - X_s^{t,x_m,\pi_m,\alpha}\big|\big)\,ds + \delta_T^R\big(\big|X_T^{t,x,\pi,\alpha} - X_T^{t,x_m,\pi_m,\alpha}\big|\big)\bigg] \notag \\
&\quad \ + \sup_{\alpha\in\Ac} \bigg( \int_t^T \delta_s^R\big(\Wc_{\text{\tiny$2$}}\big(\P_{X_s^{t,\xi,\alpha}},\P_{X_s^{t,\xi_m,\alpha}}\big)\big)\,ds + \delta_T^R\big(\Wc_{\text{\tiny$2$}}\big(\P_{X_T^{t,\xi,\alpha}},\P_{X_T^{t,\xi_m,\alpha}}\big)\big)\bigg) \notag \\
&\quad \ + C(1 + |x|^p + |x_m|^p)\,\sup_{\alpha\in\Ac}\,\P(E_\alpha^c), \label{V_x-x'}
\end{align}
for some positive constant $C$, depending only on $\bar R$, $T$, the constants $L$, $p$ in Assumption {\bf (A1)}(ii), and the maximum $\max_{0\leq r\leq \bar R}h(r)$, where the function $h$ was introduced in Assumption {\bf (A1)}(ii). Recalling that $\Wc_{\text{\tiny$2$}}(\P_{X_s^{t,\xi,\alpha}},\P_{X_s^{t,\xi_m,\alpha}})\leq\E[|X_s^{t,\xi,\alpha}-X_s^{t,\xi_m,\alpha}|^2]$ and $\delta_s^R$ is nondecreasing, we find
\begin{equation}\label{delta_W2}
\delta_s^R\big(\Wc_{\text{\tiny$2$}}\big(\P_{X_s^{t,\xi,\alpha}},\P_{X_s^{t,\xi_m,\alpha}}\big)\big) \ \leq \ \delta_s^R\Big(\E\big[\big|X_s^{t,\xi,\alpha}-X_s^{t,\xi_m,\alpha}\big|^2\big]^{1/2}\Big).
\end{equation}
Now, recall the standard estimate
\begin{equation}\label{delta_W2_bis_bis}
\sup_{\alpha\in\Ac}\E\big[\big|X_s^{t,\xi,\alpha}-X_s^{t,\xi_m,\alpha}\big|^2\big]^{1/2} \ \leq \ \hat c\,\E\big[|\xi - \xi_m|^2\big]^{1/2},
\end{equation}
for some positive constant $\hat c$, depending only on $T$ and $L$. Therefore, from \eqref{delta_W2} we obtain
\begin{equation}\label{delta_W2_bis}
\delta_s^R\big(\Wc_{\text{\tiny$2$}}\big(\P_{X_s^{t,\xi,\alpha}},\P_{X_s^{t,\xi_m,\alpha}}\big)\big) \ \leq \ \delta_s^R\Big(\hat c\,\E\big[|\xi - \xi_m|^2\big]^{1/2}\Big).
\end{equation}
On the other hand, from the concavity of $\delta_s^R$, we get
\begin{equation}\label{delta_x-x'}
\E\big[\delta_s^R\big(\big|X_s^{t,x,\pi,\alpha} - X_s^{t,x_m,\pi_m,\alpha}\big|\big)\big] \ \leq \ \delta_s^R\big(\E\big[\big|X_s^{t,x,\pi,\alpha} - X_s^{t,x_m,\pi_m,\alpha}\big|\big]\big).
\end{equation}
By standard arguments, we have
\[
\sup_{\alpha\in\Ac} \E\Big[\sup_{s\in[t,T]}\big|X_s^{t,x,\pi,\alpha} - X_s^{t,x_m,\pi_m,\alpha}\big|\Big] \ \leq \ c\,\Big(|x - x_m| + \sup_{\alpha\in\Ac}\sup_{s\in[t,T]}\Wc_{\text{\tiny$2$}}\big(\P_{X_s^{t,\xi,\alpha}},\P_{X_s^{t,\xi_m,\alpha}}\big)\Big),
\]
where $c$ is a positive constant, depending only on $T$ and $L$. Therefore, by \eqref{delta_W2_bis_bis}, we obtain
\begin{equation}\label{Estimate_x-x'}
\sup_{\alpha\in\Ac} \E\Big[\sup_{s\in[t,T]}\big|X_s^{t,x,\pi,\alpha} - X_s^{t,x_m,\pi_m,\alpha}\big|\Big] \ \leq \ c\,\Big(|x - x_m| + \hat c\,\E\big[|\xi - \xi_m|^2\big]^{1/2}\Big).
\end{equation}
Since $\delta_s^R$ is nondecreasing, from \eqref{delta_x-x'} and \eqref{Estimate_x-x'}, we find
\begin{equation}\label{sup_alpha_delta_x-x'}
\sup_{\alpha\in\Ac}\E\big[\delta_s^R\big(\big|X_s^{t,x,\pi,\alpha} - X_s^{t,x_m,\pi_m,\alpha}\big|\big)\big] \ \leq \ \delta_s^R\Big(c\,|x - x_m| + c\,\hat c\,\E\big[|\xi - \xi_m|^2\big]^{1/2}\Big).
\end{equation}
Concerning $\P(E_\alpha^c)$, we have
\begin{align}
\P(E_\alpha^c) \ &\leq \ \P\Big(\sup_{s\in[t,T]}|X_s^{t,x,\pi,\alpha}|>R\Big) + \P\Big(\sup_{s\in[t,T]}|X_s^{t,x_m,\pi_m,\alpha}|>R\Big) \label{P_E_alpha^c} \\
&\leq \ \frac{1}{R^2}\E\Big[\sup_{s\in[t,T]}|X_s^{t,x,\pi,\alpha}|^2\Big] + \frac{1}{R^2}\E\Big[\sup_{s\in[t,T]}|X_s^{t,x_m,\pi_m,\alpha}|^2\Big] \ \leq \ \frac{c_0}{R^2} \big(1 + |x|^2 + |x_m|^2\big), \notag
\end{align}
for some positive constant $c_0$, depending only on $T$, $L$, $\bar R$. In conclusion, plugging \eqref{delta_W2_bis}-\eqref{sup_alpha_delta_x-x'}-\eqref{P_E_alpha^c} into \eqref{V_x-x'}, we get
\begin{align}\label{V-Vm}
&|V(t,x,\pi) - V(t,x_m,\pi_m)| \notag \\
&\leq \ \int_t^T \delta_s^R\Big(c\,|x - x_m| + c\,\hat c\,\E\big[|\xi - \xi_m|^2\big]^{1/2}\Big)\,ds + \delta_T^R\Big(c\,|x - x_m| + c\,\hat c\,\E\big[|\xi - \xi_m|^2\big]^{1/2}\Big) \notag \\
&\quad \ + \int_t^T \delta_s^R\Big(\hat c\,\E\big[|\xi - \xi_m|^2\big]^{1/2}\Big)\,ds + \delta_T^R\Big(\hat c\,\E\big[|\xi - \xi_m|^2\big]^{1/2}\Big) \notag \\
&\quad \ + \frac{c_0 C}{R^2}\big(1 + |x|^2 + |x_m|^2\big)\big(1 + |x|^p + |x_m|^p\big).
\end{align}
Taking the $\limsup_{m\rightarrow\infty}$ in the above inequality, we find
\[
\limsup_{m\rightarrow\infty} |V(t,x,\pi) - V(t,x_m,\pi_m)| \ \leq \ \frac{c_0 C}{R^2}\big(1 + 2|x|^2\big)\big(1 + 2|x|^p\big).
\]
Letting $R\rightarrow\infty$, we deduce that $\limsup_{n\rightarrow\infty}|V(t,x,\pi)-V(t,x_m,\pi_m)|=0$, therefore \eqref{ContinuityV} holds.
\ep

\vspace{3mm}

We end this section showing that the value function $V\colon[0,T]\times\R^n\times\mathscr P_{\text{\tiny$2$}}(\R^n)\rightarrow\R$ given by \eqref{Value} is the disintegration of the value function $V_{\textup{\tiny MKV}}\colon[0,T]\times L^2(\Omega,\Gc,\P;\R^n)\rightarrow\R$ given by:
\begin{equation}\label{V_MKV_App}
V_{\textup{\tiny MKV}}(t,\xi) \ = \ \sup_{\alpha\in\Ac_{\text{\tiny$\xi$}}} \E\bigg[\int_t^T f\big(s,X_s^{t,\xi,\alpha},\P_{\text{\tiny$X_s^{t,\xi,\alpha}$}}^{\text{\tiny$\xi$}},\alpha_s\big)\,ds + g\big(X_T^{t,\xi,\alpha},\P_{\text{\tiny$X_T^{t,\xi,\alpha}$}}^{\text{\tiny$\xi$}}\big)\bigg],
\end{equation}
for every $(t,\xi)\in[0,T]\times L^2(\Omega,\Gc,\P;\R^n)$, where $\Ac_{\text{\tiny$\xi$}}$ denotes the set of $A$-valued $(\Fc_s^B\vee\sigma(\xi))$-progressive processes, $(X_s^{t,\xi,\alpha})_{s\in[t,T]}$ is the solution to the following equation:
\[
dX_s^{t,\xi,\alpha} \ = \ b\big(s,X_s^{t,\xi,\alpha},\P_{\text{\tiny$X_s^{t,\xi,\alpha}$}}^{\text{\tiny$\xi$}},\alpha_s\big)\,ds + \sigma\big(s,X_s^{t,\xi,\alpha},\P_{\text{\tiny$X_s^{t,\xi,\alpha}$}}^{\text{\tiny$\xi$}},\alpha_s\big)\,dB_s, \quad X_t^{t,\xi,\alpha} \ = \ \xi,
\]
for all $s\in[t,T]$, with $\alpha\in\Ac_{\text{\tiny$\xi$}}$, and $\P_{\text{\tiny$X_s^{t,\xi,\alpha}$}}^{\text{\tiny$\xi$}}$ denotes the regular conditional distribution of the random variable $X_s^{t,\xi,\alpha}\colon\Omega\rightarrow\R^n$ with respect to 
$\sigma(\xi)$, whose existence is guaranteed for instance by Theorem 6.3 in \cite{Kallenberg}.
\begin{Proposition}\label{P:V_MKV=V}
Under Assumptions {\bf (A1)} and {\bf (A2)}, for every $(t,\xi)\in[0,T]\times L^2(\Omega,\Gc,\P;\R^n)$, with $\pi=\P_{\text{\tiny$\xi$}}$ under $\P$, we have
\[
V_{\textup{\tiny MKV}}(t,\xi) \ = \ \E\big[V(t,\xi,\pi)\big],
\]
or, equivalently,
\[
V_{\textup{\tiny MKV}}(t,\xi) \ = \ \int_{\R^n} V(t,x,\pi)\,\pi(dx).
\]
\end{Proposition}
\textbf{Proof.}
Fix $t\in[0,T]$. Recall from Proposition \ref{P:VCont} that the map $(x,\pi)\mapsto V(t,x,\pi)$ is continuous on $\R^n\times\mathscr P_{\text{\tiny$2$}}(\R^n)$. Proceeding as in the proof of Proposition \ref{P:VCont}, we can also prove that the map $\xi\mapsto V_{\textup{\tiny MKV}}(t,\xi)$ is continuous on $L^2(\Omega,\Gc,\P;\R^n)$. As a consequence, it is enough to prove the Proposition for $\xi\in L^2(\Omega,\Gc,\P;\R^n)$ taking only a finite number of values, the general result being proved by approximation. In other words, we suppose that
\[
\xi \ = \ \sum_{k=0}^K x_k\,1_{E_k},
\]
for some $K\in\N$, $x_k\in\R^n$, $E_k\in\sigma(\xi)$, with $(E_k)_{k=1,\ldots,K}$ being a partition of $\Omega$. Notice that $\alpha\in\Ac{\text{\tiny$\xi$}}$ if and only if
\begin{equation}\label{alpha=alpha_kE_k}
\alpha \ = \ \sum_{k=0}^K \alpha_k\,1_{E_k},
\end{equation}
for some $\alpha_k\in\Ac$. We also observe that
\[
X_s^{t,\xi,\alpha} \ = \ \sum_{k=0}^K X_s^{t,x_k,\alpha_k}\,1_{E_k}, \qquad\qquad \P_{\text{\tiny$X_s^{t,\xi,\alpha}$}}^{\text{\tiny$\xi$}} \ = \ \sum_{k=0}^K \P_{\text{\tiny$X_s^{t,x_k,\alpha_k}$}}\,1_{E_k}.
\]
Then, the stochastic processes $(X_s^{t,\xi,\alpha})_{s\in[t,T]}$ and $(X_s^{t,x_1,\delta_{x_1},\alpha_1}\,1_{E_1}+\cdots+X_s^{t,x_K,\delta_{x_K},\alpha_K}\,1_{E_K})_{s\in[t,T]}$ are indistinguishable, since they solve the same equation. Therefore
\begin{align}\label{Proof_alpha=alpha_kE_k}
V_{\textup{\tiny MKV}}(t,\xi) &= \sup_{\alpha\in\Ac_{\text{\tiny$\xi$}}} \E\bigg[\int_t^T f\big(s,X_s^{t,\xi,\alpha},\P_{\text{\tiny$X_s^{t,\xi,\alpha}$}}^{\text{\tiny$\xi$}},\alpha_s\big)\,ds + g\big(X_T^{t,\xi,\alpha},\P_{\text{\tiny$X_T^{t,\xi,\alpha}$}}^{\text{\tiny$\xi$}}\big)\bigg] \\
&= \sup_{\alpha\in\Ac_{\text{\tiny$\xi$}}} \E\bigg[\sum_{k=0}^K\bigg(\int_t^T f\big(s,X_s^{t,x_k,\delta_{x_k},\alpha_k},\P_{\text{\tiny$X_s^{t,x_k,\alpha_k}$}},(\alpha_k)_s\big)\,ds + g\big(X_T^{t,x_k,\delta_{x_k},\alpha_k},\P_{\text{\tiny$X_T^{t,x_k,\alpha_k}$}}\big)\bigg)\,1_{E_k}\bigg]. \notag
\end{align}
Since $\xi$ is independent of  $X^{t,x_k,\delta_{x_K},\alpha_k}$ and of $\alpha_k$, we can write the last quantity in \eqref{Proof_alpha=alpha_kE_k} as
\[
V_{\textup{\tiny MKV}}(t,\xi) = \sup_{\alpha\in\Ac_{\text{\tiny$\xi$}}} \E\bigg[\sum_{k=0}^K\E\bigg[\int_t^T f\big(s,X_s^{t,x_k,\delta_{x_k},\alpha_k},\P_{\text{\tiny$X_s^{t,x_k,\alpha_k}$}},(\alpha_k)_s\big)\,ds + g\big(X_T^{t,x_k,\delta_{x_k},\alpha_k},\P_{\text{\tiny$X_T^{t,x_k,\alpha_k}$}}\big)\bigg]\,1_{E_k}\bigg].
\]
From \eqref{alpha=alpha_kE_k}, we conclude that
\begin{align*}
V_{\textup{\tiny MKV}}(t,\xi) &= \E\bigg[\sum_{k=0}^K \sup_{\alpha_k\in\Ac} \E\bigg[\int_t^T f\big(s,X_s^{t,x_k,\delta_{x_k},\alpha_k},\P_{\text{\tiny$X_s^{t,x_k,\alpha_k}$}},(\alpha_k)_s\big)\,ds + g\big(X_T^{t,x_k,\delta_{x_k},\alpha_k},\P_{\text{\tiny$X_T^{t,x_k,\alpha_k}$}}\big)\bigg]\,1_{E_k}\bigg] \\
&= \E\bigg[\sum_{k=0}^K V(t,x_k,\delta_{x_k})\,1_{E_k}\bigg] = \E\big[V(t,\xi,\pi)\big].
\end{align*}
\ep

\section{The randomized McKean-Vlasov control problem}
\label{S:RandomizedProblem}

Following Definition 3.2.3 in \cite{Krylov80}, we define on $\Ac$ the metric $\tilde\rho$ given by:
\begin{equation}\label{KrylovMetric}
\tilde\rho(\alpha,\beta) \ := \ \E\bigg[\int_0^T \rho(\alpha_t,\beta_t)\,dt\bigg],
\end{equation}
where we recall that $\rho$ is a metric on $A$ satisfying $\rho<1$. Notice that convergence with respect to $\tilde\rho$ is equivalent to convergence in $d\P\,dt$-measure. We also observe that $(\Ac,\tilde\rho)$ is a metric space (identifying processes $\alpha$ and $\beta$ which are equal $d\P\,dt$-a.e. on $\Omega\times[0,T]$). Moreover, since $A$ is a Polish space, it turns out that $(\Ac,\tilde\rho)$ is also a Polish space (separability follows from Lemma 3.2.6 in \cite{Krylov80}, completeness follows from the completeness of $A$ and the fact that a $\tilde\rho$-limit of $\F^B$-progressive processes is still $\F^B$-progressive). We denote by $\Bc(\Ac)$ the family of Borel subsets of $\Ac$.

Following \cite{Krylov80}, we introduce the following subset of admissible control processes.

\begin{Definition}
For every $t\in[0,T]$, let $(E_\ell^t)_{\ell\geq1}\in\Fc$ be a countable class of subsets of $\Omega$ which generates $\sigma(B_s,\,s\in[0,t])$. Fix a countable dense subset $(a_m)_{m\geq1}$ of $A$. Fix also, for every integer $k\geq1$, a subdivision $I_k:=\{0=:t_0<t_1<\ldots<t_k:=T\}$ of the interval $[0,T]$, with the diameter  $\max_{i=1,\ldots,k}(t_i-t_{i-1})$ of the subdivision $I_k$ going to zero as $k\rightarrow\infty$. Then, we denote
\begin{align*}
\Ac_{\textup{\tiny{step}}} \ := \ \Big\{\alpha\in\Ac\colon &\text{there exist $k\geq1$, $M\geq1$, $L\geq1$, such that, for every $i=0,\ldots,k-1$,} \\
&\alpha_{t_i}\colon\Omega\rightarrow(a_m)_{m=1,\ldots,M}\text{, with }\alpha_{t_i}\text{ constant on the sets of the partition} \\
&\text{generated by }E_1^{t_i},\ldots,E_L^{t_i}\text{, and, for every }t\in[0,T], \\
&\alpha_t = \alpha_{t_0}\,1_{[t_0,t_1)}(t) + \cdots + \alpha_{t_{k-1}}\,1_{[t_{k-1},t_k)}(t) + \alpha_{t_k}\,1_{\{t_k\}}(t)\Big\}.
\end{align*}
\end{Definition}

\begin{Remark}\label{R:Astep}
{\rm
Notice that $\Ac_{\textup{\tiny{step}}}$ depends (even if we omit to write explicitly this dependence) on the two sequences $(a_m)_{m\geq1}$ and $(I_k)_{k\geq1}$, which are supposed to be fixed throughout the paper. The set $\Ac_{\textup{\tiny{step}}}$, with $\alpha_{t_i}$ being $\sigma(B_s,\,s\in[0,t_i])$-measurable, is introduced in the proof of Lemma 3.2.6 in \cite{Krylov80}, where it is proved that it is dense in $\Ac$ with respect to the metric $\tilde\rho$ defined in \eqref{KrylovMetric}. It can be shown (proceeding as in the proof of Lemma \ref{L:Stability}) that the map $\alpha\mapsto J(t,x,\pi,\alpha)$ is continuous with respect to $\tilde\rho$, so that we could define $V(t,x,\pi)$ in the following equivalent way:
\begin{equation}\label{V_step}
V(t,x,\pi) \ = \ \sup_{\alpha\in\Ac_{\textup{\tiny{step}}}}J(t,x,\pi,\alpha).
\end{equation}
Finally, we observe that $\Ac_{\textup{\tiny{step}}}$ is a countable set, so that it is a Borel subset of $\Ac$, namely $\Ac_{\textup{\tiny{step}}}\in\Bc(\Ac)$.
\ep}
\end{Remark}

Now, in order to implement the randomization method, it is better to reformulate the original McKean-Vlasov control problem as follows. Let ${\boldsymbol\Ac}_{\textup{\tiny{step}}}$ be the following set:
\[
\boldsymbol\Ac_{\textup{\tiny{step}}} \ := \ \big\{\boldsymbol\alpha\colon[0,T]\rightarrow\Ac_{\textup{\tiny{step}}}\colon\,\text{$\boldsymbol\alpha$ is Borel-measurable, c\`adl\`ag, and piecewise constant}\big\}.
\]
It is easy to see that, for every $\boldsymbol\alpha\in\boldsymbol\Ac_{\textup{\tiny{step}}}$, the stochastic process $((\boldsymbol\alpha_s)_s)_{s\in[0,T]}$ is an element of $\Ac$. Vice versa, for every element $\hat\alpha\in\Ac_{\textup{\tiny{step}}}$, there exists $\hat{\boldsymbol\alpha}\in\boldsymbol\Ac_{\textup{\tiny{step}}}$ such that $((\hat{\boldsymbol\alpha}_s)_s)_{s\in[0,T]}$ coincides with $\hat\alpha$ (take $\hat{\boldsymbol\alpha}_s=\hat\alpha$, for every $s\in[0,T]$). Hence, by \eqref{V_step},
\[
V(t,x,\pi) \ = \ \sup_{\boldsymbol\alpha\in\boldsymbol\Ac_{\textup{\tiny{step}}}} J\big(t,x,\pi,((\boldsymbol\alpha_s)_s)_{s\in[0,T]}\big).
\]
On the right-hand side of the above identity we have an optimization problem with class of admissible control processes given by $\{((\boldsymbol\alpha_s)_s)_{s\in[0,T]}\colon\boldsymbol\alpha$ $\in\boldsymbol\Ac_{\textup{\tiny{step}}}\}$. We now randomize this latter control problem.

Consider another complete probability space $(\Omega^1,\Fc^1,\P^1)$. We denote by $\E^1$ the $\P^1$-expected value. We suppose that a Poisson random measure $\mu$ on $\R_+\times\Ac$ is defined on $(\Omega^1,\Fc^1,\P^1)$. The random measure $\mu$ has compensator $\lambda(d\alpha)\,dt$, for some finite positive measure $\lambda$ on $\Ac$, with full topological support given by $\Ac_{\textup{\tiny{step}}}$. We denote $\tilde\mu(dt\,d\alpha):=\mu(dt\,d\alpha)-\lambda(d\alpha)\,dt$ the compensated martingale measure associated to $\mu$. We introduce $\F^\mu=(\Fc_t^\mu)_{t\geq0}$, which is the $\P^1$-completion of the filtration generated by $\mu$, given by:
\[
\Fc_t^\mu \ = \ \sigma\big(\mu((0,s]\times\Ac')\colon s\in[0,t],\,\Ac'\subset\Ac_{\textup{\tiny{step}}}\big)\vee\Nc^1,
\]
for all $t\geq0$, where $\Nc^1$ is the class of $\P^1$-null sets of $\Fc^1$. We also denote $\Pc(\F^\mu)$ the predictable $\sigma$-algebra on $\Omega^1\times\R_+$ corresponding to $\F^\mu$.

We recall that $\mu$ is associated to a marked point process $(T_n,\Ac_n)_{n\geq1}$ on $\R_+\times \Ac$ by the formula $\mu=\sum_{n\ge 1}\delta_{(T_n,\Ac_n)}$, where $\delta_{(T_n,\Ac_n)}$ is the Dirac measure concentrated at the random point $(T_n,\Ac_n)$. We recall that every $T_n$ is an $\F^\mu$-stopping time and every $\Ac_n$ is $\Fc_{T_n}^\mu$-measurable.

Let $\bar\Omega=\Omega\times\Omega^1$, and let $\bar\Fc$ be the $\P\otimes\P^1$-completion of $\Fc\otimes\Fc^1$, and $\bar\P$ the extension of $\P\otimes\P^1$ to $\bar\Fc$. We denote by $\bar\Gc$, $\bar B$, $\bar\mu$ the canonical extensions of $\Gc$, $B$, $\mu$, to $\bar\Omega$, given by: $\bar\Gc:=\{G\times\Omega^1\colon G\in\Gc\}$, $\bar B(\omega,\omega^1):=B(\omega)$, $\bar\mu(\omega,\omega^1;dt\,d\alpha):=\mu(\omega^1;dt\,d\alpha)$. Let $\bar\F^B=(\bar\Fc_t^B)_{t\geq 0}$ (resp. $\bar\F^\mu=(\bar\Fc_t^\mu)_{t\geq0}$) denote the $\bar\P$-completion of the filtration generated by $\bar B$ (resp. $\bar\mu$). Notice that $\bar\Fc_\infty^B$ and $\bar\Fc_\infty^\mu$ are independent.

Let $\bar\F^{B,\mu}=(\bar\Fc_t^{B,\mu})_{t\geq0}$ denote the $\bar\P$-completion of the filtration generated by $\bar B$ and $\bar\mu$. Notice that $\bar B$ is a Brownian motion with respect to $\bar\F^{B,\mu}$ and the $\bar\F^{B,\mu}$-compensator of $\bar\mu$ is given by $\lambda(d\alpha)\,dt$. We define the $A$-valued piecewise constant process $\bar I=(\bar I_t)_{t\geq0}$ on $(\bar\Omega,\bar\Fc,\bar\P)$ as follows:
\begin{equation}\label{barI}
\bar I_t(\omega,\omega^1) \ = \ \sum_{n\ge 0} (\Ac_n(\omega^1))_{t\wedge T}(\omega)\,1_{[T_n(\omega^1),T_{n+1}(\omega^1))}(t), \qquad \text{for all }t\geq0,
\end{equation}
where $T_0:=0$ and $\Ac_0:=\bar\alpha$, for some deterministic and arbitrary control process $\bar\alpha\in\Ac_{\textup{\tiny{step}}}$, which will remain fixed throughout the paper. Notice that $\bar I$ is $\bar\F^{B,\mu}$-adapted.

Randomizing the control in \eqref{StateEq1}-\eqref{StateEq2}, we are led to consider the following equations on $(\bar\Omega,\bar\Fc,\bar\P)$, for every $(t,x,\bar\xi)\in[0,T]\times\R^n\times L^2(\bar\Omega,\bar\Gc,\bar\P;\R^n)$, with $\pi=\P_{\text{\tiny$\bar\xi$}}$ under $\bar\P$:
\begin{align}
d\bar X_s^{t,\bar\xi} \ &= \ b\big(s,\bar X_s^{t,\bar\xi},\P_{\text{\tiny$\bar X_s^{t,\bar\xi}$}}^{\text{\tiny$\bar\Fc_s^\mu$}},\bar I_s\big)\,ds + \sigma\big(s,\bar X_s^{t,\bar\xi},\P_{\text{\tiny$\bar X_s^{t,\bar\xi}$}}^{\text{\tiny$\bar\Fc_s^\mu$}},\bar I_s\big)\,d\bar B_s, \qquad & \bar X_t^{t,\bar\xi} \ = \ \bar\xi, \label{StateEq1Rand} \\
d\bar X_s^{t,x,\pi} \ &= \ b\big(s,\bar X_s^{t,x,\pi},\P_{\text{\tiny$\bar X_s^{t,\bar\xi}$}}^{\text{\tiny$\bar\Fc_s^\mu$}},\bar I_s\big)\,ds + \sigma\big(s,\bar X_s^{t,x,\pi},\P_{\text{\tiny$\bar X_s^{t,\bar\xi}$}}^{\text{\tiny$\bar\Fc_s^\mu$}},\bar I_s\big)\,d\bar B_s, & \bar X_t^{t,x,\pi} \ = \ x, \label{StateEq2Rand}
\end{align}
for all $s\in[t,T]$, where $\P_{\text{\tiny$\bar X_s^{t,\bar\xi}$}}^{\text{\tiny$\bar\Fc_s^\mu$}}$ denotes the regular conditional distribution of the random variable $\bar X_s^{t,\bar\xi}\colon\bar\Omega\rightarrow\R^n$ with respect to $\bar\Fc_s^\mu$, whose existence is guaranteed for instance by Theorem 6.3 in \cite{Kallenberg}. Notice that $\P_{\text{\tiny$\bar X_s^{t,\bar\xi}$}}^{\text{\tiny$\bar\Fc_s^\mu$}}$ depends on $\xi$ only through its law $\pi$, so that equation \eqref{StateEq2Rand} depends only on $\pi$. Under Assumption {\bf (A1)}, it follows by standard arguments that there exists a unique (up to indistinguishability) pair $(\bar X_s^{t,\xi},\bar X_s^{t,x,\pi})_{s\in[t,T]}$ of continuous $(\bar\Fc_s^{B,\mu}\vee\bar\Gc)_s$-adapted processes solution to equations \eqref{StateEq1Rand}-\eqref{StateEq2Rand}, satisfying
\begin{equation}\label{EstimateRand}
\bar\E\Big[\sup_{s\in[t,T]}\big(\big|\bar X_s^{t,\bar\xi}\big|^2 + \big|\bar X_s^{t,x,\pi}\big|^q\big)\Big] \ < \ \infty,
\end{equation}
for all $q\geq1$, where $\bar\E$ denotes the $\bar\P$-expected value. Moreover, $(\bar X_s^{t,x,\pi})_{s\in[t,T]}$ is $\bar\F^{B,\mu}$-adapted. 

We now prove two technical results concerning the process $(\P_{\text{\tiny$\bar X_s^{t,\bar\xi}$}}^{\text{\tiny$\bar\Fc_s^\mu$}})_{s\in[t,T]}$. In particular, the first result (Lemma \ref{L:hatP}) concerns a particular version of $(\P_{\text{\tiny$\bar X_s^{t,\bar\xi}$}}^{\text{\tiny$\bar\Fc_s^\mu$}})_{s\in[t,T]}$, which will be used in the proof of Lemma \ref{L:P}. This latter proves the existence of another version of $(\P_{\text{\tiny$\bar X_s^{t,\bar\xi}$}}^{\text{\tiny$\bar\Fc_s^\mu$}})_{s\in[t,T]}$, which will be used throughout the paper.

\begin{Lemma}\label{L:hatP}
Under Assumption {\bf (A1)}, for every $(t,\pi)\in[0,T]\times\mathscr P_{\text{\tiny$2$}}(\R^n)$, there exists a $\mathscr P_{\text{\tiny$2$}}(\R^n)$-valued $\F^\mu$-predictable stochastic process $(\hat\P_s^{t,\pi})_{s\in[t,T]}$ which is a version of $(\P_{\text{\tiny$\bar X_s^{t,\bar\xi}$}}^{\text{\tiny$\bar\Fc_s^\mu$}})_{s\in[t,T]}$, with $\bar\xi\in L^2(\bar\Omega,\bar\Gc,\bar\P;\R^n)$ such that $\pi=\P_{\text{\tiny$\bar\xi$}}$ under $\bar\P$. For all $s\in[t,T]$, $\hat\P_s^{t,\pi}$ is given by
\begin{equation}\label{RCPD}
\hat\P_s^{t,\pi}(\omega^1)[\varphi] \ = \ \E\big[\varphi\big(\bar X_s^{t,\bar\xi}(\cdot,\omega^1)\big)\big],
\end{equation}
for every $\omega^1\in\Omega^1$ and $\varphi\in\mathscr B_{\text{\tiny$2$}}(\R^n)$.
\end{Lemma}
\textbf{Proof.}
See Appendix \ref{App:Proofs}.
\ep

\begin{Lemma}\label{L:P}
Under Assumption {\bf (A1)}, for every $t\in[0,T]$, there exists a measurable map $\P_\cdot^{t,\cdot}\colon(\Omega^1\times[t,T]\times\mathscr P_{\text{\tiny$2$}}(\R^n),$ $\Fc^1\otimes\Bc([t,T])\otimes\Bc(\mathscr P_{\text{\tiny$2$}}(\R^n)))\rightarrow(\mathscr P_{\text{\tiny$2$}}(\R^n),\Bc(\mathscr P_{\text{\tiny$2$}}(\R^n)))$ such that
\[
\P_s^{t,\pi} \ = \ \P_{\text{\tiny$\bar X_s^{t,\bar\xi}$}}^{\text{\tiny$\bar\Fc_s^\mu$}},
\]
$\P^1$-a.s., for every $s\in[t,T]$, $\pi\in\mathscr P_{\text{\tiny$2$}}(\R^n)$, where $\bar\xi\in L^2(\bar\Omega,\bar\Gc,\bar\P;\R^n)$ has law $\pi$ under $\bar\P$. In other words, for every $s\in[t,T]$ and $\pi\in\mathscr P_{\text{\tiny$2$}}(\R^n)$, $(\P_s^{t,\pi})_{s\in[t,T]}$ is a version of $(\P_{\text{\tiny$\bar X_s^{t,\bar\xi}$}}^{\text{\tiny$\bar\Fc_s^\mu$}})_{s\in[t,T]}$.
\end{Lemma}
\textbf{Proof.}
See Appendix \ref{App:Proofs}.
\ep

\vspace{3mm}

From now on, we will always suppose that $(\P_{\text{\tiny$\bar X_s^{t,\bar\xi}$}}^{\text{\tiny$\bar\Fc_s^\mu$}})_{s\in[t,T]}$ stands for the stochastic process $(\P_s^{t,\pi})_{s\in[t,T]}$ introduced in Lemma \ref{L:P}.

Let us now formulate the randomized McKean-Vlasov control problem. An admissible control is a $\Pc(\F^\mu)\otimes\Bc(\Ac)$-measurable map $\nu\colon\Omega^1\times\R_+\times\Ac\rightarrow(0,\infty)$, which is both bounded away from zero and bounded from above: $0<\inf_{\Omega^1\times\R_+\times\Ac}\nu\leq\sup_{\Omega^1\times\R_+\times\Ac}\nu<\infty$. We denote by $\Vc$ the set of admissible controls. Given $\nu\in\Vc$, we define $\P^\nu$ on $(\Omega^1,\Fc^1)$ as $d\P^\nu=\kappa_T^\nu\,d\P^1$, where $\kappa^\nu=(\kappa_t^\nu)_{t\in[0,T]}$ is the Dol\'eans exponential process on $(\Omega^1,\Fc^1,\P^1)$ defined as
\begin{align*}
\kappa_t^\nu \ &= \ \Ec_t\bigg(\int_0^\cdot\int_\Ac (\nu_s(\alpha) - 1)\,\tilde\mu(ds\,d\alpha)\bigg) \\
&= \ \exp\bigg(\int_0^t\int_\Ac \ln\nu_s(\alpha)\,\mu(ds\,d\alpha) - \int_0^t\int_\Ac(\nu_s(\alpha)-1)\,\lambda(d\alpha)\,ds\bigg), \qquad \text{for all }t\in[0,T]. \notag
\end{align*}
Notice that $\kappa^\nu$ is an $\F^\mu$-martingale under $\P^1$, so that $\P^\nu$ is a probability measure on $(\Omega^1,\Fc^1)$. We denote by $\E^\nu$ the $\P^\nu$-expected value. Observe that, by the Girsanov theorem, the $\F^\mu$-compensator of $\mu$ under $\P^\nu$ is given by $\nu_t(\alpha)\,\lambda(d\alpha)\,dt$. Let $\bar\P^\nu$ denote the extension of $\P\otimes\P^\nu$ to $(\bar\Omega,\bar\Fc)$. Then $d\bar\P^\nu=\bar\kappa_T^\nu d\bar\P$, where $\bar\kappa_t^\nu(\omega,\omega^1):=\kappa_t^\nu(\omega^1)$, for all $t\in[0,T]$. Using again the Girsanov theorem, we see that the $\bar\F^{B,\mu}$-compensator of $\bar\mu$ under $\bar\P^\nu$ is $\bar\nu_t(\alpha)\,\lambda(d\alpha)\,dt$, where $\bar\nu_t(\omega,\omega^1,\alpha):=\nu_t(\omega^1,\alpha)$ is the canonical extension of $\nu$ to $\bar\Omega\times\R_+\times A$.

Notice that a $\bar\Gc$-measurable $\bar\xi\colon\bar\Omega\rightarrow\R^n$ has law $\pi$ under $\bar\P$ if and only if it has the same law under $\bar\P^\nu$. In particular, $\bar\xi\in L^2(\bar\Omega,\bar\Gc,\bar\P;\R^n)$ if and only if $\bar\xi\in L^2(\bar\Omega,\bar\Gc,\bar\P^\nu;\R^n)$. As a consequence, the following generalization of estimate \eqref{EstimateRand} holds ($\bar\E^\nu$ denotes the $\bar\P^\nu$-expected value):
\[
\sup_{\nu\in\Vc}\bar\E^\nu\Big[\sup_{s\in[t,T]}\big(\big|\bar X_s^{t,\bar\xi}\big|^2 + \big|\bar X_s^{t,x,\pi}\big|^q\big)\Big] \ < \ \infty,
\]
for all $q\geq1$, for every $(t,x,\bar\xi)\in[0,T]\times\R^n\times L^2(\bar\Omega,\bar\Gc,\bar\P;\R^n)$, with $\pi=\P_{\text{\tiny$\bar\xi$}}$ under $\bar\P$ (or, equivalently, under $\bar\P^\nu$).

Let $(t,x,\bar\xi)\in[0,T]\times\R^n\times L^2(\bar\Omega,\bar\Gc,\bar\P;\R^n)$, with $\pi=\P_{\text{\tiny$\xi$}}$ under $\bar\P$, and $\nu\in\Vc$, then the gain functional for the randomized McKean-Vlasov control problem is given by:
\[
J^\Rc(t,x,\pi,\nu) \ = \ \bar\E^\nu\bigg[\int_t^T f\big(s,\bar X_s^{t,x,\pi},\P_{\text{\tiny$\bar X_s^{t,\bar\xi}$}}^{\text{\tiny$\bar\Fc_s^\mu$}},\bar I_s\big)\,ds + g\big(\bar X_T^{t,x,\pi},\P_{\text{\tiny$\bar X_T^{t,\bar\xi}$}}^{\text{\tiny$\bar\Fc_T^\mu$}}\big)\bigg].
\]
As for the functional \eqref{FunctionalJ}, the quantity $J^\Rc(t,x,\pi,\nu)$ is defined for every $(t,x,\pi,\nu)\in[0,T]\times\R^n\times\mathscr P_{\text{\tiny$2$}}(\R^n)\times\Vc$, since by assumption $\mathscr P_{\text{\tiny$2$}}(\R^n)=\{\P_{\text{\tiny$\xi$}}\colon\xi\in L^2(\bar\Omega,\bar\Gc,\bar\P;\R^n)\}$. Then, we can define the value function of the randomized McKean-Vlasov control problem as
\begin{equation}\label{ValueRand}
V^\Rc(t,x,\pi) \ = \ \sup_{\nu\in\Vc} J^\Rc(t,x,\pi,\nu),
\end{equation}
for all $(t,x,\pi)\in[0,T]\times\R^n\times\mathscr P_{\text{\tiny$2$}}(\R^n)$.

\begin{Remark}\label{R:hatV}
{\rm
Let $\hat\Vc$ be the set of $\Pc(\F^{\mu})\otimes\Bc(\Ac)$-measurable maps $\hat\nu\colon\Omega^1\times\R_+\times\Ac\rightarrow(0,\infty)$, which are bounded from above $\sup_{\Omega^1\times\R_+\times\Ac}\hat\nu<\infty$, but not necessarily bounded away from zero. For every $(t,x,\pi)\in[0,T]\times\R^n\times\mathscr P_{\text{\tiny$2$}}(\R^n)$, we define
\[
\hat V^\Rc(t,x,\pi) \ = \ \sup_{\hat\nu\in\hat\Vc} J^\Rc(t,x,\pi,\hat\nu)
\]
In \cite{BCFP15} the randomized control problem is formulated over $\hat\Vc$. Here we considered $\Vc$ because this set is more convenient for the proof of Theorem \ref{Thm:Equivalence}. However, notice that
\begin{equation}\label{equivalence_hat}
V^\Rc(t,x,\pi) \ = \ \hat V^\Rc(t,x,\pi).
\end{equation}
Indeed, clearly we have $\Vc\subset\hat\Vc$, so that $V^\Rc(t,x,\pi)\leq\hat V^\Rc(t,x,\pi)$. On the other hand, let $\hat\nu\in\hat\Vc$ and define $\nu^\eps=\hat\nu\vee\eps$, for every $\eps\in(0,1)$. Observe that $\nu^\eps\in\Vc$ and $\bar\kappa_T^{\nu^\eps}$ converges pointwise $\bar\P$-a.s. to $\bar\kappa_T^{\hat\nu}$. Then, it is easy to see that
\[
J^\Rc(t,x,\pi,\nu^\eps) \ = \ \bar\E\bigg[\bar\kappa_T^{\nu^\eps}\bigg(\int_t^T f\big(s,\bar X_s^{t,x,\pi},\P_{\text{\tiny$\bar X_s^{t,\bar\xi}$}}^{\text{\tiny$\bar\Fc_s^\mu$}},\bar I_s\big)\,ds + g\big(\bar X_T^{t,x,\pi},\P_{\text{\tiny$\bar X_T^{t,\bar\xi}$}}^{\text{\tiny$\bar\Fc_T^\mu$}}\big)\bigg)\bigg] \ \overset{\eps\rightarrow0^+}{\longrightarrow} \ J^\Rc(t,x,\pi,\hat\nu).
\]
This implies that $J^\Rc(t,x,\pi,\hat\nu)\leq\sup_{\nu\in\Vc}J^\Rc(t,x,\pi,\nu)$, from which we get the other inequality $\hat V^\Rc(t,x,\pi)\leq V^\Rc(t,x,\pi)$, and identity \eqref{equivalence_hat} follows.
\ep
}
\end{Remark}

We can now prove one of the main results of the paper, namely the equivalence of the two value functions $V$ and $V^\Rc$.

\begin{Theorem}\label{Thm:Equivalence}
Under Assumption {\bf (A1)}, the value function $V$ in \eqref{Value} of the McKean-Vlasov control problem coincides with the value function $V^\Rc$ in \eqref{ValueRand} of the randomized problem:
\[
V(t,x,\pi) \ = \ V^\Rc(t,x,\pi),
\]
for all $(t,x,\pi)\in[0,T]\times\R^n\times\mathscr P_{\text{\tiny$2$}}(\R^n)$.
\end{Theorem}

\begin{Remark}
{\rm
As an immediate consequence of Theorem \ref{Thm:Equivalence}, we see that $V^\Rc$ does not depend on $a_0$ and $\lambda$, since $V$ does not depend on them.
\ep
}
\end{Remark}

\noindent\textbf{Proof (of Theorem \ref{Thm:Equivalence}).}
Fix $(t,x,\xi)\in[0,T]\times\R^n\times L^2(\Omega,\Gc,\P;\R^n)$, with $\pi=\P_{\text{\tiny$\xi$}}$ under $\P$. Set $\bar\xi(\omega,\omega^1):=\xi(\omega)$, then $\bar\xi\in L^2(\bar\Omega,\bar\Gc,\bar\P;\R^n)$ and $\pi=\P_{\text{\tiny$\bar\xi$}}$ under $\bar\P$. We split the proof of the equality $V(t,x,\pi)=V^\Rc(t,x,\pi)$ into three steps, that we now summarize:
\begin{itemize}
\item[I)] In step I we prove that the value of the randomized problem does not change if we formulate the randomized McKean-Vlasov control problem on a new probability space.
\item[II)] Step II is devoted to the proof of the first inequality $V(t,x,\pi)\geq V^\Rc(t,x,\pi)$.
\begin{itemize}
\item[1)] In order to prove it, we construct in substep 1 a new probability space $(\check\Omega,\check\Fc,\check\P)$ for the randomized problem, which is a product space of $(\Omega,\Fc,\P)$ and a canonical space supporting the Poisson random measure. Step I guarantees that the value of the new randomized problem is still given by $V^\Rc(t,x,\pi)$.
\item[2)] In substep 2 we prove that the value of the original McKean-Vlasov control problem is still equal to $V(t,x,\pi)$ if we enlarge the class of admissible controls, taking all $\check\alpha\colon\check\Omega\times[0,T]\rightarrow A$ which are progressive with respect to the filtration $\check\F^{B,\mu_\infty}$. The new class of admissible controls is denoted $\check\Ac^{B,\mu_\infty}$.
\item[3)] In substep 3 we conclude the proof of the inequality $V(t,x,\pi)\geq V^\Rc(t,x,\pi)$, proving that for every $\check\nu\in\check\Vc$ there exists $\check\alpha^{\check\nu}\in\check\Ac^{B,\mu_\infty}$ such that $\check J^\Rc(t,x,\pi,\check\nu)=\check J(t,x,\pi,\check\alpha^{\check\nu})$. From substep 2, we immediately deduce that $V(t,x,\pi)\geq V^\Rc(t,x,\pi)$.
\end{itemize}
\item[III)] Step III is devoted to the proof of the other inequality $V(t,x,\pi)\leq V^\Rc(t,x,\pi)$. In few words, we prove that the set $\{\check\alpha^{\check\nu}\colon\check\nu\in\check\Vc\}$ is dense in $\check\Ac^{B,\mu_\infty}$ with respect to the distance $\tilde\rho$ in \eqref{KrylovMetric}. Then, the claim follows from the stability Lemma \ref{L:Stability}.
\end{itemize}

\noindent\textbf{Step I.} \emph{Value of the randomized McKean-Vlasov control problem.} Consider another probabilistic setting for the randomized problem, defined starting from $(\Omega,\Fc,\P)$, along the same lines as in Section \ref{S:RandomizedProblem}, where the objects $(\Omega^1,\Fc^1,\P^1)$, $(\bar\Omega,\bar\Fc,\bar\P)$, $\bar\Gc$, $\bar B$, $\bar\mu$, $T_n$, $\Ac_n$, $\bar I$, $\bar X^{t,\bar\xi}$, $\bar X^{t,x,\pi}$, $\Vc$, $J^\Rc(t,x,\pi,\nu)$, $V^\Rc(t,x,\pi)$ are replaced respectively by $(\check\Omega^1,\check\Fc^1,\check\P^1)$, $(\check\Omega,\check\Fc,\check\P)$, $\check\Gc$, $\check B$, $\check\mu$, $\check T_n$, $\check\Ac_n$, $\check I$, $\check X^{t,\check\xi}$, $\check X^{t,x,\pi}$, $\check\Vc$, $\check J^\Rc(t,x,\pi,\check\nu)$, $\check V^\Rc(t,x,\pi)$, with $\check\xi(\omega,\check\omega^1):=\xi(\omega)$, so that $\check\xi\in L^2(\check\Omega,\check\Gc,\check\P;\R^n)$ and $\pi=\P_{\text{\tiny$\check\xi$}}$ under $\check\P$.

We claim that $V^\Rc(t,x,\pi)=\check V^\Rc(t,x,\pi)$. Let us prove $V^\Rc(t,x,\pi)\leq\check V^\Rc(t,x,\pi)$, the other inequality can be proved in a similar way. We begin noting that $V^\Rc(t,x,\pi)\leq\check V^\Rc(t,x,\pi)$ follows if we prove that for every $\nu\in\Vc$ there exists $\check\nu\in\check\Vc$ such that $J^\Rc(t,x,\pi,\nu)=\check J^\Rc(t,x,\pi,\check\nu)$. Observe that
\[
J^\Rc(t,x,\pi,\nu) \ = \ \bar\E\bigg[\bar\kappa_T^\nu\bigg(\int_t^T f\big(s,\bar X_s^{t,x,\pi},\P_{\text{\tiny$\bar X_s^{t,\bar\xi}$}}^{\text{\tiny$\bar\Fc_s^\mu$}},\bar I_s\big)\,ds + g\big(\bar X_T^{t,x,\pi},\P_{\text{\tiny$\bar X_T^{t,\bar\xi}$}}^{\text{\tiny$\bar\Fc_T^\mu$}}\big)\bigg)\bigg].
\]
The quantity $J^\Rc(t,x,\pi,\nu)$ depends only on the joint law of $\bar\kappa_T^\nu$, $\bar X_\cdot^{t,x,\pi}$, $\P_{\text{\tiny$\bar X_\cdot^{t,\bar\xi}$}}^{\text{\tiny$\bar\Fc_\cdot^\mu$}}$, $\bar I_\cdot$ under $\bar \P$, which in turn depends on the joint law of $\bar B$, $\bar\mu$, $\bar\nu$ under $\bar\P$.

Recall that $\bar\nu_t(\omega,\omega^1,\alpha):=\nu_t(\omega^1,\alpha)$ and $\nu$ is $\Pc(\F^\mu)\otimes\Bc(\Ac)$-measurable. Then, we can suppose, using a monotone class argument, that $\nu$ is given by
\[
\nu_s(\alpha) \ = \ k(\alpha)1_{(T_n,T_{n+1}]}(s)\Psi(s,T_1,\ldots,T_n,\Ac_1,\ldots,\Ac_n),
\]
for some bounded and positive Borel-measurable maps $k$ and $\Psi$. We then see that $\check\nu$ defined by
\[
\check\nu_s(\alpha) \ := \ k(\alpha)1_{(\check T_n,\check T_{n+1}]}(s)\Psi(s,\check T_1,\ldots,\check T_n,\check\Ac_1,\ldots,\check \Ac_n)
\]
is such that $J^\Rc(t,x,\pi,\nu)=\check J^\Rc(t,x,\pi,\check\nu)$.

\vspace{2mm}

\noindent\textbf{Step II.} \emph{Proof of the inequality $V(t,x,\pi)\geq V^\Rc(t,x,\pi)$.} We shall exploit Proposition 4.1 in \cite{BCFP15}, for which we need to introduce a specific probabilistic setting for the randomized problem.

\vspace{1mm}

\noindent\textbf{Substep 1.} \emph{Canonical probabilistic setting for the randomized McKean-Vlasov control problem.}
Recall that the Polish space $\Ac$ can be countable or uncountable, and in this latter case it is Borel-isomorphic to $\R$ (see Corollary 7.16.1 in \cite{BertShreve78}). Then, in both cases, it can be proved (see the beginning of Section 4.1 in \cite{BCFP15}) that there exists a surjective measurable map $\iota\colon\R\rightarrow\Ac$ and a finite positive measure $\lambda'$ on $(\R,\Bc(\R))$ with full topological support, such that $\lambda=\lambda'\circ\iota^{-1}$ and $\lambda'$ is diffuse, namely $\lambda'(\{r\})=0$ for every $r\in\R$.

Now, consider the canonical probability space $(\Omega',\Fc',\P')$ of a marked point process on $\R_+\times\R$ associated to a Poisson random measure with compensator $\lambda'(dr)\,dt$. In other words, $\omega'\in\Omega'$ is a double sequence $\omega'=(t_n,r_n)_{n\geq1}\subset(0,\infty)\times\R$, with $t_n<t_{n+1}\nearrow\infty$. We denote by $(T_n',R_n')_{n\geq1}$ the canonical marked point process defined as $(T_n'(\omega'),R_n'(\omega'))=(t_n,r_n)$, and by $\zeta'=\sum_{n\geq1}\delta_{(T_n',R_n')}$ the canonical random measure. $\Fc'$ is the $\sigma$-algebra generated by the sequence $(T_n',R_n')_{n\geq1}$. $\P'$ is the unique probability on $\Fc'$ under which $\zeta'$ has compensator $\lambda'(dr)\,ds$. Finally, we complete $(\Omega',\Fc',\P')$ and, to simplify the notation, we still denote its completion by $(\Omega',\Fc',\P')$.

Set $\Ac_n'=\iota(R_n')$ and $\mu'=\sum_{n\geq1}\delta_{(T_n',\Ac_n')}$. Then $\mu'$ is a Poisson random measure on $(\Omega',\Fc',\P')$ with compensator $\lambda(d\alpha)\,ds$. Proceeding along the same lines as in Section \ref{S:RandomizedProblem}, we define, starting from $(\Omega,\Fc,\P)$ and $(\Omega',\Fc',\P')$, a new setting for the randomized problem where the objects $(\Omega^1,\Fc^1,\P^1)$, $(\bar\Omega,\bar\Fc,\bar\P)$, $\bar\Gc$, $\bar B$, $\bar\mu$, $\bar\F^B=(\bar\Fc_s^B)_{s\geq0}$, $\F^\mu=(\Fc_s^\mu)_{s\geq0}$, $\bar\F^{B,\mu}=(\bar\Fc_s^{B,\mu})_{s\geq0}$, $(T_n,\Ac_n)_{n\geq1}$, $\bar I$, $\bar X^{t,\bar\xi}$, $\bar X^{t,x,\pi}$, $\Vc$, $\P^\nu$, $\bar\P^\nu$, $J^\Rc(t,x,\pi,\nu)$, $V^\Rc(t,x,\pi)$ are replaced respectively by $(\Omega',\Fc',\P')$, $(\check\Omega,\check\Fc,\check\P)$, $\check\Gc$, $\check B$, $\check\mu$, $\check\F^B=(\check\Fc_s^B)_{s\geq0}$, $\F^{\mu'}=(\Fc_s^{\mu'})_{s\geq0}$, $\check\F^{B,\mu}=(\check\Fc_s^{B,\mu})_{s\geq0}$, $(\check T_n,\check\Ac_n)_{n\geq1}$, $\check I$, $\check X^{t,\check\xi}$, $\check X^{t,x,\pi}$, $\check\Vc$, $\P^{\check\nu}$, $\check\P^{\check\nu}$, $\check J^\Rc(t,x,\pi,\check\nu)$, $\check V^\Rc(t,x,\pi)$, with $\check\xi(\omega,\omega'):=\xi(\omega)$, so that $\check\xi\in L^2(\check\Omega,\check\Gc,\check\P;\R^n)$ and $\pi=\P_{\text{\tiny$\check\xi$}}$ under $\check\P$.

\vspace{1mm}

\noindent\textbf{Substep 2.} \emph{Value of the original McKean-Vlasov control problem.}
$\check\F^{B,\mu_\infty}=(\check\Fc_s^{B,\mu_\infty})_{s\geq0}$ be the $\check\P$-completion of the filtration $(\Fc_s^B\otimes\Fc')_{s\geq0}$, and $\check\Fc'$ the canonical extension of $\Fc'$ to $\check\Omega$. We define the set $\check\Ac^{B,\mu_\infty}$ of all $\check\F^{B,\mu_\infty}$-progressive processes $\check\alpha\colon\check\Omega\times[0,T]\rightarrow A$. For every $\check\alpha\in\check\Ac^{B,\mu_\infty}$, we denote $(\check X_s^{t,\check\xi,\check\alpha},\check X_s^{t,x,\pi,\check\alpha})_{s\in[t,T]}$ the unique continuous $(\check\Fc_s^{B,\mu_\infty}\vee\check\Gc)_s$-adapted solution to the following system of equations:
\begin{align}
d\check X_s^{t,\check\xi,\check\alpha} \ &= \ b\big(s,\check X_s^{t,\check\xi,\check\alpha},\P_{\text{\tiny$\check X_s^{t,\check\xi,\check\alpha}$}}^{\text{\tiny$\check\Fc_s^\mu$}},\check\alpha_s\big)\,ds + \sigma\big(s,\check X_s^{t,\check\xi,\check\alpha},\P_{\text{\tiny$\check X_s^{t,\check\xi,\check\alpha}$}}^{\text{\tiny$\check\Fc_s^\mu$}},\check\alpha_s\big)\,d\check B_s, \quad &\check X_t^{t,\check\xi,\check\alpha} = \ \check\xi, \label{StateEq1_mu_infty} \\
d\check X_s^{t,x,\pi,\check\alpha} \ &= \ b\big(s,\check X_s^{t,x,\pi,\check\alpha},\P_{\text{\tiny$\check X_s^{t,\check\xi,\check\alpha}$}}^{\text{\tiny$\check\Fc_s^\mu$}},\check\alpha_s\big)\,ds + \sigma\big(s,\check X_s^{t,x,\pi,\check\alpha},\P_{\text{\tiny$\check X_s^{t,\check\xi,\check\alpha}$}}^{\text{\tiny$\check\Fc_s^\mu$}},\check\alpha_s\big)\,d\check B_s, &\check X_t^{t,x,\pi,\check\alpha} = \ x, \label{StateEq2_mu_infty}
\end{align}
for all $s\in[t,T]$, where $\P_{\text{\tiny$\check X_s^{t,\check\xi,\check\alpha}$}}^{\text{\tiny$\check\Fc_s^\mu$}}$ denotes the regular conditional distribution of the random variable $\check X_s^{t,\check\xi,\check\alpha}\colon\check\Omega\rightarrow\R^n$ with respect to $\check\Fc_s^\mu$. We also define ($\check\E$ denotes the $\check\P$-expected value)
\[
\check J(t,x,\pi,\check\alpha) \ = \ \check\E\bigg[\int_t^T f\big(s,\check X_s^{t,x,\pi,\check\alpha},\P_{\text{\tiny$\check X_s^{t,\check\xi,\check\alpha}$}}^{\text{\tiny$\check\Fc_s^\mu$}},\check\alpha_s\big)\,ds + g\big(\check X_T^{t,x,\pi,\check\alpha},\P_{\text{\tiny$\check X_T^{t,\check\xi,\check\alpha}$}}^{\text{\tiny$\check\Fc_T^\mu$}}\big)\bigg],
\]
and 
\[
\check V(t,x,\pi) \ = \ \sup_{\check\alpha\in\check\Ac^{B,\mu_\infty}} \check J(t,x,\pi,\check\alpha).
\]
Let us prove that $V(t,x,\pi)=\check V(t,x,\pi)$. 

The inequality $V(t,x,\pi)\leq\check V(t,x,\pi)$ is obvious. Indeed, every $\alpha\in\Ac$ admits an obvious extension $\check\alpha(\omega,\omega'):=\alpha(\omega)$ to $\check\Omega$. Notice that $\check\alpha\in\check\Ac^{B,\mu_\infty}$. We also observe that $\check X_s^{t,\check\xi,\check\alpha}(\omega,\omega')=X_s^{t,\xi,\alpha}(\omega)$, for $\check\P$-almost every $(\omega,\omega')\in\check\Omega$. Therefore $\P_{\text{\tiny$\check X_s^{t,\check\xi,\check\alpha}$}}^{\text{\tiny$\check\Fc_s^\mu$}}$ is equal $\check\P$-a.s. to $\P_{\text{\tiny$X_s^{t,\xi,\alpha}$}}$. Then, $\check X_s^{t,x,\pi,\check\alpha}(\omega,\omega')=X_s^{t,x,\pi,\alpha}(\omega)$, for $\check\P$-almost every $(\omega,\omega')\in\check\Omega$. As a consequence, we see that $J(t,x,\pi,\alpha)=\check J(t,x,\pi,\check\alpha)$.

To prove the other inequality, let $\tilde\alpha\in\check\Ac^{B,\mu_\infty}$. Then, there exists an $A$-valued $(\Fc_s^B\otimes\Fc')_{s\geq0}$-progressive process $\check\alpha\colon\check\Omega\times[0,T]\rightarrow A$ satisfying $\check\alpha=\tilde\alpha$, $d\check\P\,ds$-a.e., so that $\check J(t,x,\pi,\check\alpha)=\check J(t,x,\pi,\tilde\alpha)$. Moreover, for every $\omega'\in\Omega'$ the process $\alpha^{\omega'}$, given by
$\alpha_s^{\omega'}(\omega):=\check\alpha_s(\omega,\omega')$, is $\F^B$-progressive.

Now, for every $\omega'\in\Omega'$, consider the solution $(X_s^{t,\xi,\alpha^{\omega'}},X_s^{t,x,\pi,\alpha^{\omega'}})_{s\in[t,T]}$ to \eqref{StateEq1}-\eqref{StateEq2} with $\alpha$ replaced by $\alpha^{\omega'}$, namely
\begin{align*}
dX_s^{t,\xi,\alpha^{\omega'}} \ &= \ b\big(s,X_s^{t,\xi,\alpha^{\omega'}},\P_{\text{\tiny$X_s^{t,\xi,\alpha^{\omega'}}$}},\alpha_s^{\omega'}\big)\,ds + \sigma\big(s,X_s^{t,\xi,\alpha^{\omega'}},\P_{\text{\tiny$X_s^{t,\xi,\alpha^{\omega'}}$}},\alpha_s^{\omega'}\big)\,dB_s, \\
dX_s^{t,x,\pi,\alpha^{\omega'}} \ &= \ b\big(s,X_s^{t,x,\pi,\alpha^{\omega'}},\P_{\text{\tiny$X_s^{t,\xi,\alpha^{\omega'}}$}},\alpha_s^{\omega'}\big)\,ds + \sigma\big(s,X_s^{t,x,\pi,\alpha^{\omega'}},\P_{\text{\tiny$X_s^{t,\xi,\alpha^{\omega'}}$}},\alpha_s^{\omega'}\big)\,dB_s.
\end{align*}
On the other hand, since $(\check X_s^{t,\check\xi,\check\alpha},\check X_s^{t,x,\pi,\check\alpha})_{s\in[t,T]}$ is the solution to \eqref{StateEq1_mu_infty}-\eqref{StateEq2_mu_infty}, we have, for $\P'$-a.e. $\omega'\in\Omega'$,
\begin{align*}
d\check X_s^{t,\check\xi,\check\alpha}(\cdot,\omega') \ &= \ b\big(s,\check X_s^{t,\check\xi,\check\alpha}(\cdot,\omega'),\P_{\text{\tiny$\check X_s^{t,\check\xi,\check\alpha}$}}^{\text{\tiny$\check\Fc_s^\mu$}}(\cdot,\omega'),\check\alpha_s(\cdot,\omega')\big)\,ds \\
&\quad \ + \sigma\big(s,\check X_s^{t,\check\xi,\check\alpha}(\cdot,\omega'),\P_{\text{\tiny$\check X_s^{t,\check\xi,\check\alpha}$}}^{\text{\tiny$\check\Fc_s^\mu$}}(\cdot,\omega'),\check\alpha_s(\cdot,\omega')\big)\,dB_s, \\
d\check X_s^{t,x,\pi,\check\alpha}(\cdot,\omega') \ &= \ b\big(s,\check X_s^{t,x,\pi,\check\alpha}(\cdot,\omega'),\P_{\text{\tiny$\check X_s^{t,\check\xi,\check\alpha}$}}^{\text{\tiny$\check\Fc_s^\mu$}}(\cdot,\omega'),\check\alpha_s(\cdot,\omega')\big)\,ds \\
&\quad \ + \sigma\big(s,\check X_s^{t,x,\pi,\check\alpha}(\cdot,\omega'),\P_{\text{\tiny$\check X_s^{t,\check\xi,\check\alpha}$}}^{\text{\tiny$\check\Fc_s^\mu$}}(\cdot,\omega'),\check\alpha_s(\cdot,\omega')\big)\,dB_s.
\end{align*}
Notice that, for $\P'$-a.e. $\omega'\in\Omega'$ we have that $\P_{\text{\tiny$\check X_s^{t,\check\xi,\check\alpha}$}}^{\text{\tiny$\check\Fc_s^\mu$}}(\cdot,\omega')$ is equal $\P$-a.s. to $\P_{\text{\tiny$\check X_s^{t,\check\xi,\check\alpha}(\cdot,\omega')$}}$, the law under $\P$ of the random variable $\check X_s^{t,\check\xi,\check\alpha}(\cdot,\omega')\colon\Omega\rightarrow\R^n$ 

Recalling the identity $\alpha_s^{\omega'}=\check\alpha_s(\cdot,\omega')$, we see that, for $\P'$-a.e. $\omega'\in\Omega'$, $(X_s^{t,\xi,\alpha^{\omega'}},X_s^{t,x,\pi,\alpha^{\omega'}})_{s\in[t,T]}$ and $(\check X_s^{t,\check\xi,\check\alpha}(\cdot,\omega'),\check X_s^{t,x,\pi,\check\alpha}(\cdot,\omega'))_{s\in[t,T]}$ solve the same system of equations. Then, by pathwise uniqueness, for $\P'$-a.e. $\omega'\in\Omega'$, we have $X_s^{t,\xi,\alpha^{\omega'}}(\omega)=\check X_s^{t,\check\xi,\check\alpha}(\omega,\omega')$ and $X_s^{t,x,\pi,\alpha^{\omega'}}(\omega)=\check X_s^{t,x,\pi,\check\alpha}(\omega,\omega')$, for all $s\in[t,T]$, $\P(d\omega)$-almost surely. Therefore, by Fubini's theorem,
\begin{align*}
\check J(t,x,\pi,\check\alpha) \ &= \ \int_{\Omega'} \E\bigg[\int_t^Tf\big(s,X_s^{t,x,\pi,\alpha^{\omega'}},\P_{\text{\tiny$X_s^{t,\xi,\alpha^{\omega'}}$}},\alpha_s^{\omega'}\big)\,ds+g\big(X_T^{t,x,\pi,\alpha^{\omega'}},\P_{\text{\tiny$X_T^{t,\xi,\alpha^{\omega'}}$}}\big)
\bigg]\,\P'(d\omega') \\
&= \ \int_{\Omega'} J(t,x,\pi,\alpha^{\omega'})\,\P'(d\omega') \ \leq \ V(t,x,\pi).
\end{align*}
Recalling that $\check J(t,x,\pi,\tilde\alpha)=\check J(t,x,\pi,\check\alpha)$, we deduce that $\check J(t,x,\pi,\tilde\alpha)\leq V(t,x,\pi)$. Taking the supremum over $\tilde\alpha\in\check\Ac^{B,\mu_\infty}$, we conclude that $\check V(t,x,\pi)\leq V(t,x,\pi)$.

\vspace{1mm}

\noindent\textbf{Substep 3.} \emph{Proof of the inequality $V(t,x,\pi)\geq V^\Rc(t,x,\pi)$.}
Let $\check\nu\in\check\Vc$. By Lemma 4.3 in \cite{BCFP15} there exists a sequence $(\check T_n^{\check\nu},\check\Ac_n^{\check\nu})_{n\geq1}$ on $(\Omega',\Fc',\P')$ such that:
\begin{itemize}
\item $(\check T_n^{\check\nu},\check\Ac_n^{\check\nu})$ takes values in $(0,\infty)\times\Ac$;
\item $\check T_n^{\check\nu}<\check T_{n+1}^{\check\nu}\nearrow\infty$;
\item $\check T_n^{\check\nu}$ is an $\F^{\mu'}$-stopping time and $\check\Ac_n^{\check\nu}$ is $\Fc_{\check T_n^{\check\nu}}^{\mu'}$-measurable;
\item the law of $(\check T_n^{\check\nu},\check\Ac_n^{\check\nu})_{n\geq1}$ under $\P'$ coincides with the law of $(\check T_n,\check\Ac_n)_{n\geq1}$ under $\P^{\check\nu}$.
\end{itemize}
Let $\check\alpha^{\check\nu}\colon\check\Omega\times[0,T]\rightarrow A$ be given by ($\bar\alpha$ was introduced in \eqref{barI})
\[
\check\alpha_s^{\check\nu}(\omega,\omega') \ = \ \bar\alpha_s(\omega)\,1_{[0,\check T_1^{\check\nu}(\omega'))}(s) + \sum_{n\ge 1} (\check\Ac_n^{\check\nu}(\omega'))_{s\wedge T}(\omega)\,1_{[\check T_n^{\check\nu}(\omega'),\check T_{n+1}^{\check\nu}(\omega'))}(s).
\]
Notice that $\check\alpha^{\check\nu}\in\check\Ac^{B,\mu_\infty}$. For every $n\geq1$, set $\check\alpha_{n,s}(\omega,\omega'):=(\check\Ac_n(\omega'))_s(\omega)$ and $\check\alpha_{n,s}^{\check\nu}(\omega,\omega'):=(\check\Ac_n^{\check\nu}(\omega'))_s(\omega)$, for all $s\in[0,T]$. Notice that the law of $(\check\alpha_{n,s})_{s\in[0,T]}$ under $\check\P^{\check\nu}$ coincides with the law of $(\check\alpha_{n,s}^{\check\nu})_{s\in[0,T]}$ under $\check\P$ (to see this, we can suppose, by an approximation argument, that the $\Ac$-valued random variables $\check\Ac_n$ and $\check\Ac_n^{\check\nu}$ take only a finite number of values). It follows that the law of $\check I$ under $\check\P^{\check\nu}$ coincides with the law of $\check\alpha^{\check\nu}$ under $\check\P$.  

More generally, for every $n\geq1$, the law of $(\check\xi,\check B,\check\alpha_{n,\cdot})$ under $\check\P^{\check\nu}$ is equal to the law of $(\check\xi,\check B,\check\alpha_{n,\cdot}^{\check\nu})$ under $\check\P$. Therefore, the law of $(\check\xi,\check B,\check I)$ under $\check\P^{\check\nu}$ coincides with the law of $(\check\xi,\check B,\check\alpha^{\check\nu})$ under $\check\P$. This implies that the law of $(\check X^{t,\check\xi},\check X^{t,x,\pi},\check I)$ under $\check\P^{\check\nu}$ is equal to the law of $(\check X^{t,\check\xi,\check\alpha^{\check\nu}},\check X^{t,x,\pi,\check\alpha^{\check\nu}},\check\alpha^{\check\nu})$ under $\check\P$. It follows that $\check J^\Rc(t,x,\pi,\check\nu)=\check J(t,x,\pi,\check\alpha^{\check\nu})$. In particular, we have
\[
\sup_{\check\nu\in\check\Vc}\check J^\Rc(t,x,\pi,\check\nu) \ = \ \sup_{\substack{\check\alpha^{\check\nu}\\\check\nu\in\check\Vc}}\check J(t,x,\pi,\check\alpha^{\check\nu}).
\]
Since the left-hand side is equal to $\check V^\Rc(t,x,\pi)$, while the right-hand side is clearly less than or equal to $\check V(t,x,\pi)$, we get $\check V^\Rc(t,x,\pi)\leq\check V(t,x,\pi)$. Recalling from step I that $V^\Rc(t,x,\pi)=\check V^\Rc(t,x,\pi)$ and from substep 2 that $\check V(t,x,\pi)=V(t,x,\pi)$, we conclude $V^\Rc(t,x,\pi)\leq V(t,x,\pi)$.

\vspace{2mm}

\noindent\textbf{Step III.} \emph{Proof of the inequality $V(t,x,\pi)\leq V^\Rc(t,x,\pi)$.} The proof of this step is based on Proposition A.1 in \cite{BCFP15} (notice, however, that we will need to use some results from the proof of this Proposition, not only from its statement). More precisely, the set $\Omega$ appearing in Proposition A.1 of \cite{BCFP15} is the empty set $\Omega=\emptyset$ in our context, so that the product probability space $(\tilde\Omega,\tilde\Fc,\Q)$ coincides with $(\Omega',\Fc',\P')$, which is some suitably defined probability space (see Appendix A in \cite{BCFP15} for the definition of $(\Omega',\Fc',\P')$; here, we do not need to know the structure of $(\Omega',\Fc',\P')$). Fix $\hat\alpha\in\Ac$ and denote by $\boldsymbol\alpha\colon[0,T]\rightarrow\Ac$ the map $\boldsymbol\alpha_s=\hat\alpha$, for every $s\in[0,T]$. By Proposition A.1 in \cite{BCFP15} we have that, for every $\ell\in\N\backslash\{0\}$, there exists a marked point process $( T_n^\ell, \Ac_n^\ell)_{n\geq1}$ on $(\Omega',\Fc',\P')$ such that ($\bar\alpha$ was introduced in \eqref{barI})
\[
 T_0^\ell \ = \ 0, \qquad\quad  \Ac_0^\ell \ = \ \bar\alpha, \qquad\quad  \Ic_s^\ell(\omega') \ = \ \sum_{n\geq0}  \Ac_n^\ell(\omega')\,1_{[ T_n^\ell(\omega'), T_{n+1}^\ell(\omega'))}(s), \quad \text{for all }s\geq0
\]
and
\begin{equation}\label{KrilovMetric<1/m}
\E'\bigg[\int_0^T\tilde\rho( \Ic_s^\ell,\boldsymbol\alpha_s)\,ds\bigg] \ < \ \frac{1}{\ell},
\end{equation}
where $\E'$ denotes the $\P'$-expected value. Set $\mu_\ell=\sum_{n\geq1}\delta_{( T_n^\ell, \Ac_n^\ell)}$ the random measure associated to $( T_n^\ell, \Ac_n^\ell)_{n\geq 1}$, and denote $\F^{\mu_\ell}=(\Fc_s^{\mu_\ell})_{s\geq0}$ the filtration generated by $\mu_\ell$. Then, by Proposition A.1 of \cite{BCFP15} we have that the $\F^{\mu_\ell}$-compensator of $\mu_\ell$ under $\P'$ is given by $\nu_s^\ell(\alpha)\,\lambda(d\alpha)\,ds$ for some $\Pc(\F^{\mu_\ell})\otimes\Bc(\Ac)$-measurable map $\nu^\ell\colon\Omega'\times\R_+\times\Ac\rightarrow\R_+$ satisfying
\begin{equation}\label{BoundsNu}
0 \ < \ \inf_{\Omega'\times[0,T]\times\Ac}\nu^\ell \ \leq \ \sup_{\Omega'\times[0,T]\times\Ac}\nu^\ell \ < \ \infty.
\end{equation}
Noting that the definition of $\nu^\ell$ on $\Omega'\times(T,\infty)\times\Ac$ is not relevant in order to guarantee \eqref{KrilovMetric<1/m}, we can assume that $\nu^\ell\equiv1$ on $\Omega'\times(T,\infty)\times\Ac$.

Observe that
\[
\E'\bigg[\int_0^T\tilde\rho( \Ic_s^\ell,\boldsymbol\alpha_s)\,ds\bigg] \ = \ \sum_{n\geq0} \E'\bigg[ 1_{\{ T_n^\ell< T\}} \int_{ T_n^\ell}^{ T_{n+1}^\ell\wedge T}\E\bigg[\int_0^T \rho((\Ac_n^\ell)_r,\hat\alpha_r)\,dr\bigg] ds\bigg] \ < \ \frac{1}{\ell}.
\]
On the other hand, let
\[
\tilde I_s^\ell(\omega,\omega') \ = \ \sum_{n\geq0} ( \Ac_n^\ell(\omega'))_{s\wedge T}(\omega)\,1_{[ T_n^\ell(\omega'), T_{n+1}^\ell(\omega'))}(s), \qquad \text{for all }s\geq0.
\]
Our aim is to prove that
\begin{align}\label{tilderhoQ,0}
&\tilde\rho^\Q(\tilde I^\ell,\hat\alpha) \ := \ \E'\bigg[\E\bigg[\int_0^T\rho(\tilde I_r^\ell,\hat\alpha_r)\,dr\bigg]\bigg] \ \overset{\ell\rightarrow\infty}{\longrightarrow} \ 0.
\end{align}
\textbf{Digression.} \emph{Estimate for the series $\sum_{n\geq0}\P'( T_n^\ell<T)$.} We recall from the proof of Proposition A.1 in \cite{BCFP15} that the sequence $( T_n^\ell)_{n\geq0}$ is the disjoint union of $(R_n^m)_{n\geq1}$ and $(T_n^k)_{n\geq0}$ (we refer to the proof of Proposition A.1 in \cite{BCFP15} for all unexplained notations), namely
\begin{equation}\label{seriesT_n^ell}
\sum_{n\geq0}\P'\big( T_n^\ell<T\big) \ = \ \sum_{n\geq1}\P'\big(R_n^m<T\big) + \sum_{n\geq0}\P'\big(T_n^k<T\big).
\end{equation}
We also recall that $T_n^k-T_{n-1}^k$ has an exponential distribution with parameter $k^{-1}\lambda(\Ac)$. Then, it is easy to prove by induction on $n$, the estimate
\begin{equation}\label{seriesT_n^k}
\P'\big(T_n^k<T\big) \ \leq \ \big(1 - e^{-k^{-1}\lambda(\Ac) T}\big)^n.
\end{equation}
On the other hand, concerning the sequence $(R_n^m)_{n\geq1}$, we begin noting that since $\boldsymbol\alpha$ is constant and identically equal to $\hat\alpha$, the sequence of deterministic times $(t_n)_{n\geq0}$ appearing in the proof of Proposition A.1 in \cite{BCFP15} can be taken as follows: $t_0=0$, $t_1\in(0,\frac{1}{3\ell}\wedge T)$, and $t_n=T+n-2$ for every $n\geq2$. Therefore $R_n^m\geq T$ for all $n\geq2$, while $R_1^m=t_1+V_1^m$, where $V_1^m$ is an exponential random variable with parameter $\lambda_{1m}>m$. In particular, we have
\begin{equation}\label{seriesR_n^m}
\P'\big(R_1^m<T\big) \ = \ \P'\big(V_1^m< T - t_1\big) \ = \ 1 - e^{-\lambda_{1m}(T - t_1)} \ \leq \ 1.
\end{equation}
Plugging \eqref{seriesT_n^k} and \eqref{seriesR_n^m} into \eqref{seriesT_n^ell}, we obtain
\begin{equation}\label{seriesT_n^ell_2}
\sum_{n\geq0}\P'\big( T_n^\ell<T\big) \ \leq \ 1 + \sum_{n\geq0} \big(1 - e^{-k^{-1}\lambda(\Ac) T}\big)^n \ \leq \ 1 + e^{k^{-1}\lambda(\Ac) T} \ \leq \ 1 + e^{\lambda(\Ac) T}.
\end{equation}
\textbf{Continuation of the proof of Step III.} We can now prove \eqref{tilderhoQ,0}. In particular, we have, using \eqref{seriesT_n^ell_2},
\begin{align*}
&\tilde\rho^\Q(\tilde I^\ell,\hat\alpha) \ = \ \E'\bigg[\E\bigg[\int_0^T\rho(\tilde I_r^\ell,\hat\alpha_r)\,dr\bigg]\bigg] \ = \ \sum_{n\geq0}\E'\bigg[1_{\{ T_n^\ell<T\}} \E\bigg[\int_{ T_n^\ell}^{ T_{n+1}^\ell\wedge T} \rho((\Ac_n^\ell)_r,\hat\alpha_r)\,dr\bigg] \bigg] \notag \\
&= \ \sum_{n\geq0}\E'\bigg[1_{\{ T_n^\ell< T\}} \frac{1}{ T_{n+1}^\ell\wedge T -  T_n^\ell}\int_{ T_n^\ell}^{ T_{n+1}^\ell\wedge T}\E\bigg[\int_{ T_n^\ell}^{ T_{n+1}^\ell\wedge T} \rho((\Ac_n^\ell)_r,\hat\alpha_r)\,dr\bigg]ds \bigg] \notag \\
&= \ \sum_{n\geq0}\E'\bigg[1_{\{ T_{n+1}^\ell\wedge T -  T_n^\ell\geq1/\sqrt{\ell}\}} 1_{\{ T_n^\ell< T\}} \frac{1}{ T_{n+1}^\ell\wedge T -  T_n^\ell} \int_{ T_n^\ell}^{ T_{n+1}^\ell\wedge T}\E\bigg[\int_{ T_n^\ell}^{ T_{n+1}^\ell\wedge T} \!\!\!\!\!\!\!\! \rho((\Ac_n^\ell)_r,\hat\alpha_r)\,dr\bigg]ds \bigg] \notag \\
&\quad \ + \sum_{n\geq0}\E'\bigg[1_{\{ T_{n+1}^\ell\wedge T -  T_n^\ell<1/\sqrt{\ell}\}} 1_{\{ T_n^\ell< T\}} \frac{1}{ T_{n+1}^\ell\wedge T -  T_n^\ell} \int_{ T_n^\ell}^{ T_{n+1}^\ell\wedge T}\E\bigg[\int_{ T_n^\ell}^{ T_{n+1}^\ell\wedge T} \!\!\!\!\!\!\!\! \rho((\Ac_n^\ell)_r,\hat\alpha_r)\,dr\bigg]ds \bigg] \notag \\
&\leq \ \sqrt{\ell}\,\sum_{n\geq0}\E'\bigg[ 1_{\{ T_n^\ell< T\}} \int_{ T_n^\ell}^{ T_{n+1}^\ell\wedge T}\E\bigg[\int_0^T \rho((\Ac_n^\ell)_r,\hat\alpha_r)\,dr\bigg]ds \bigg] + \frac{1}{\sqrt{\ell}} \sum_{n\geq0} \P'\big( T_n^\ell<T\big) \notag \\
&= \ \sqrt{\ell}\,\E'\bigg[\int_0^T\tilde\rho(\Ic_s^\ell,\boldsymbol\alpha_s)\,ds\bigg] + \frac{1}{\sqrt{\ell}} \sum_{n\geq0} \P'\big( T_n^\ell<T\big) \notag \\
&\leq \ \sqrt{\ell}\,\E'\bigg[\int_0^T\tilde\rho(\Ic_s^\ell,\boldsymbol\alpha_s)\,ds\bigg] + \frac{1 + e^{\lambda(\Ac)T}}{\sqrt{\ell}} \ \leq \frac{2 + e^{\lambda(\Ac)T}}{\sqrt{\ell}},
\end{align*}
which yields \eqref{tilderhoQ,0}.

We consider now the product probability space $(\Omega\times\Omega',\Fc\otimes\Fc',\P\otimes\P')$, which we still denote $(\tilde\Omega,\tilde\Fc,\Q)$ (by an abuse of notation, since according to Proposition A.1 in \cite{BCFP15}, $(\tilde\Omega,\tilde\Fc,\Q)$ coincides with $(\Omega',\Fc',\P')$). We complete the probability space $(\tilde\Omega,\tilde\Fc,\Q)$ and, to simplify the notation, we still denote by $(\tilde\Omega,\tilde\Fc,\Q)$ its completion. Let $\tilde\xi$, $\tilde B$, $\tilde\nu^\ell$  be the canonical extensions of $\xi$, $B$, $\nu^\ell$ to $\tilde\Omega$. On the other hand, we still denote by $\mu_\ell$ the extension of $\mu_\ell$ to $\tilde\Omega$. We denote by $\tilde\mu_\ell(ds\,d\alpha)=\mu_\ell(ds\,d\alpha)-\tilde\nu_s^\ell(\alpha)\lambda(d\alpha)ds$ the compensated martingale measure associated to $\mu^\ell$. We also denote by $\tilde\F^{B,\mu_\ell}=(\tilde\Fc_s^{B,\mu_\ell})_{s\geq0}$ (resp. $\tilde\F^{\mu_\ell}=(\tilde\Fc_s^{\mu_\ell})_{s\geq0}$) the $\Q$-completion of the filtration generated by $\tilde B$ and $\mu_\ell$ (resp. $\mu_\ell$). For every $\ell\in\N\backslash\{0\}$, we define the Dol\'eans exponential
\[
\tilde\kappa_s^\ell \ = \ \Ec_s\bigg(\int_0^\cdot\int_A (\tilde\nu_r^\ell(\alpha)^{-1} - 1)\,\tilde\mu_\ell(dr\, d\alpha)\bigg), \qquad \text{for all }s\in[0,T].
\]
By \eqref{BoundsNu} we see that $(\tilde\kappa_s^\ell)_{s\in[0,T]}$ is an $\tilde\F^{B,\mu_\ell}$-martingale under $\Q$, so that we can
define on $(\tilde\Omega,\tilde\Fc)$ a probability $\tilde\P_\ell$ equivalent to $\Q$ by $d\tilde\P_\ell$ $=$ $\tilde\kappa_T^\ell\,d\Q$. By the Girsanov theorem, $\mu_\ell$ has $\tilde\F^{B,\mu_\ell}$-compensator given by $\lambda(d\alpha)\,ds$ under $\tilde\P_\ell$. Moreover, $\tilde B$ remains a Brownian motion under $\tilde\P_\ell$, and $\pi=\P_{\tilde\xi}$ under $\tilde\P_\ell$.

Let $\tilde\Gc$ be the canonical extension of $\Gc$ to $\tilde\Omega$ and denote $(\tilde X_s^{t,\tilde\xi,\ell},\tilde X_s^{t,x,\pi,\ell})_{s\in[t,T]}$ the unique continuous $(\tilde\Fc_s^{B,\mu_\ell}\vee\tilde\Gc)$-adapted solution to equations \eqref{StateEq1Rand}-\eqref{StateEq2Rand} on $(\tilde\Omega,\tilde\Fc,\tilde\P_\ell)$ with $\bar\xi$, $\bar B$, $\bar I$, $\bar\Fc_s^\mu$ replaced by $\tilde\xi$, $\tilde B$, $\tilde I^\ell$, $\tilde\Fc_s^{\mu_\ell}$. Finally, we define in an obvious way the following objects: $\tilde\Vc_\ell$, $\tilde\P_\ell^{\tilde\nu}$, $\tilde\E_\ell^{\tilde\nu}$, $\tilde J_\ell^\Rc(t,x,\pi,\tilde\nu)$, $\tilde V_\ell^\Rc(t,x,\pi)$.

For every $\ell$ we have constructed a new probabilistic setting for the randomized problem, where the objects $(\Omega^1,\Fc^1,\P^1)$, $(\Omega,\Fc,\P)$, $\bar\Gc$, $\bar B$, $\bar\mu$, $\bar I$, $\bar X^{t,\bar\xi}$, $\bar X^{t,x,\pi}$, $\Vc$, $J^\Rc(t,x,\pi,\nu)$, $V^\Rc(t,x,\pi)$ are replaced respectively by $(\Omega',\Fc',\P')$, $(\tilde\Omega,\tilde\Fc,\tilde\P_\ell)$, $\tilde\Gc$, $\tilde B$, $\mu^\ell$, $\tilde I^\ell$, $\tilde X^{t,\tilde\xi,\ell}$, $\tilde X^{t,x,\pi,\ell}$, $\tilde\Vc_\ell$, $\tilde J_\ell^\Rc(t,x,\pi,\tilde\nu)$, $\tilde V_\ell^\Rc(t,x,\pi)$.

Now, let us prove that $\tilde J_\ell^\Rc(t,x,\pi,\tilde\nu^\ell)  \rightarrow J(t,x,\pi,\hat\alpha)$ as $\ell\rightarrow\infty$. To this end, notice that $\tilde\P_\ell^{\tilde\nu^\ell}\equiv\Q$. Therefore $\tilde J_\ell^\Rc(t,x,\pi,\tilde\nu^\ell)$ can be written in terms of $\E^\Q$ as follows:
\[
\tilde J_\ell^\Rc(t,x,\pi,\tilde\nu^\ell) \ = \ \E^\Q\bigg[\int_t^T f\big(s,\tilde X_s^{t,x,\pi,\ell},\P_{\text{\tiny$\tilde X_s^{t,\tilde\xi,\ell}$}}^{\text{\tiny$\tilde\Fc_s^{\mu_\ell}$}},\tilde I_s^\ell\big)\,ds + g\big(\tilde X_T^{t,x,\pi,\ell},\P_{\text{\tiny$\tilde X_T^{t,\tilde\xi,\ell}$}}^{\text{\tiny$\tilde\Fc_T^{\mu_\ell}$}}\big)\bigg].
\]
On the other hand, let $\tilde\F^B=(\tilde\Fc_s^B)_{s\geq0}$ be the $\Q$-completion of the filtration generated by $\tilde B$, and $\tilde\alpha$ the canonical extension of $\hat\alpha$ to $\tilde\Omega$. Then, we denote by $(\tilde X_s^{t,\tilde\xi,\tilde\alpha},\tilde X_s^{t,x,\pi,\tilde\alpha})_{s\in[t,T]}$ the unique continuous $(\tilde\Fc_s^B\vee\tilde\Gc)$-adapted solution to equations \eqref{StateEq1}-\eqref{StateEq2} on $(\tilde\Omega,\tilde\Fc,\Q)$ with $\xi$, $B$, $\alpha$ replaced by $\tilde\xi$, $\tilde B$, $\tilde\alpha$. Notice that $(\tilde X_s^{t,\tilde\xi,\tilde\alpha},\tilde X_s^{t,x,\pi,\tilde\alpha})_{s\in[t,T]}$ coincides with the obvious extension of $(X_s^{t,\xi,\hat\alpha},X_s^{t,x,\pi,\hat\alpha})_{s\in[t,T]}$ to $\tilde\Omega$. Hence, we have
\[
J(t,x,\pi,\hat\alpha) \ = \ \E^\Q\bigg[\int_t^T f\big(s,\tilde X_s^{t,x,\pi,\tilde\alpha},\P_{\text{\tiny$\tilde X_s^{t,\tilde\xi,\tilde\alpha}$}}^{\text{\tiny$\tilde\Fc_s^{\mu_\ell}$}},\tilde\alpha_s\big)\,ds + g\big(\tilde X_T^{t,x,\pi,\tilde\alpha},\P_{\text{\tiny$\tilde X_T^{t,\tilde\xi,\tilde\alpha}$}}^{\text{\tiny$\tilde\Fc_T^{\mu_\ell}$}}\big)\bigg].
\]
Then, it follows that $\tilde J_\ell^\Rc(t,x,\pi,\tilde\nu^\ell)  \rightarrow J(t,x,\pi,\hat\alpha)$ as $\ell\rightarrow\infty$. Indeed, this is a direct consequence of Lemma \ref{L:Stability}, with $\tilde\F^{\mu_0}:=(\{\emptyset,\tilde\Omega\})_{s\geq0}$ being the trivial filtration, $\tilde\F^\ell:=(\tilde\Fc_s^{B,\mu_\ell}\vee\tilde\Gc)_{s\geq0}$ for every $\ell\geq1$, $\tilde\F^0:=(\tilde\Fc_s^B\vee\tilde\Gc)_{s\geq0}$, $\tilde I^0:=\tilde\alpha$, $\tilde X^{t,\tilde\xi,0}:=\tilde X^{t,\tilde\xi,\tilde\alpha}$, and $\tilde X^{t,x,\pi,0}:=\tilde X^{t,x,\pi,\tilde\alpha}$. 

We conclude that for every $\eps>0$ there exists some $L_\eps\in\N$ such that, for every $\ell>L_\eps$, we have
\[
J(t,x,\pi,\hat\alpha) - \eps \ \leq \ \tilde J_\ell^\Rc(t,x,\pi,\tilde\nu^\ell) \ \leq \ \sup_{\tilde\nu\in\tilde\Vc_\ell}\tilde J_\ell^\Rc(t,x,\pi,\tilde\nu) \ =: \ \tilde V_\ell^\Rc(t,x,\pi) \ \overset{\overset{\text{Step I}}{\downarrow}}{=} \ V^\Rc(t,x,\pi).
\]
From the arbitrariness of $\eps$, we see that $J(t,x,\pi,\hat\alpha)\leq V^\Rc(t,x,\pi)$. The claim follows taking the supremum over $\hat\alpha\in\Ac$.
\ep

\begin{Remark}\label{R:Vc1,t}
{\rm
Let $\Vc_{1,t}\subset\Vc$ be the set of $\nu\in\Vc$ such that $\nu\equiv1$ on $\Omega\times[0,t)\times\Ac$. Then
\begin{equation}\label{equivalence_1,t}
V(t,x,\pi) \ = \ \sup_{\nu\in\Vc_{1,t}} J^\Rc(t,x,\pi,\nu),
\end{equation}
for all $(t,x,\pi)\in[0,T]\times\R^n\times\mathscr P_{\text{\tiny$2$}}(\R^n)$. Indeed, by step II of the proof of Theorem \ref{Thm:Equivalence}, we have $V(t,x,\pi)\geq V^\Rc(t,x,\pi)\geq\sup_{\nu\in\Vc_{1,t}} J^\Rc(t,x,\pi,\nu)$. Let us prove the other inequality. We begin noting that in Lemma \ref{L:Stability}, the convergence $\E^\Q[\int_t^T\tilde\rho(\tilde I_s^\ell,\tilde I_s^0)\,ds]\rightarrow0$ as $\ell\rightarrow\infty$ is needed, rather than $\E^\Q[\int_0^T\tilde\rho(\tilde I_s^\ell,\tilde I_s^0)\,ds]\rightarrow0$. In other words, the behavior of $(\tilde I_s^\ell)_{s\in[0,T]}$ on the interval $[0,t)$ is not relevant. Therefore, proceeding as in step III of the proof of Theorem \ref{Thm:Equivalence}, we see that we can take $\tilde\nu^\ell\equiv1$ on $\tilde\Omega\times[0,t)\times\Ac$, in order to guarantee the convergence $\E^\Q[\int_t^T\tilde\rho(\tilde I_s^\ell,\hat\alpha_s)\,ds]\rightarrow0$ as $\ell\rightarrow\infty$. Then, from the same proof of Lemma \ref{L:Stability}, we conclude that $\tilde J^\Rc(t,x,\pi,\tilde\nu^\ell)\rightarrow J(t,x,\pi,\hat\alpha)$ as $\ell\rightarrow\infty$. This implies the validity of the other inequality $V(t,x,\pi)\leq\sup_{\nu\in\Vc_{1,t}}J^\Rc(t,x,\pi)$ and proves \eqref{equivalence_1,t}.
\ep
}
\end{Remark}

\section{Feynman-Kac representation: randomized equation}
\label{S:Feynman-Kac}

In the present section we introduce, for every $(t,x,\bar\xi)\in[0,T]\times\R^n\times L^2(\bar\Omega,\bar\Gc,\bar\P;\R^n)$, a forward-backward stochastic differential system of equations, which provides a probabilistic representation for the value $V(t,x,\pi)$, with $\pi=\P_{\text{\tiny$\xi$}}$ under $\bar\P$. In other words, we derive a nonlinear Feynman-Kac formula for the value function $V$ in \eqref{Value} of the McKean-Vlasov control problem.

We firstly introduce the following spaces, for every $t\in[0,T]$.
\begin{itemize}
\item $\Sc^2(t,T)$, the set of real-valued c\`adl\`ag $\F^\mu$-adapted processes $Y$ $=$ $(Y_s)_{s\in[t,T]}$, with $Y\colon\Omega^1\times[t,T]\rightarrow\R$, satisfying $\|Y\|_{\Sc^2(t,T)}^2$ $:=$ $\E^1\big[\sup_{t\leq s\leq T}|Y_s|^2\big]$ $<$ $\infty$.
\item $L_{\tilde\mu}^2(t,T)$, the set of real-valued $\Pc(\F^\mu)\otimes\Bc(\Ac)$-measurable maps $U=(U_s(\alpha))_{s\in[t,T],\,\alpha\in\Ac}$, with $U\colon\Omega^1\times[t,T]\times\Ac\rightarrow\R$, satisfying $\|U\|_{L_{\tilde\mu}^2(t,T)}^2$ $:=$ $\E^1\big[\int_t^T\int_\Ac|U_s(\alpha)|^2\lambda(d\alpha)\,ds\big]$ $<$ $\infty$.
  \item $\Kc^2(t,T)$, the set of nondecreasing $\F^\mu$-predictable processes $K=(K_s)_{s\in[t,T]}$, with $K\colon\Omega^1\times[t,T]\rightarrow\R_+$, satisfying $K\in\Sc^2(t,T)$ and $K_t=0$.
\end{itemize}

Given $(t,x,\bar\xi)\in[0,T]\times\R^n\times L^2(\bar\Omega,\bar\Gc,\bar\P;\R^n)$, with $\pi=\P_{\text{\tiny$\xi$}}$ under $\bar\P$, consider on $(\Omega^1,\Fc^1,\F^\mu,\P^1)$ the following backward stochastic differential equation with constrained jumps over $[t,T]$:
\begin{equation}\label{BSDE}
\begin{cases}
\vspace{2mm} \dis Y_s \ = \ \E\big[g(\bar X_T^{t,x,\pi},\P_{\text{\tiny$\bar X_T^{t,\bar\xi}$}}^{\text{\tiny$\bar\Fc_T^\mu$}})\big]  + \int_s^T \E\big[f(r,\bar X_r^{t,x,\pi},\P_{\text{\tiny$\bar X_r^{t,\bar\xi}$}}^{\text{\tiny$\bar\Fc_r^\mu$}},\bar I_r)\big]\,dr + K_T - K_s \\
\vspace{2mm}
\qquad\;\;\;\dis - \int_s^T\int_A U_r(\alpha)\,\mu(dr\,d\alpha), \qquad\qquad\, s\in[t,T], \\
\dis U_s(\alpha) \ \leq \ 0, \qquad\qquad\qquad\quad\;\, d\P^1\,ds\,\lambda(d\alpha)\text{-a.e. on }\Omega^1\times[t,T]\times\Ac.
\end{cases}
\end{equation}
Notice that $\E\big[g(\bar X_T^{t,x,\pi},\P_{\text{\tiny$\bar X_T^{t,\bar\xi}$}}^{\text{\tiny$\bar\Fc_T^\mu$}})\big]$, as well as $\E\big[f(r,\bar X_r^{t,x,\pi},\P_{\text{\tiny$\bar X_r^{t,\bar\xi}$}}^{\text{\tiny$\bar\Fc_r^\mu$}},\bar I_r)\big]$, is a random variable on $(\Omega^1,\Fc^1,\P^1)$.

Equations \eqref{barI}-\eqref{StateEq1Rand}-\eqref{StateEq2Rand}-\eqref{BSDE} constitute a forward-backward stochastic differential system of equations. We also observe that equation \eqref{BSDE} depends on $\xi$ only through its law $\pi=\P_{\text{\tiny$\xi$}}$. We now prove that there exists a unique solution $(Y^{t,x,\pi},U^{t,x,\pi},K^{t,x,\pi})\in\Sc^2(t,T)\times L_{\tilde\mu}^2(t,T)\times\Kc^2(t,T)$ to \eqref{BSDE}, which is minimal in the following sense: if $(\bar Y,\bar U,\bar K)\in\Sc^2(t,T)\times L_{\tilde\mu}^2(t,T)\times\Kc^2(t,T)$ is another solution to \eqref{BSDE}, then the inequality $Y^{t,x,\pi}\leq\bar Y$ holds on $\Omega^1\times[t,T]$, up to a $\P^1$-evanescent set.

\begin{Theorem}\label{Thm:Feynman-Kac}
Under Assumption {\bf (A1)}, for every $(t,x,\bar\xi)\in[0,T]\times\R^n\times L^2(\bar\Omega,\bar\Gc,\bar\P;\R^n)$, with $\pi=\P_{\text{\tiny$\xi$}}$ under $\bar\P$, there exists a unique minimal solution $(Y^{t,x,\pi},U^{t,x,\pi},K^{t,x,\pi})\in\Sc^2(t,T)\times L_{\tilde\mu}^2(t,T)\times\Kc^2(t,T)$ to \eqref{BSDE}, with $Y_t^{t,x,\pi}$ equal $\P^1$-a.s. to a constant. In addition, $V$ admits the Feynman-Kac representation
\begin{equation}\label{Feynman-Kac}
V(t,x,\pi) \ = \ Y_t^{t,x,\pi}
\end{equation}
$\P^1$-a.s., for all $(t,x,\pi)\in[0,T]\times\R^n\times\mathscr P_{\text{\tiny$2$}}(\R^n)$. Moreover, we have
\begin{align}\label{DPPBSDE}
Y_t^{t,x,\pi} \ &= \ \sup_{\nu\in\Vc} \E^\nu\bigg[\int_t^s \E\big[f(r,\bar X_r^{t,x,\pi},\P_{\text{\tiny$\bar X_r^{t,\bar\xi}$}}^{\text{\tiny$\bar\Fc_r^\mu$}},\bar I_r)\big]\,dr + Y_s^{t,x,\pi}\bigg] \\
&= \ \sup_{\nu\in\Vc} \bar\E^\nu\bigg[\int_t^s f(r,\bar X_r^{t,x,\pi},\P_{\text{\tiny$\bar X_r^{t,\bar\xi}$}}^{\text{\tiny$\bar\Fc_r^\mu$}},\bar I_r)\,dr + Y_s^{t,x,\pi}\bigg], \notag
\end{align}
$\P^1$-a.s., for all $s\in[t,T]$.
\end{Theorem}
\textbf{Proof.}
\emph{Existence and uniqueness of the minimal solution to \eqref{BSDE}.} Fix $(t,x,\bar\xi)\in[0,T]\times\R^n\times L^2(\bar\Omega,\bar\Gc,\bar\P;\R^n)$, with $\pi=\P_{\text{\tiny$\xi$}}$ under $\bar\P$. Consider, for every $n\in\N$, the following unconstrained backward stochastic differential equation on $[t,T]$:
\begin{align}\label{BSDE_n}
Y_s \ &= \ \E\big[g(\bar X_T^{t,x,\pi},\P_{\text{\tiny$\bar X_T^{t,\bar\xi}$}}^{\text{\tiny$\bar\Fc_T^\mu$}})\big]  + \int_s^T \E\big[f(r,\bar X_r^{t,x,\pi},\P_{\text{\tiny$\bar X_r^{t,\bar\xi}$}}^{\text{\tiny$\bar\Fc_r^\mu$}},\bar I_r)\big]\,dr + n\int_s^T\int_\Ac (U_r(\alpha))_+\,\lambda(d\alpha)\,dr \notag \\
&\quad \ - \int_s^T\int_\Ac U_r(\alpha)\,\mu(dr\,d\alpha).
\end{align}
By Lemma 2.4 in \cite{tang_li}, there exists a unique solution $(Y^{n,t,x,\pi},U^{n,t,x,\pi})\in\Sc^2(t,T)\times L_{\tilde\mu}^2(t,T)$ to the above equation. 

For every $n\in\N$, let $\hat\Vc^n$ denote the set of $\Pc(\F^\mu)\otimes\Bc(\Ac)$-measurable maps $\hat\nu\colon\Omega^1\times\R_+\times\Ac\rightarrow(0,n]$, which are not necessarily bounded away from zero. Then, let us prove the following formula: 
\begin{equation}\label{dual_n}
Y_{\bar t}^{n,t,x,\pi} \ = \ \esssup_{\hat\nu\in\hat\Vc^n} \E^{\hat\nu}\bigg[\int_{\bar t}^s \E\big[f(r,\bar X_r^{t,x,\pi},\P_{\text{\tiny$\bar X_r^{t,\bar\xi}$}}^{\text{\tiny$\bar\Fc_r^\mu$}},\bar I_r)\big]\,dr + Y_s^{n,t,x,\pi}\bigg|\Fc_{\bar t}^\mu\bigg],
\end{equation}
for all $\bar t,s\in[t,T]$, with $\bar t\leq s$. Let $\hat\nu\in\hat\Vc$ (see Remark \ref{R:hatV} for the definition of $\hat\Vc$). Then, considering \eqref{BSDE_n} between $\bar t$ and $s$, and taking the $\P^{\hat\nu}$-conditional expectation with respect to $\Fc_{\bar t}^\mu$, we obtain
\begin{align}\label{Proof_dual_n}
Y_{\bar t}^{n,t,x,\pi} \ &= \ \E^{\hat\nu}\bigg[\int_{\bar t}^s \E\big[f(r,\bar X_r^{t,x,\pi},\P_{\text{\tiny$\bar X_r^{t,\bar\xi}$}}^{\text{\tiny$\bar\Fc_r^\mu$}},\bar I_r)\big]\,dr + Y_s^{n,t,x,\pi} \\
&\quad \ + \int_{\bar t}^s\int_\Ac\big[n(U_r(\alpha)^{n,t,x,\pi})_+ - U_r^{n,t,x,\pi}(\alpha)\hat\nu_r(\alpha)\big]\,\lambda(d\alpha)\,dr\bigg|\Fc_{\bar t}^\mu\bigg]. \notag
\end{align}
Since $\hat\nu_r(\alpha)\in(0,n]$, the last term inside the expectation is nonnegative. Therefore
\begin{equation}\label{FirstIneqYn}
Y_{\bar t}^{n,t,x,\pi} \ \geq \ \esssup_{\hat\nu\in\hat\Vc^n} \E^{\hat\nu}\bigg[\int_{\bar t}^s \E\big[f(r,\bar X_r^{t,x,\pi},\P_{\text{\tiny$\bar X_r^{t,\bar\xi}$}}^{\text{\tiny$\bar\Fc_r^\mu$}},\bar I_r)\big]\,dr + Y_s^{n,t,x,\pi}\bigg|\Fc_{\bar t}^\mu\bigg].
\end{equation}
To prove the other inequality, define, for every $\eps\in(0,n]$, the map $\hat\nu^{n,\eps}$ as
\[
\hat\nu_r^{n,\eps}(\alpha) \ = \ n\,1_{\{U_r^{n,t,x,\pi}(\alpha)\geq0\}} + \eps\,1_{\{-1\leq U_r^{n,t,x,\pi}(\alpha)<0\}} + \frac{\eps}{|U_r^{n,t,x,\pi}(\alpha)|}\,1_{\{U_r^{n,t,x,\pi}(\alpha)<-1\}},
\]
on $\Omega^1\times[t,T]\times\Ac$, and $\hat\nu^{n,\eps}\equiv1$ on $\Omega^1\times([0,t)\cup(T,\infty))\times\Ac$. Notice that $\hat\nu^{n,\eps}$ belongs to $\hat\Vc^n$, and it is not necessarily bounded away from zero. Taking $\hat\nu$ equal to $\hat\nu^{n,\eps}$ in \eqref{Proof_dual_n}, we obtain
\begin{align}\label{Yn_nu^esp}
Y_{\bar t}^{n,t,x,\pi} \ &\leq \ \E^{\hat\nu^{n,\eps}}\bigg[\int_{\bar t}^s \E\big[f(r,\bar X_r^{t,x,\pi},\P_{\text{\tiny$\bar X_r^{t,\bar\xi}$}}^{\text{\tiny$\bar\Fc_r^\mu$}},\bar I_r)\big]\,dr + Y_s^{n,t,x,\pi}\bigg|\Fc_{\bar t}^\mu\bigg] + \eps(T-t)\lambda(\Ac) \\
&\leq \ \esssup_{\hat\nu\in\hat\Vc^n} \bigg[\int_{\bar t}^s \E\big[f(r,\bar X_r^{t,x,\pi},\P_{\text{\tiny$\bar X_r^{t,\bar\xi}$}}^{\text{\tiny$\bar\Fc_r^\mu$}},\bar I_r)\big]\,dr + Y_s^{n,t,x,\pi}\bigg|\Fc_{\bar t}^\mu\bigg] + \eps(T-t)\lambda(\Ac). \notag
\end{align}
From the arbitrariness of $\eps$ we get the reverse inequality of \eqref{FirstIneqYn}, from which we deduce the validity of \eqref{dual_n}. In particular, when $s=T$ in \eqref{dual_n}, we obtain
\begin{equation}\label{dual_n_t}
Y_{\bar t}^{n,t,x,\pi} \ = \ \esssup_{\hat\nu\in\hat\Vc^n} \E^{\hat\nu}\bigg[\int_{\bar t}^T \E\big[f(r,\bar X_r^{t,x,\pi},\P_{\text{\tiny$\bar X_r^{t,\bar\xi}$}}^{\text{\tiny$\bar\Fc_r^\mu$}},\bar I_r)\big]\,dr + \E\big[g(\bar X_T^{t,x,\pi},\P_{\text{\tiny$\bar X_T^{t,\bar\xi}$}}^{\text{\tiny$\bar\Fc_T^\mu$}})\big]\bigg|\Fc_{\bar t}^\mu\bigg],
\end{equation}
for all $\bar t\in[t,T]$. Then, it is easy to see that the following estimate holds:
\begin{equation}\label{Yn<infty}
\sup_n Y_{\bar t}^{n,t,x,\pi} \ < \ \infty, \qquad \text{for all }\bar t\in[t,T].
\end{equation}
Hence, the existence and uniqueness of the minimal solution to equation \eqref{BSDE} follows from Theorem 2.1 in \cite{KP15} (apart from the fact that $K_t^{t,x,\pi}=0$, as required in the definition of $\Kc^2(t,T)$, which will be proved later). Indeed, \eqref{BSDE} can be seen as an equation on the entire interval $[0,T]$, with terminal condition $\E\big[g(\bar X_T^{t,x,\pi},\P_{\text{\tiny$\bar X_T^{t,\bar\xi}$}}^{\text{\tiny$\bar\Fc_T^\mu$}})\big]$ and generator $\E\big[f(r,\bar X_r^{t,x,\pi},\P_{\text{\tiny$\bar X_r^{t,\bar\xi}$}}^{\text{\tiny$\bar\Fc_r^\mu$}},\bar I_r)\big]1_{[t,T]}(r)$. Assumption (H0) in \cite{KP15} holds under Assumption {\bf (A1)}. Moreover, Assumption (H1) in \cite{KP15} is imposed only to guarantee the validity of \eqref{Yn<infty}, which in our case follows directly from formula \eqref{dual_n_t}, since $f$ does not depend on $Y^{n,t,x,\pi}$, $U^{n,t,x,\pi}$. It only remains to prove that $K_t^{t,x,\pi}=0$. This is clearly true if we show that $Y_t^{t,x,\pi}$ is equal $\P^1$-a.s. to a constant (as a matter of fact, if $Y_t^{t,x,\pi}$ is equal $\P^1$-a.s. to a constant, then, by uniqueness, $Y_s^{t,x,\pi}=Y_t^{t,x,\pi}$ on $[0,t]$, so that $K_s^{t,x,\pi}$ is also constant on $[0,t]$, and, in particular, equal to $K_0^{t,x,\pi}=0$). This latter property is proved below. Finally, for later use, we notice that, according to Theorem 2.1 in \cite{KP15}, the sequence $(Y_{\bar t}^{n,t,x,\pi})_{n\geq0}$ is nondecreasing (this is a direct consequence of formula \eqref{dual_n_t}, since $\hat\Vc^n\subset\hat\Vc^{n+1}$) and converges pointwise $\P^1$-a.s. to $Y_{\bar t}^{t,x,\pi}$, for all $\bar t\in[t,T]$.

\vspace{2mm}

\noindent\emph{Proof of \eqref{Feynman-Kac}, in particular $Y_t^{t,x,\pi}$ is equal $\P^1$-a.s. to a constant.} Notice that $Y_t^{t,x,\pi}$ is $\Fc_t^\mu$-measurable, therefore it is not a priori clear that it is $\P^1$-a.s. a constant. For every $n\in\N$, consider \eqref{dual_n} with $\bar t=t$ and $s=T$:
\[
Y_t^{n,t,x,\pi} \ = \ \esssup_{\hat\nu\in\hat\Vc^n} \E^{\hat\nu}\bigg[\int_t^T \E\big[f(r,\bar X_r^{t,x,\pi},\P_{\text{\tiny$\bar X_r^{t,\bar\xi}$}}^{\text{\tiny$\bar\Fc_r^\mu$}},\bar I_r)\big]\,dr + \E\big[g(\bar X_T^{t,x,\pi},\P_{\text{\tiny$\bar X_T^{t,\bar\xi}$}}^{\text{\tiny$\bar\Fc_T^\mu$}})\big]\bigg|\Fc_t^\mu\bigg].
\]
Letting $n\rightarrow\infty$, recalling that $Y_t^{n,t,x,\pi}\nearrow Y_t^{t,x,\pi}$, $\P^1$-a.s., and noting that $\hat\Vc_n\subset\hat\Vc_{n+1}\subset\cup_n\hat\Vc_n=\hat\Vc$, we obtain
\begin{equation}\label{Ydual}
Y_t^{t,x,\pi} \ = \ \esssup_{\hat\nu\in\hat\Vc} \E^{\hat\nu}\bigg[\int_t^T \E\big[f(r,\bar X_r^{t,x,\pi},\P_{\text{\tiny$\bar X_r^{t,\bar\xi}$}}^{\text{\tiny$\bar\Fc_r^\mu$}},\bar I_r)\big]\,dr + \E\big[g(\bar X_T^{t,x,\pi},\P_{\text{\tiny$\bar X_T^{t,\bar\xi}$}}^{\text{\tiny$\bar\Fc_T^\mu$}})\big]\bigg|\Fc_t^\mu\bigg].
\end{equation}
Reasoning as in Remark \ref{R:hatV}, we can show that the right-hand side of \eqref{Ydual} does not change if we take the supremum over $\Vc$. In other words, \eqref{Ydual} can be equivalently written as follows:
\begin{equation}\label{YesssupVc}
Y_t^{t,x,\pi} \ = \ \esssup_{\nu\in\Vc} \E^\nu\bigg[\int_t^T \E\big[f(r,\bar X_r^{t,x,\pi},\P_{\text{\tiny$\bar X_r^{t,\bar\xi}$}}^{\text{\tiny$\bar\Fc_r^\mu$}},\bar I_r)\big]\,dr + \E\big[g(\bar X_T^{t,x,\pi},\P_{\text{\tiny$\bar X_T^{t,\bar\xi}$}}^{\text{\tiny$\bar\Fc_T^\mu$}})\big]\bigg|\Fc_t^\mu\bigg].
\end{equation}
From Corollary \ref{C:EquivConditional} it follows that the right-hand side of \eqref{YesssupVc} is equal $\P^1$-a.s. to $V(t,x,\pi)$, which yields $Y_t^{t,x,\pi}=V(t,x,\pi)$, $\P^1$-a.s..

\vspace{2mm}

\noindent\emph{Proof of formula \eqref{DPPBSDE}.} Let $\nu\in\Vc$. Consider \eqref{BSDE} between $t$ and $s$, and take the expectation with respect to $\E^\nu$, then (recalling that $K^{t,x,\pi}$ is nondecreasing and $U^{t,x,\pi}$ is nonpositive)
\begin{align}\label{Proof_dual}
Y_t^{t,x,\pi} \ &\geq \ \E^\nu\bigg[\int_t^s \E\big[f(r,\bar X_r^{t,x,\pi},\P_{\text{\tiny$\bar X_r^{t,\bar\xi}$}}^{\text{\tiny$\bar\Fc_r^\mu$}},\bar I_r)\big]\,dr + Y_s^{t,x,\pi}\bigg].
\end{align}
From the arbitrariness of $\nu\in\Vc$, we get the first inequality. To prove the reverse inequality, considering \eqref{Yn_nu^esp} with $\bar t=t$, and taking the expectation $\E^{\hat\nu^{n,\eps}}$, we obtain
\begin{align*}
\E^{\hat\nu^{n,\eps}}\big[Y_t^{n,t,x,\pi}\big] \ &\leq \ \E^{\hat\nu^{n,\eps}}\bigg[\int_t^s \E\big[f(r,\bar X_r^{t,x,\pi},\P_{\text{\tiny$\bar X_r^{t,\bar\xi}$}}^{\text{\tiny$\bar\Fc_r^\mu$}},\bar I_r)\big]\,dr + Y_s^{n,t,x,\pi}\bigg] + \eps(T-t)\lambda(\Ac) \\
&\leq \ \sup_{\hat\nu\in\hat\Vc} \E^{\hat\nu}\bigg[\int_t^s \E\big[f(r,\bar X_r^{t,x,\pi},\P_{\text{\tiny$\bar X_r^{t,\bar\xi}$}}^{\text{\tiny$\bar\Fc_r^\mu$}},\bar I_r)\big]\,dr + Y_s^{t,x,\pi}\bigg] + \eps(T-t)\lambda(\Ac) \\
&= \ \sup_{\nu\in\Vc} \E^\nu\bigg[\int_t^s \E\big[f(r,\bar X_r^{t,x,\pi},\P_{\text{\tiny$\bar X_r^{t,\bar\xi}$}}^{\text{\tiny$\bar\Fc_r^\mu$}},\bar I_r)\big]\,dr + Y_s^{t,x,\pi}\bigg] + \eps(T-t)\lambda(\Ac), 
\end{align*}
where the last equality can be proved arguing as in Remark \ref{R:hatV}. From the definition of $\hat\nu^{n,\eps}$, we see that $\kappa_t^{\hat\nu^{n,\eps}}=1$, therefore $\E^{\hat\nu^{n,\eps}}[Y_t^{n,t,x,\pi}]=\E^1[Y_t^{n,t,x,\pi}]$. Hence
\[
\E^1\big[Y_t^{n,t,x,\pi}\big] \ \leq \ \sup_{\nu\in\Vc} \E^\nu\bigg[\int_t^s \E\big[f(r,\bar X_r^{t,x,\pi},\P_{\text{\tiny$\bar X_r^{t,\bar\xi}$}}^{\text{\tiny$\bar\Fc_r^\mu$}},\bar I_r)\big]\,dr + Y_s^{t,x,\pi}\bigg] + \eps(T-t)\lambda(\Ac).
\]
Recall that the sequence $(Y_t^{n,t,x,\pi})_{n\geq0}$ is nondecreasing and converges pointwise $\P^1$-a.s. to $Y_t^{t,x,\pi}$. In particular, $Y_t^{0,t,x,\pi}\leq Y_t^{n,t,x,\pi}\leq Y_t^{t,x,\pi}$, for every $n\in\N$. Therefore, letting $n\rightarrow\infty$ and using Lebesgue's dominated convergence theorem, we obtain
\[
Y_t^{t,x,\pi} \ = \ \E^1\big[Y_t^{t,x,\pi}\big] \ \leq \ \sup_{\nu\in\Vc} \E^\nu\bigg[\int_t^s \E\big[f(r,\bar X_r^{t,x,\pi},\P_{\text{\tiny$\bar X_r^{t,\bar\xi}$}}^{\text{\tiny$\bar\Fc_r^\mu$}},\bar I_r)\big]\,dr + Y_s^{t,x,\pi}\bigg] + \eps(T-t)\lambda(\Ac).
\]
Sending $\eps\rightarrow0$, we get
\[
Y_t^{t,x,\pi} \ \leq \ \sup_{\nu\in\Vc} \E^\nu\bigg[\int_t^s \E\big[f(r,\bar X_r^{t,x,\pi},\P_{\text{\tiny$\bar X_r^{t,\bar\xi}$}}^{\text{\tiny$\bar\Fc_r^\mu$}},\bar I_r)\big]\,dr + Y_s^{t,x,\pi}\bigg],
\]
which, together with \eqref{Proof_dual}, gives formula \eqref{DPPBSDE} and concludes the proof.
\ep

\section{Randomized dynamic programming principle}
\label{S:DPP}

The present section is devoted to the proof of the dynamic programming principle for $V$ in the randomized framework. Firstly, we prove the flow properties of $\bar X^{t,\bar\xi}$ and $\bar X^{t,x,\pi}$. These in turn imply the identification $\E\big[V(s,\bar X_s^{t,x,\pi},\P_{\text{\tiny$\bar X_s^{t,\bar\xi}$}}^{\text{\tiny$\bar\Fc_s^\mu$}})\big]=Y_s^{t,x,\pi}$, $\P^1$ -a.s., for all $s\in[t,T]$. Then, \eqref{DPPBSDE} allows to derive the randomized dynamic programming principle for $V$.

\subsection{Flow properties}
\label{SubS:Flow}

We begin considering the solution to system \eqref{StateEq1Rand}-\eqref{StateEq2Rand} with more general initial conditions. More precisely, concerning equation \eqref{StateEq1Rand}, for every $(t,\bar\eta)\in[0,T]\times L^2(\bar\Omega,\bar\Fc_t^{B,\mu}\vee\bar\Gc,\bar\P;\R^n)$, consider the following equation:
\begin{equation}
d\bar X_s^{t,\bar\eta} \ = \ b\big(s,\bar X_s^{t,\bar\eta},\P_{\text{\tiny$\bar X_s^{t,\bar\eta}$}}^{\text{\tiny$\bar\Fc_s^\mu$}},\bar I_s\big)\,ds + \sigma\big(s,\bar X_s^{t,\bar\eta},\P_{\text{\tiny$\bar X_s^{t,\bar\eta}$}}^{\text{\tiny$\bar\Fc_s^\mu$}},\bar I_s\big)\,d\bar B_s, \qquad \bar X_t^{t,\bar\eta} \ = \ \bar\eta,\label{StateEq1Rand_general}
\end{equation}
for all $s\in[t,T]$. Concerning equation \eqref{StateEq2Rand}, we begin recalling that $(\P_{\text{\tiny$\bar X_s^{t,\bar\eta}$}}^{\text{\tiny$\bar\Fc_s^\mu$}})_{s\in[t,T]}$ stands for the stochastic process $(\P_s^{t,\pi})_{s\in[t,T]}$ introduced in Lemma \ref{L:P}, with $\pi=\P_{\text{\tiny$\bar\eta$}}$ under $\bar\P$. In the sequel, when considering equation \eqref{StateEq2Rand}, it is more convenient to adopt the notation $\P_s^{t,\pi}$ instead of $\P_{\text{\tiny$\bar X_s^{t,\bar\eta}$}}^{\text{\tiny$\bar\Fc_s^\mu$}}$. For every $(t,\bar\eta)\in[0,T]\times L^2(\bar\Omega,\bar\Fc_t^{B,\mu}\vee\bar\Gc,\bar\P;\R^n)$ and $\bar\Pi\colon\bar\Omega\rightarrow\mathscr P_{\text{\tiny$2$}}(\R^n)$, with $\bar\Pi$ measurable with respect to $\bar\Fc_t^\mu$ and such that $\bar\E[\|\bar\Pi\|_{\text{\tiny$2$}}^2]<\infty$, consider the following equation:
\begin{equation}
d\bar X_s^{t,\bar\eta,\bar\Pi} \ = \ b\big(s,\bar X_s^{t,\bar\eta,\bar\Pi},\P_s^{t,\bar\Pi},\bar I_s\big)\,ds + \sigma\big(s,\bar X_s^{t,\bar\eta,\bar\Pi},\P_s^{t,\bar\Pi},\bar I_s\big)\,d\bar B_s, \qquad \bar X_t^{t,\bar\eta,\bar\Pi} \ = \ \bar\eta, \label{StateEq2Rand_general}
\end{equation}
for all $s\in[t,T]$, where
\begin{equation}\label{P^Pi}
\P_s^{t,\bar\Pi}(\bar\omega) \ := \ \P_s^{t,\bar\Pi(\bar\omega)}(\omega^1), \qquad \text{for all }(\bar\omega,s)=(\omega,\omega^1,s)\in\bar\Omega\times[t,T].
\end{equation}
Notice that, thanks to Lemma \ref{L:P}, the stochastic process $(\P_s^{t,\bar\Pi})_{s\in[t,T]}$ is well-defined. In particular, for every $s\in[t,T]$, $\P_s^{t,\bar\Pi}$ is $\bar\Fc_s^\mu$-measurable. Under Assumption {\bf (A1)}, we have the following result, whose standard proof is not reported.

\begin{Lemma}\label{L:ExistUniqX}
Under Assumption {\bf (A1)}, for every $(t,\bar\eta)\in[0,T]\times L^2(\bar\Omega,\bar\Fc_t^{B,\mu}\vee\bar\Gc,\bar\P;\R^n)$ and $\bar\Pi\colon\bar\Omega\rightarrow\mathscr P_{\text{\tiny$2$}}(\R^n)$, with $\bar\Pi$ measurable with respect to $\bar\Fc_t^\mu$ and such that $\bar\E[\|\bar\Pi\|_{\text{\tiny$2$}}^2]<\infty$, there exists a unique (up to indistinguishability) pair $(\bar X_s^{t,\bar\eta},\bar X_s^{t,\bar\eta,\bar\Pi})_{s\in[t,T]}$ of continuous $(\bar\Fc_s^{B,\mu,t}\vee\bar\Gc\vee\sigma(\bar\eta,\bar\Pi))_s$-adapted processes solution to equations \eqref{StateEq1Rand_general}-\eqref{StateEq2Rand_general}, satisfying
\[
\bar\E\Big[\sup_{s\in[t,T]}\big(\big|\bar X_s^{t,\bar\eta}\big|^2 + \big|\bar X_s^{t,\bar\eta,\bar\Pi}\big|^2\big)\Big] \ < \ \infty.
\]
Moreover, there exists a positive constant $C$ such that
\begin{equation}\label{EstimateMKV}
\bar\E\Big[\sup_{s\in[t,T]}\big|\bar X_s^{t,\bar\eta,\bar\Pi} - \bar X_s^{t,\bar\eta',\bar\Pi'}\big|^2\Big] \ \leq \ C\big(\bar\E[|\bar\eta - \bar\eta'|^2] + \bar\E[\Wc_{\text{\tiny$2$}}(\bar\Pi,\bar\Pi')^2]\big),
\end{equation}
for every $t\in[0,T]$, $\bar\eta,\bar\eta'\in L^2(\bar\Omega,\bar\Fc_t^{B,\mu}\vee\bar\Gc,\bar\P;\R^n)$, and any $\bar\Pi,\bar\Pi'\colon\bar\Omega\rightarrow\mathscr P_{\text{\tiny$2$}}(\R^n)$, with $\bar\Pi,\bar\Pi'$ measurable with respect to $\bar\Fc_t^\mu$ and such that $\bar\E[\|\bar\Pi\|_{\text{\tiny$2$}}^2],\bar\E[\|\bar\Pi'\|_{\text{\tiny$2$}}^2]<\infty$.
\end{Lemma}
\textbf{Proof.}
The proof of the existence and uniqueness of $(\bar X_s^{t,\bar\eta},\bar X_s^{t,\bar\eta,\bar\Pi})_{s\in[t,T]}$ is standard under Assumption {\bf (A1)}, and can be done as usual by a fixed point argument. Concerning estimate \eqref{EstimateMKV}, the proof can be done proceeding as in Lemma 3.1 in \cite{BLPR14}.
\ep

\begin{Remark}
{\rm
When in equation \eqref{StateEq2Rand_general} the random variables $\bar\eta$ and $\bar\Pi$ are equal $\bar\P$-a.s. to some $x\in\R^n$ and $\pi\in\mathscr P_{\text{\tiny$2$}}(\R^n)$, respectively, then $(\bar X_s^{t,\bar\eta,\bar\Pi})_{s\in[t,T]}$ coincides (up to indistinguishability) with the stochastic process $(\bar X_s^{t,x,\pi})_{s\in[t,T]}$ defined in Section \ref{S:RandomizedProblem}. Indeed, $(\bar X_s^{t,\bar\eta,\bar\Pi})_{s\in[t,T]}$ and $(\bar X_s^{t,x,\pi})_{s\in[t,T]}$ solve the same equation, therefore the claim follows from the uniqueness of the solution.
\ep}
\end{Remark}

\begin{Remark}\label{R:X_Partition}
{\rm
Suppose that $\bar\eta$ and $\bar\Pi$ in Lemma \ref{L:ExistUniqX} takes only a finite number of values, namely
\[
\bar\eta \ = \ \sum_{k=0}^K x_k\,1_{E_k}, \qquad\qquad\quad \bar\Pi \ = \ \sum_{k=0}^K \pi_k\,1_{E_k},
\]
for some $K\in\N$, $x_k\in\R^n$, $\pi_k\in\mathscr P_{\text{\tiny$2$}}(\R^n)$, $E_k\in\bar\Fc_t^{B,\mu}\vee\bar\Gc$, with $(E_k)_{k=1,\ldots,K}$ being a partition of $\bar\Omega$. Then, by definition of $\P_s^{t,\bar\Pi}$ (formula \eqref{P^Pi}), we have $\P_s^{t,\bar\Pi}=\P_s^{t,\pi_0}\,1_{E_0}+\cdots+\P_s^{t,\pi_K}\,1_{E_K}$. Therefore, the stochastic processes $(\bar X_s^{t,x_0,\pi_0}\,1_{E_0}+\cdots+\bar X_s^{t,x_K,\pi_K}\,1_{E_K})_{s\in[t,T]}$ and $(\bar X_s^{t,\bar\eta,\bar\Pi})_{s\in[t,T]}$ are indistinguishable, since they solve the same stochastic differential equation.
\ep}
\end{Remark}

\begin{Lemma}\label{L:Flow}
Under Assumption {\bf (A1)}, for every $(t,s,x,\bar\xi)\in[0,T]\times[0,T]\times\R^n\times L^2(\bar\Omega,\bar\Gc,\bar\P;\R^n)$, with $t\leq s$ and $\pi=\P_{\text{\tiny$\bar\xi$}}$ under $\bar\P$, we have the flow properties:
\begin{align}
\bar X_r^{s,\bar X_s^{t,\bar\xi}} \ = \ \bar X_r^{t,\bar\xi}, \label{FlowX1} \\
\bar X_r^{s,\bar X_s^{t,x,\pi},\P_{\text{\tiny$\bar X_s^{t,\bar\xi}$}}^{\text{\tiny$\bar\Fc_s^\mu$}}} \ = \ \bar X_r^{t,x,\pi}, \label{FlowX2}
\end{align}
for all $r\in[s,T]$, $\bar\P$-almost surely.
\end{Lemma}
\textbf{Proof.} \emph{Flow property \eqref{FlowX1}.} Consider the process $(\bar X_r^{s,\bar X_s^{t,\bar\xi}})_{r\in[s,T]}$ solution to equation \eqref{StateEq1Rand_general} with initial conditions $t=s$ and $\bar\eta=\bar X_s^{t,\bar\xi}$. Since $(\bar X_r^{t,\bar\xi})_{r\in[s,T]}$ solves the same equation, by pathwise uniqueness we deduce that
$(\bar X_r^{s,\bar X_s^{t,\bar\xi}})_{r\in[s,T]}$ and $(\bar X_r^{t,\bar\xi})_{r\in[s,T]}$ are indistinguishable, namely \eqref{FlowX1} holds.

\vspace{2mm}

\emph{Flow property \eqref{FlowX2}.} Recall that $(\P_{\text{\tiny$\bar X_s^{t,\bar\xi}$}}^{\text{\tiny$\bar\Fc_s^\mu$}})_{s\in[t,T]}$ stands for the stochastic process $(\P_s^{t,\pi})_{s\in[t,T]}$ introduced in Lemma \ref{L:P}. In the present proof it is more convenient to adopt the notation $\P_s^{t,\pi}$ instead of $\P_{\text{\tiny$\bar X_s^{t,\bar\xi}$}}^{\text{\tiny$\bar\Fc_s^\mu$}}$. Notice that, by \eqref{FlowX1}, we have $\P_r^{t,\pi}=\P_r^{s,\P_s^{t,\pi}}$, for all $r\in[s,T]$, $\bar\P$-almost surely. Therefore
\begin{align*}
\bar X_r^{t,x,\pi} \ &= \ \bar X_s^{t,x,\pi} + \int_s^r b\big(u,\bar X_u^{t,x,\pi},\P_u^{t,\pi},\bar I_u\big)\,du + \int_s^r \sigma\big(u,\bar X_u^{t,x,\pi},\P_u^{t,\pi},\bar I_u\big)\,d\bar B_u \notag \\
&= \ \bar X_s^{t,x,\pi} + \int_s^r b\big(u,\bar X_u^{t,x,\pi},\P_u^{s,\P_s^{t,\pi}},\bar I_u\big)\,du + \int_s^r \sigma\big(u,\bar X_u^{t,x,\pi},\P_u^{s,\P_s^{t,\pi}},\bar I_u\big)\,d\bar B_u,
\end{align*}
for all $r\in[s,T]$, $\bar\P$-a.s.. On the other hand, consider the process $(\bar X_r^{s,\bar X_s^{t,x,\pi},\P_s^{t,\pi}})_{r\in[s,T]}$ solution to equation \eqref{StateEq2Rand_general} with initial conditions $t=s$, $\bar\eta=\bar X_s^{t,x,\pi}$, $\bar\Pi=\P_s^{t,\pi}$. Then, we see that $(\bar X_r^{s,\bar X_s^{t,x,\pi},\P_s^{t,\pi}})_{r\in[s,T]}$ and $(\bar X_r^{t,x,\pi})_{r\in[s,T]}$ solve the same equation. It follows that they are indistinguishable, namely \eqref{FlowX2} holds.
\ep

\subsection{Randomized dynamic programming principle}

We begin proving the following identification result between $V$ and $Y^{t,x,\pi}$.

\begin{Lemma}\label{Identification}
Under Assumptions {\bf (A1)} and {\bf (A2)}, for every $(t,x,\bar\xi)\in[0,T]\times\R^n\times L^2(\bar\Omega,\bar\Gc,\bar\P;\R^n)$, with $\pi=\P_{\text{\tiny$\bar\xi$}}$ under $\bar\P$, we have
\[
\E\big[V(s,\bar X_s^{t,x,\pi},\P_{\text{\tiny$\bar X_s^{t,\bar\xi}$}}^{\text{\tiny$\bar\Fc_s^\mu$}})\big] \ = \ Y_s^{t,x,\pi},
\]
$\P^1$-a.s., for all $s\in[t,T]$.
\end{Lemma}
\textbf{Proof.}
Fix $(t,s,x,\bar\xi)\in[0,T]\times[0,T]\times\R^n\times L^2(\bar\Omega,\bar\Gc,\bar\P;\R^n)$, with $t\leq s$ and $\pi=\P_{\text{\tiny$\bar\xi$}}$ under $\bar\P$. Using the same notations as in the proof of Theorem \ref{Thm:Feynman-Kac}, let us consider, for every $n\in\N$, formula \eqref{dual_n} with $\bar t$ and $s$ replaced respectively by $s$ and $T$:
\[
Y_s^{n,t,x,\pi} \ = \ \esssup_{\hat\nu\in\hat\Vc^n} \E^{\hat\nu}\bigg[\int_s^T \E\big[f(r,\bar X_r^{t,x,\pi},\P_{\text{\tiny$\bar X_r^{t,\bar\xi}$}}^{\text{\tiny$\bar\Fc_r^\mu$}},\bar I_r)\big]\,dr + \E\big[g(\bar X_T^{t,x,\pi},\P_{\text{\tiny$\bar X_T^{t,\bar\xi}$}}^{\text{\tiny$\bar\Fc_T^\mu$}})\big]\bigg|\Fc_s^\mu\bigg].
\]
Letting $n\rightarrow\infty$, we obtain
\[
Y_s^{t,x,\pi} \ = \ \esssup_{\hat\nu\in\hat\Vc} \E^{\hat\nu}\bigg[\int_s^T \E\big[f(r,\bar X_r^{t,x,\pi},\P_{\text{\tiny$\bar X_r^{t,\bar\xi}$}}^{\text{\tiny$\bar\Fc_r^\mu$}},\bar I_r)\big]\,dr + \E\big[g(\bar X_T^{t,x,\pi},\P_{\text{\tiny$\bar X_T^{t,\bar\xi}$}}^{\text{\tiny$\bar\Fc_T^\mu$}})\big]\bigg|\Fc_s^\mu\bigg].
\]
Reasoning as in Remark \ref{R:hatV}, we can show that the right-hand side of \eqref{Ydual} does not change if we take the supremum over $\Vc$. In other words, \eqref{Ydual} can be equivalently written as follows:
\[
Y_s^{t,x,\pi} \ = \ \esssup_{\nu\in\Vc} \E^\nu\bigg[\int_s^T \E\big[f(r,\bar X_r^{t,x,\pi},\P_{\text{\tiny$\bar X_r^{t,\bar\xi}$}}^{\text{\tiny$\bar\Fc_r^\mu$}},\bar I_r)\big]\,dr + \E\big[g(\bar X_T^{t,x,\pi},\P_{\text{\tiny$\bar X_T^{t,\bar\xi}$}}^{\text{\tiny$\bar\Fc_T^\mu$}})\big]\bigg|\Fc_s^\mu\bigg].
\]
Then, we see that the claim follows if we prove the following equality: $\P^1$-a.s.
\begin{equation}\label{EquivFlow_Proof}
\E\big[V(s,\bar X_s^{t,x,\pi}\!,\P_{\text{\tiny$\bar X_s^{t,\bar\xi}$}}^{\text{\tiny$\bar\Fc_s^\mu$}})\big] = \esssup_{\nu\in\Vc} \E^\nu\!\bigg[\!\int_s^T \!\!\! \E\big[f(r,\bar X_r^{t,x,\pi}\!,\P_{\text{\tiny$\bar X_r^{t,\bar\xi}$}}^{\text{\tiny$\bar\Fc_r^\mu$}},\bar I_r)\big]dr + \E\big[g(\bar X_T^{t,x,\pi}\!,\P_{\text{\tiny$\bar X_T^{t,\bar\xi}$}}^{\text{\tiny$\bar\Fc_T^\mu$}})\big]\bigg|\Fc_s^\mu\!\bigg].
\end{equation}
As in the proof of Lemma \ref{L:Flow}, it is more convenient to adopt the notation $\P_s^{t,\pi}$ instead of $\P_{\text{\tiny$\bar X_s^{t,\bar\xi}$}}^{\text{\tiny$\bar\Fc_s^\mu$}}$ (recall that $(\P_{\text{\tiny$\bar X_s^{t,\bar\xi}$}}^{\text{\tiny$\bar\Fc_s^\mu$}})_{s\in[t,T]}$ stands for the stochastic process $(\P_s^{t,\pi})_{s\in[t,T]}$ introduced in Lemma \ref{L:P}). Then, from the flow properties \eqref{FlowX1} and \eqref{FlowX2}, we have
\begin{equation}\label{EquivFlow_Proof2}
Y_s^{t,x,\pi}\! = \esssup_{\nu\in\Vc} \E^\nu\!\bigg[\!\int_s^T \!\!\!\!\E\big[f(r,\bar X_r^{s,\bar X_s^{t,x,\pi},\P_s^{t,\pi}}\!,\P_r^{s,\P_s^{t,\pi}}\!,\bar I_r)\big]dr   + \E\big[g(\bar X_T^{s,\bar X_s^{t,x,\pi},\P_s^{t,\pi}}\!,\P_T^{s,\P_s^{t,\pi}}\!)\big]\bigg|\Fc_s^\mu\bigg].
\end{equation}
Now, notice that $\bar X_s^{t,x,\pi}\in L^2(\bar\Omega,\bar\Fc_s^{B,\mu},\bar\P;\R^n)$, so that it is the $L^2$-limit (and also pointwise $\bar\P$-a.s.) of a sequence $(\bar X_m)_{m\geq0}\subset L^2(\bar\Omega,\bar\Fc_s^{B,\mu},\bar\P;\R^n)$, where each $\bar X_m$ takes only a finite number of values. Similarly, $\P_s^{t,\pi}$ is a random variable $\P_s^{t,\pi}\colon\bar\Omega\rightarrow\mathscr P_{\text{\tiny$2$}}(\R^n)$ such that $\bar\E[\|\P_s^{t,\pi}\|_{\text{\tiny$2$}}^2]<\infty$. Therefore, by Lemma \ref{L:W2Sequence} there exists a sequence $(\P_m)_{m\geq0}$ of $\bar\Fc_s^{B,\mu}$-measurable maps $\P_m\colon\bar\Omega\rightarrow\mathscr P_{\text{\tiny$2$}}(\R^n)$, with $\bar\E[\|\P_m\|_{\text{\tiny$2$}}^2]<\infty$ and each $\P_m$ takes only a finite number values, such that $\bar\E[\Wc_{\text{\tiny$2$}}(\P_m,\P_s^{t,\pi})^2]\rightarrow0$ as $m$ goes to infinity (and also $\Wc_{\text{\tiny$2$}}(\P_m,\P_s^{t,\pi})\rightarrow0$ pointwise $\bar\P$-a.s.). In particular, for every $m\geq0$, we have
\[
\bar X_m \ = \ \sum_{k=0}^{K_m} x_{m,k}\,1_{E_{m,k}}, \qquad\qquad\quad \P_m \ = \ \sum_{k=0}^{K_m} \pi_{m,k}\,1_{E_{m,k}},
\]
for some $K_m\in\N$, $x_{m,k}\in\R^n$, $\pi_{m,k}\in\mathscr P_{\text{\tiny$2$}}(\R^n)$, $E_{m,k}\in\bar\Fc_s^{B,\mu}$, with $(E_{m,k})_k$ being a partition of $\bar\Omega$. For every $m\geq0$, consider the process $(\bar X_r^{s,\bar X_m,\P_m})_{r\in[s,T]}$, solution to equation \eqref{StateEq2Rand_general} with initial conditions $t=s$, $\bar\eta=\bar X_m$, $\bar\Pi=\P_m$. Recall from Remark \ref{R:X_Partition}, we have that the stochastic processes $(\bar X_r^{s,\bar X_m,\P_m})_{r\in[s,T]}$ and $(\sum_{k=0}^{K_m} \bar X_r^{s,x_{m,k},\pi_{m,k}}\,1_{E_{m,k}})_{r\in[s,T]}$ are indistinguishable.

Notice that, for every $\nu\in\Vc$, we have, from Corollary \ref{C:EquivConditional}, $\P^1$-a.s.,
\begin{align}
&\E\big[V(s,\bar X_m,\P_m)\big] \ = \ \sum_{k=0}^{K_m} \E\big[V(s,x_{m,k},\pi_{m,k})\,1_{E_{m,k}}\big] \notag \\
&= \ \sum_{k=0}^{K_m} \E\bigg[1_{E_{m,k}} \esssup_{\nu\in\Vc}\E^\nu\bigg[\int_s^T \E\big[f\big(r,\bar X_r^{s,x_{m,k},\pi_{m,k}},\P_r^{s,\pi_{m,k}},\bar I_r\big)\big]\,dr \notag \\
&\quad \ + \E\big[g\big(\bar X_T^{s,x_{m,k},\pi_{m,k}},\P_T^{s,\pi_{m,k}}\big)\big]\bigg|\Fc_s^\mu\bigg]\bigg] \notag \\
&= \ \esssup_{\nu\in\Vc}\E^\nu\bigg[\int_s^T \E\big[f\big(r,\bar X_r^{s,\bar X_m,\P_m},\P_r^{s,\P_m},\bar I_r\big)\big]\,dr + \E\big[g\big(\bar X_T^{s,\bar X_m,\P_m},\P_T^{s,\P_m}\big)\big]\bigg|\Fc_s^\mu\bigg]. \label{VY_m=esssup}
\end{align}
From the continuity of the map $(y,\gamma)\mapsto V(s,y,\gamma)$ stated in Proposition \ref{P:VCont}, and the growth condition \eqref{GrowthV}, we see that
\begin{equation}\label{VY_m-->VX}
\E\big[V(s,\bar X_m,\P_m)\big] \ \underset{\P^1\text{-a.s.}}{\overset{m\rightarrow\infty}{\longrightarrow}} \ \E\big[V(s,\bar X_s^{t,x,\pi},\P_s^{t,\pi})\big].
\end{equation}
On the other hand, using estimate \eqref{EstimateMKV} and proceeding as in the proof of inequality \eqref{V-Vm} in Proposition \ref{P:VCont}, we can prove the following convergence:
\begin{align}
&\esssup_{\nu\in\Vc} \E^\nu\bigg[\int_s^T \E\big[f\big(r,\bar X_r^{s,\bar X_m,\P_m},\P_r^{s,\P_m},\bar I_r\big)\big]\,dr + \E\big[g\big(\bar X_T^{s,\bar X_m,\P_m},\P_T^{s,\P_m}\big)\big]\bigg|\Fc_s^\mu\bigg] \label{esssupY^m-->X} \\
&\underset{\P^1\text{-a.s.}}{\overset{m\rightarrow\infty}{\longrightarrow}} \ \esssup_{\nu\in\Vc} \E^\nu\bigg[\int_s^T \E\big[f\big(r,\bar X_r^{s,\bar X_s^{t,x,\pi},\P_s^{t,\pi}},\P_r^{s,\P_s^{t,\pi}},\bar I_r\big)\big]\,dr + \E\big[g\big(\bar X_T^{s,\bar X_s^{t,x,\pi},\P_s^{t,\pi}},\P_T^{s,\P_s^{t,\pi}}\big)\big]\bigg|\Fc_s^\mu\bigg]. \notag
\end{align}
Hence, by \eqref{VY_m-->VX} and \eqref{esssupY^m-->X}, together with equalities \eqref{EquivFlow_Proof2} and \eqref{VY_m=esssup}, we see that \eqref{EquivFlow_Proof} holds, therefore the claim follows.
\ep

\vspace{3mm}

We can now state the main result of this section.

\begin{Theorem}\label{thm:rdpp}
Suppose that Assumptions {\bf (A1)} and {\bf (A2)} hold. Then, for every $(t,s,x,\bar\xi)\in[0,T]\times[0,T]\times\R^n\times L^2(\bar\Omega,\bar\Gc,\bar\P;\R^n)$, with $t\leq s$ and $\pi=\P_{\text{\tiny$\bar\xi$}}$ under $\bar\P$, we have
\[
V(t,x,\pi) \ = \ \sup_{\nu\in\Vc}\E^\nu\bigg[\int_t^s \E\big[f(r,\bar X_r^{t,x,\pi},\P_{\text{\tiny$\bar X_r^{t,\bar\xi}$}}^{\text{\tiny$\bar\Fc_r^\mu$}},\bar I_r)\big]\,dr + \E\big[V(s,\bar X_s^{t,x,\pi},\P_{\text{\tiny$\bar X_s^{t,\bar\xi}$}}^{\text{\tiny$\bar\Fc_s^\mu$}})\big]\bigg].
\]
\end{Theorem}
\textbf{Proof.}
Fix $(t,s,x,\bar\xi)\in[0,T]\times[0,T]\times\R^n\times L^2(\bar\Omega,\bar\Gc,\bar\P;\R^n)$, with $t\leq s$ and $\pi=\P_{\text{\tiny$\bar\xi$}}$ under $\bar\P$. Recall that by \eqref{DPPBSDE} we have, $\P^1$-a.s.,
\[
Y_t^{t,x,\pi} \ = \ \sup_{\nu\in\Vc} \E^\nu\bigg[\int_t^s \E\big[f(r,\bar X_r^{t,x,\pi},\P_{\text{\tiny$\bar X_r^{t,\bar\xi}$}}^{\text{\tiny$\bar\Fc_r^\mu$}},\bar I_r)\big]\,dr + Y_s^{t,x,\pi}\bigg].
\]
Then, the claim follows from Lemma \ref{Identification}.
\ep

\begin{Remark}\label{R:HJB}
{\rm
\emph{Hamilton-Jacobi-Bellman equation for $V$ and $V_{\textup{\tiny MKV}}$.} Let us derive, in a formal way, the dynamic programming equation for the value function $V$. We proceed as usual, starting from the dynamic programming principle of Theorem \ref{thm:rdpp} and applying It\^o's formula (see the Appendix in \cite{CarmDel14}) to the process $V(s,\bar X_s^{t,x,\pi},\P_{\text{\tiny$\bar X_s^{t,\bar\xi}$}}^{\text{\tiny$\bar\Fc_s^\mu$}})$, supposing that $V$ is smooth enough. Then, it is easy to see that the Hamilton-Jacobi-Bellman equation for $V$ takes the following form (see Section 6 of \cite{car12} for the definition of $\partial_\pi$):
\beqs
\partial _t V(t,x,\pi) +  \sup_{a \in A} \Big\{ f(t,x,\pi,a) +  b(t,x,\pi,a).\Big[ \partial_x V(t,x,\pi) + \int_{\R^n} \partial_\pi V(t,x',\pi)(x) \pi(dx') \Big]  & & \\
\;\;\; + \; \frac{1}{2}{\rm tr}(\sigma\sigma\trans(t,x,\pi,a)\Big[ \partial_x^2 V(t,x,\pi) + \int_{\R^n} \partial_x\partial_\pi V(t,x',\pi)(x) \pi(dx') \Big] \Big\} &=& 0,
\enqs
for all $(t,x,\pi)\in[0,T)\times\R^n\times\mathscr P_{\text{\tiny$2$}}(\R^n)$, with terminal condition
\[
V(T,x,\pi) \ = \ g(x,\pi), \qquad \text{for all }(x,\pi)\in\R^n\times\mathscr P_{\text{\tiny$2$}}(\R^n).
\]
We can also derive the Hamilton-Jacobi-Bellman equation for the value function $V_{\textup{\tiny MKV}}$ defined by \eqref{V_MKV_App}. From Proposition \ref{P:V_MKV=V}, we have 
\beqs
V_{\textup{\tiny MKV}}(t,\xi) &=& \E[V(t,\xi,\pi)] \; = \; \int_{\R^n} V(t,x,\pi)\,\pi(dx), 
\enqs
for all $(t,\xi)\in[0,T]\times L^2(\Omega,\Gc,\P;\R^n)$, with $\pi=\P_{\text{\tiny$\xi$}}$ under $\P$. From the above formula we see that $V_{\textup{\tiny MKV}}$ depends on $\xi$ only through its law $\pi$. In other words, $V_{\textup{\tiny MKV}}(t,\xi)=V_{\textup{\tiny MKV}}(t,\xi')$ whenever $\xi$ and $\xi'$ have the same law $\pi$. Then, by an abuse of notation, we suppose that $V_{\textup{\tiny MKV}}$ is defined on $[0,T]\times\mathscr P_{\text{\tiny$2$}}(\R^n)$ with $V_{\textup{\tiny MKV}}(t,\pi)$ given by $V_{\textup{\tiny MKV}}(t,\xi)$, for some $\xi$ such that $\pi=\P_{\text{\tiny$\xi$}}$. Now, recalling the definition of the derivative $\partial_\pi$, we obtain 
\beqs
\partial_t V_{\textup{\tiny MKV}}(t,\pi) &=&  \E[ \partial_t V(t,\xi,\pi)], \\
\partial_\pi V_{\textup{\tiny MKV}}(t,\pi)(x) &=& \partial_x V(t,x,\pi) +  \E[ \partial_\pi V(t,\xi,\pi)(x) ], \\
\partial_x \partial_\pi V_{\textup{\tiny MKV}}(t,\pi)(x) &=& \partial_x^2 V(t,x,\pi) +  \E[ \partial_x\partial_\pi V(t,\xi,\pi)(x) ].
\enqs
Integrating with respect to $\pi$ in the Hamilton-Jacobi-Bellman equation of $V$, we obtain the following dynamic programming equation for  $V_{\textup{\tiny MKV}}$: 
\beq
\partial_t V_{\textup{\tiny MKV}}(t,\pi) + 
\int_{\R^n} \sup_{a\in A} \big[ f(t,x,\pi,a) + b(t,x,\pi,a). \partial_\pi V_{\textup{\tiny MKV}}(t,\pi)(x) & & \label{HJB_V_MKV} \\
\;\;\;\;\; + \;  \frac{1}{2}{\rm tr}\big(\sigma\sigma\trans(t,x,\pi,a) \partial_x \partial_\pi V_{\textup{\tiny MKV}}(t,\pi)(x) \big) \big] \pi(dx) &=& 0, \notag
\enq
for all $(t,\pi)\in[0,T]\times\mathscr P_{\text{\tiny$2$}}(\R^n)$, with terminal condition
\[
V_{\textup{\tiny MKV}}(T,\pi) \ = \ \int_{\R^n} g(x,\pi)\,\pi(dx), \qquad \text{for all }\pi\in\mathscr P_{\text{\tiny$2$}}(\R^n).
\]
Notice that if the supremum inside the integral in \eqref{HJB_V_MKV} is attained at some $\hat a(x)$, for some map $\hat a\colon\R^n\rightarrow A$ Lipschitz continuous in $x$, then the above equation can be written as (we denote by $L(\R^n;A)$ the set of Lipschitz continuous maps from $\R^n$ into $A$)
\beqs
\partial_t V_{\textup{\tiny MKV}}(t,\pi) + 
\sup_{\tilde\alpha\in L(\R^n;A)} \int_{\R^n} \big[ f(t,x,\pi,\tilde\alpha(x)) + b(t,x,\pi,\tilde\alpha(x)). \partial_\pi V_{\textup{\tiny MKV}}(t,\pi)(x) & & \\
\;\;\;\;\; + \; \frac{1}{2}{\rm tr}\big(\sigma\sigma\trans(t,x,\pi,\tilde\alpha(x)) \partial_x \partial_\pi V_{\textup{\tiny MKV}}(t,\pi)(x) \big) \big] \pi(dx) &=& 0.
\enqs
This latter is the Hamilton-Jacobi-Bellman equation obtained in \cite{PhamWei} under the assumption that the optimization in the McKean-Vlasov control problem is performed only over the class of Lipschitz continuous closed-loop controls.
\ep}
\end{Remark}

\appendix

\renewcommand\thesection{}

\section{}

\renewcommand\thesection{\Alph{subsection}}

\renewcommand\thesubsection{\Alph{subsection}}

\subsection{Some convergence results with respect to the 2-Wasserstein metric $\Wc_{\text{\tiny$2$}}$}
\label{App:Wass}

\setcounter{Theorem}{0}
\setcounter{Definition}{0}
\setcounter{Proposition}{0}
\setcounter{Assumption}{0}
\setcounter{Lemma}{0}
\setcounter{Corollary}{0}
\setcounter{Remark}{0}
\setcounter{Example}{0}
\setcounter{equation}{0}

\begin{Lemma}[Skorohod's representation theorem for $\Wc_{\text{\tiny$2$}}$-convergence]\label{L:W2Skorohod}
Let $(\pi_m)_m$ be a sequence in $\mathscr P_{\text{\tiny$2$}}(\R^n)$ such that $\Wc_{\text{\tiny$2$}}(\pi_m,\pi)\rightarrow0$, for some $\pi\in\mathscr P_{\text{\tiny$2$}}(\R^n)$. Then, there exists a sequence of random variables $(\xi_m)_m\subset L^2(\Omega,\Gc,\P;\R^n)$, with $\P_{\text{\tiny$\xi_m$}}=\pi_m$, converging pointwise $\P$-a.s. and in $L^2(\Omega,\Gc,\P;\R^n)$ to some 
$\xi\in L^2(\Omega,\Gc,\P;\R^n)$, with $\P_{\text{\tiny$\xi$}}=\pi$.
\end{Lemma}
\textbf{Proof.}
By Theorem 6.9 and point (i) of Definition 6.8 in \cite{Villani}, we have that $\Wc_{\text{\tiny$2$}}(\pi_m,\pi)\rightarrow0$ is equivalent to:
\begin{equation}\label{SkorohodProof}
\pi_m \ \underset{\text{weakly}}{\overset{m\rightarrow\infty}{\longrightarrow}} \ \pi \qquad\qquad \text{ and } \qquad\qquad \int_{\R^n}|x|^2\,\pi_m(dx) \ \overset{m\rightarrow\infty}{\longrightarrow} \ \int_{\R^n}|x|^2\,\pi(dx).
\end{equation}
Then, by the classical Skorohod representation theorem for weak convergence, there exist random variables $\xi_m,\xi\in L^2(\Omega,\Gc,\P;\R^n)$, with $\P_{\text{\tiny$\xi_m$}}=\pi_m$ and $\P_{\text{\tiny$\xi$}}=\pi$, such that $\xi_m$ converges pointwise $\P$-a.s. to $\xi$. It remains to prove the convergence in $L^2(\Omega,\Gc,\P;\R^n)$. To this end, we notice that \eqref{SkorohodProof} implies $\E[|\xi_m|^2]\rightarrow\E[|\xi|^2]$. Therefore, by Theorem II.6.5 in \cite{shiryaev}, the sequence $(|\xi_m|^2)_m$ is uniformly integrable. Then, it follows that $\xi_m\rightarrow\xi$ in $L^2(\Omega,\Gc,\P;\R^n)$.
\ep

\vspace{3mm}

\begin{Lemma}\label{L:CountableDeterm}
There exists a countable convergence determining class $(\varphi_k)_{k\geq1}\subset C_{\text{\tiny$2$}}(\R^n)$ for the $\Wc_{\text{\tiny$2$}}$-convergence. In other words, given $\pi_1,\pi_2,\ldots,\pi\in\mathscr P_{\text{\tiny$2$}}(\R^n)$, we have:
\[
\Wc_{\text{\tiny$2$}}(\pi_m,\pi) \ \overset{m\rightarrow\infty}{\longrightarrow} \ 0 \quad\;\, \text{ if and only if } \quad\;\, \int_{\R^n}\varphi_k(x)\,\pi_m(dx) \ \overset{m\rightarrow\infty}{\longrightarrow} \ \int_{\R^n}\varphi_k(x)\,\pi(dx), \;\; \text{for all $k$}.
\]
\end{Lemma}
\textbf{Proof.}
Let $\pi_1,\pi_2,\ldots,\pi\in\mathscr P_{\text{\tiny$2$}}(\R^n)$. We recall from Theorem 6.9 and point (i) of Definition 6.8 in \cite{Villani} that
\[
\Wc_{\text{\tiny$2$}}(\pi_m,\pi) \ \overset{m\rightarrow\infty}{\longrightarrow} \ 0 \quad \text{ if and only if } \quad \pi_m \ \underset{\text{weakly}}{\overset{m\rightarrow\infty}{\longrightarrow}} \ \pi \text{ and } \int_{\R^n}|x|^2\,\pi_m(dx) \ \overset{m\rightarrow\infty}{\longrightarrow} \ \int_{\R^n}|x|^2\,\pi(dx).
\]
Now, it is well-known that there exists a countable convergence determining class $(\psi_h)_{h\geq1}\subset C_b(\R^n)$ (the set of real-valued continuous and bounded functions) for the weak convergence (see, for instance, Theorem 2.18 in \cite{BaCri}). In other words, we have
\[
\pi_m \ \underset{\text{weakly}}{\overset{m\rightarrow\infty}{\longrightarrow}} \ \pi \quad\;\; \text{ if and only if } \quad\;\; \int_{\R^n}\psi_h(x)\,\pi_m(dx) \ \overset{m\rightarrow\infty}{\longrightarrow} \ \int_{\R^n}\psi_h(x)\,\pi(dx), \;\; \text{for all $h$}.
\]
Then, the claim follows taking $\varphi_1(x):=|x|^2$, for every $x\in\R^n$, and $\varphi_k:=\psi_{k-1}$, for every $k\geq2$.
\ep

\begin{Lemma}\label{L:W2Sequence}
Let $(\tilde\Omega,\tilde\Fc,\tilde\P)$ be a probability space and let $\Pi\colon\tilde\Omega\rightarrow\mathscr P_{\text{\tiny$2$}}(\R^n)$ be a measurable map. Suppose that $($$\tilde\E$ denotes the $\tilde\P$-expected value$)$
\begin{equation}\label{L2normPi}
\tilde\E\big[\|\Pi\|_{\text{\tiny$2$}}^2\big] \ < \ +\infty.
\end{equation}
Then, there exists a sequence $(\Pi_m)_{m\geq1}$ of measurable maps $\Pi_m\colon\tilde\Omega\rightarrow\mathscr P_{\text{\tiny$2$}}(\R^n)$ such that:
\[
\Wc_{\text{\tiny$2$}}(\Pi_m(\tilde\omega),\Pi(\tilde\omega)) \ \overset{m\rightarrow\infty}{\longrightarrow} \ 0, \;\; \tilde\P(d\tilde\omega)\text{-a.s.,} \qquad \text{ and } \qquad \tilde\E\big[\Wc_{\text{\tiny$2$}}(\Pi_m,\Pi)^2\big] \ \overset{m\rightarrow\infty}{\longrightarrow} \ 0,
\]
where, for every $m\geq1$,
\[
\Pi_m(\tilde\omega) \ = \ \sum_{k=1}^{K_m} \pi_{m,k}\,1_{E_{m,k}}(\tilde\omega), \qquad \text{for every }\tilde\omega\in\tilde\Omega,
\]
for some finite integer $K_m\geq1$, $\pi_{m,k}\in\mathscr P_{\text{\tiny$2$}}(\R^n)$, $E_{m,k}\in\tilde\Fc$, with $(E_{m,k})_{k=1,\ldots,K_m}$ being a partition of $\tilde\Omega$.
\end{Lemma}
\textbf{Proof.}
Recall from Theorem 6.18 in \cite{Villani} that $(\mathscr P_{\text{\tiny$2$}}(\R^n),\Wc_{\text{\tiny$2$}})$ is a complete separable metric space. Then, there exists a sequence $(\pi_h)_{h\geq1}$ dense in $\mathscr P_{\text{\tiny$2$}}(\R^n)$. Now, for every $\ell,h\geq1$, define the measurable set $\bar B_{\ell,h}\in\tilde\Fc$ by
\[
\bar B_{\ell,h} \ := \ \big\{\tilde\omega\in\tilde\Omega\colon\Wc_{\text{\tiny$2$}}(\Pi(\tilde\omega),\pi_h)\leq1/\ell\big\}.
\]
We also define the disjoint measurable sets: $B_{\ell,1}:=\bar B_{\ell,1}$ and $B_{\ell,h}:=\bar B_{\ell,h}\backslash(\bar B_{\ell,1}\cup\cdots\cup\bar B_{\ell,h-1})$, for any $h\geq2$. Notice that $\tilde\Omega=\cup_{h\geq1}B_{\ell,h}$. In particular, for every $\ell\geq1$, there exists $K_\ell\geq1$ such that $\tilde\P(\cup_{h\geq K_\ell+1}B_{\ell,h})\leq1/\ell^2$. Finally, we set
\[
\bar\Pi_\ell(\tilde\omega) \ := \ \sum_{h=1}^{K_\ell} \pi_h\,1_{B_{\ell,h}\cap A_\ell}(\tilde\omega) + \delta_0\,\big(1_{(\cup_{h\geq K_\ell+1}B_{\ell,h})\cap A_\ell}(\tilde\omega) + 1_{A_\ell^c}(\tilde\omega)\big), \qquad \text{for every }\tilde\omega\in\tilde\Omega,
\]
where
\[
A_\ell \ := \ \big\{\tilde\omega\in\tilde\Omega\colon \|\Pi(\tilde\omega)\|_{\text{\tiny$2$}}^2 \ \leq \ \ell \big\}.
\]
Then, we see that (recall from \eqref{||pi||_2^2} that $\Wc_{\text{\tiny$2$}}(\delta_0,\Pi(\tilde\omega))=\|\Pi(\tilde\omega)\|_{\text{\tiny$2$}}$)
\[
\Wc_{\text{\tiny$2$}}(\bar\Pi_\ell(\tilde\omega),\Pi(\tilde\omega)) \ \leq \ \frac{1}{\ell}\,1_{(\cup_{h=1}^{K_\ell}B_{\ell,h})\cap A_\ell}(\tilde\omega) + \|\Pi(\tilde\omega)\|_{\text{\tiny$2$}}\,\big(1_{(\cup_{h\geq K_\ell+1}B_{\ell,h})\cap A_\ell}(\tilde\omega) + 1_{A_\ell^c}(\tilde\omega)\big),
\]
for all $\tilde\omega\in\tilde\Omega$. Therefore (recalling that $\tilde\P(\cup_{h\geq K_\ell+1}B_{\ell,h})\leq1/\ell^2$)
\begin{align*}
\tilde\E[\Wc_{\text{\tiny$2$}}(\bar\Pi_\ell,\Pi)^2] \ &\leq \ \frac{1}{\ell^2} + \tilde\E\big[ \|\Pi(\tilde\omega)\|_{\text{\tiny$2$}}^2\,1_{(\cup_{h\geq K_\ell+1}B_{\ell,h})\cap A_\ell}\big] + \tilde\E\big[ \|\Pi(\tilde\omega)\|_{\text{\tiny$2$}}^2\,1_{A_\ell^c}\big] \\
&\leq \ \frac{1}{\ell^2} + \ell\,\tilde\P\big((\cup_{h\geq K_\ell+1}B_{\ell,h})\cap A_\ell\big) + \tilde\E\big[ \|\Pi(\tilde\omega)\|_{\text{\tiny$2$}}^2\,1_{A_\ell^c}\big] \\
&\leq \ \frac{1}{\ell^2} + \ell\,\frac{1}{\ell^2} + \tilde\E\big[ \|\Pi(\tilde\omega)\|_{\text{\tiny$2$}}^2\,1_{A_\ell^c}\big] \ \overset{\ell\rightarrow\infty}{\longrightarrow} \ 0,
\end{align*}
where the convergence $\tilde\E\big[ \|\Pi(\tilde\omega)\|_{\text{\tiny$2$}}^2\,1_{A_\ell^c}\big]\rightarrow0$ follows from the Lebesgue dominated convergence theorem, using \eqref{L2normPi} and noting that $1_{A_\ell^c}$ converges pointwise $\tilde\P$-a.s. to zero.

Let $Y_\ell\colon\tilde\Omega\rightarrow[0,\infty)$ be the nonnegative random variable given by $Y_\ell:=\Wc_{\text{\tiny$2$}}(\bar\Pi_\ell,\Pi)$. We know that $Y_\ell\rightarrow0$, as $\ell\rightarrow\infty$, in $L^2(\tilde\Omega,\tilde\Fc,\tilde\P)$. Then, it is well-known that this implies the existence of a subsequence $(Y_{\ell_m})_{m\geq1}$ such that $Y_{\ell_m}=\Wc_{\text{\tiny$2$}}(\bar\Pi_{\ell_m},\Pi)\rightarrow0$, as $m\rightarrow\infty$, pointwise $\tilde\P$-a.s. and in $L^2(\tilde\Omega,\tilde\Fc,\tilde\P)$. Then, $(\Pi_m)_{m\geq1}$, with $\Pi_m:=\bar\Pi_{\ell_m}$, is the desired sequence.
\ep

\subsection{Proofs of Lemma \ref{L:hatP} and Lemma \ref{L:P}}
\label{App:Proofs}

\setcounter{Theorem}{0}
\setcounter{Definition}{0}
\setcounter{Proposition}{0}
\setcounter{Assumption}{0}
\setcounter{Lemma}{0}
\setcounter{Corollary}{0}
\setcounter{Remark}{0}
\setcounter{Example}{0}
\setcounter{equation}{0}

\textbf{Proof of Lemma \ref{L:hatP}.}
Recall that, by construction, the map $\bar X^{t,\bar\xi}\colon([t,T]\times\Omega\times\Omega^1,\Bc([t,T])\otimes\bar\Fc)\rightarrow(\R^n,\Bc(\R^n))$ is measurable. Therefore, up to indistinguishability, we can suppose that $\bar X^{t,\bar\xi}\colon([t,T]\times\Omega\times\Omega^1,\Bc([t,T])\otimes\Fc\otimes\Fc^1)\rightarrow(\R^n,\Bc(\R^n))$ is measurable. Since $(\bar X_s^{t,\bar\xi})_{s\in[t,T]}$ is also $(\bar\Fc_s^{B,\mu}\vee\bar\Gc)_s$-adapted, we deduce that, for every $s\in[t,T]$, the map $\bar X_s^{t,\bar\xi}\colon(\Omega\times\Omega^1,(\Gc\vee\Fc_s^B)\otimes\Fc_s^\mu)\rightarrow(\R^n,\Bc(\R^n))$ is measurable. Therefore, by estimate \eqref{EstimateRand} and Fubini's theorem, we see that, for every $\varphi\in\mathscr B_{\text{\tiny$2$}}(\R^n)$, the map
\[
\omega^1 \ \longmapsto \ \E\big[\varphi\big(\bar X_s^{t,\bar\xi}(\cdot,\omega^1)\big)\big],
\]
from $\Omega^1$ into $\R$, is $\Fc_s^\mu$-measurable. In particular, when $\varphi\in C_{\text{\tiny$2$}}(\R^N)$, the continuous process $(\E[\varphi(\bar X_s^{t,\bar\xi})])_{s\in[t,T]}$ is $\F^\mu$-predictable. Then, by Remark \ref{R:MeasWass} it follows that the process $(\hat\P_s^{t,\pi})_{s\in[t,T]}$ is $\F^\mu$-predictable.

Finally, we observe that
\[
\hat\P_s^{t,\pi}(\omega^1)[\varphi] \ = \ \E\big[\varphi\big(\bar X_s^{t,\bar\xi}(\cdot,\omega^1)\big)\big] \ = \ \bar\E\big[\varphi\big(\bar X_s^{t,\bar\xi}\big)\big|\bar\Fc_s^\mu\big](\omega^1) \ = \ \P_{\text{\tiny$\bar X_s^{t,\bar\xi}$}}^{\text{\tiny$\bar\Fc_s^\mu$}}(\omega^1)[\varphi],
\]
$\P^1(d\omega^1)$-a.s., for every $\varphi\in\mathscr B_{\text{\tiny$2$}}(\R^n)$. Let $(\varphi_k)_k\subset\mathscr B_{\text{\tiny$2$}}(\R^n)$ be a countable separating class of continuous functions, whose existence is guaranteed for instance by Theorem 2.18 in \cite{BaCri} ($\varphi_k$ can be taken even bounded). Then, there exists a unique $\P^1$-null set $N^1\in\Fc^1$ such that
\[
\hat\P_s^{t,\pi}(\omega^1)[\varphi_k] \ = \ \P_{\text{\tiny$\bar X_s^{t,\bar\xi}$}}^{\text{\tiny$\bar\Fc_s^\mu$}}(\omega^1)[\varphi_k], \qquad \text{for every }k,
\]
whenever $\omega^1\notin N^1$. Since $(\varphi_k)_k$ is separating, we conclude that $\hat\P_s^{t,\pi}$ coincides with $\P_{\text{\tiny$\bar X_s^{1,t,\bar\xi}$}}^{\text{\tiny$\bar\Fc_s^\mu$}}$ on $\Omega^1\backslash N^1$. In other words, $(\hat\P_s^{t,\pi})_{s\in[t,T]}$ is a version of $(\P_{\text{\tiny$\bar X_s^{t,\bar\xi}$}}^{\text{\tiny$\bar\Fc_s^\mu$}})_{s\in[t,T]}$.
\ep

\vspace{3mm}

\noindent\textbf{Proof of Lemma \ref{L:P}.}
Fix $t\in[0,T]$ and consider a generic $\pi\in\mathscr P_{\text{\tiny$2$}}(\R^n)$. Let $\bar\xi\in L^2(\bar\Omega,\bar\Gc,\bar\P;\R^n)$ be such that $\pi=\P_{\text{\tiny$\bar\xi$}}$ under $\bar\P$. We construct $\bar X^{t,\bar\xi}$ using Picard's iterations. More precisely, we define recursively a sequence of $\R^n$-valued processes $(\bar X^{m,t,\bar\xi})_m$ on $\bar\Omega\times[t,T]$ as follows.

\vspace{2mm}

\noindent\textbf{\emph{Recursive construction of the sequence $(\bar X^{m,t,\bar\xi})_m$.}} \emph{Definition of $\bar X^{0,t,\bar\xi}$.} We set $\bar X^{0,t,\bar\xi}\equiv0$. Defining $\hat\P^{0,t,\bar\xi}$ by formula \eqref{RCPD} with $\bar X^{0,t,\bar\xi}$ in place of $\bar X^{t,\bar\xi}$, we see that $\hat\P_s^{0,t,\pi}\equiv\delta_0$, the Dirac delta at zero, for all $s\in[t,T]$. In other words, up to a version, $(\P_{\text{\tiny$\bar X_s^{0,t,\bar\xi}$}}^{\text{\tiny$\bar\Fc_s^\mu$}})_{s\in[t,T]}$ is identically equal to $\delta_0$.

\vspace{2mm}

\noindent\emph{Definition of $\bar X^{1,t,\bar\xi}$.} The process $\bar X^{1,t,\bar\xi}$ is given by:
\[
\bar X_s^{1,t,\bar\xi} \ = \ \bar\xi + \int_t^s b\big(r,0,\delta_0,\bar I_r\big) dr + \int_t^s \sigma\big(r,0,\delta_0,\bar I_r\big) d\bar B_r,
\]
for all $s\in[t,T]$. Notice that, by construction, the map $\bar X^{1,t,\bar\xi}\colon([t,T]\times\Omega\times\Omega^1,\Bc([t,T])\otimes\bar\Fc)\rightarrow(\R^n,\Bc(\R^n))$ is measurable. Up to indistinguishability, we can suppose that $\bar X^{1,t,\bar\xi}$ $\colon([t,T]\times\Omega\times\Omega^1,\Bc([t,T])\otimes\Fc\otimes\Fc^1)\rightarrow(\R^n,\Bc(\R^n))$ is measurable. As a consequence, by Fubini's theorem, we can define the $\mathscr P_{\text{\tiny$2$}}(\R^n)$-valued $\F^\mu$-predictable stochastic process $(\hat\P_s^{1,t,\pi})_{s\in[t,T]}$ by formula \eqref{RCPD} with $\bar X^{1,t,\bar\xi}$ in place of $\bar X^{t,\bar\xi}$. Notice that $(\hat\P_s^{1,t,\pi})_{s\in[t,T]}$ is a version of $(\P_{\text{\tiny$\bar X_s^{1,t,\bar\xi}$}}^{\text{\tiny$\bar\Fc_s^\mu$}})_{s\in[t,T]}$. Moreover, from \eqref{RCPD}, we see that (using the definition of $\bar X_s^{1,t,\bar\xi}$, and the independence of $\bar\Gc$ and $\bar\Fc_\infty^B$)
\[
\hat\P_s^{1,t,\pi}(\omega^1)[\varphi] \ = \ \E\big[\varphi\big(\bar X_s^{1,t,\bar\xi}(\cdot,\omega^1)\big)\big] \ = \ \int_{\R^n} \Phi_{1,\varphi}(\omega^1,s,x)\,\pi(dx),
\]
for every $\omega^1\in\Omega^1$ and $\varphi\in\mathscr B_{\text{\tiny$2$}}(\R^n)$, where $\Phi_{1,\varphi}\colon\Omega^1\times[t,T]\times\R^n\rightarrow\R$ is measurable, with at most quadratic growth in $x$ uniformly with respect to $(\omega^1,s)$, and it is given by
\[
\Phi_{1,\varphi}(\omega^1,s,x) \ := \ \E\bigg[\varphi\bigg(x + \int_t^s b\big(r,0,\delta_0,\bar I_r(\cdot,\omega^1)\big) dr + \int_t^s \sigma\big(r,0,\delta_0,\bar I_r(\cdot,\omega^1)\big) d B_r\bigg)\bigg].
\]
Then, we see that the map $\hat\P_\cdot^{1,t,\cdot}[\varphi]\colon\Omega^1\times[t,T]\times\mathscr P_{\text{\tiny$2$}}(\R^n)\rightarrow\R$ is measurable. Indeed, when $\Phi_{1,\varphi}(\omega^1,s,x)=\ell(\omega^1,s)h(x)$, for some measurable functions $\ell$ and $h$, with $\ell$ bounded and $h$ with at most quadratic growth (namely $h\in\mathscr B_{\text{\tiny$2$}}(\R^n)$), the result follows from Remark \ref{R:MeasWass}. The general case can be proved by a monotone class argument.

Using again Remark \ref{R:MeasWass}, we conclude that the map $\hat\P_\cdot^{1,t,\cdot}\colon\Omega^1\times[t,T]\times\mathscr P_{\text{\tiny$2$}}(\R^n)\rightarrow\mathscr P_{\text{\tiny$2$}}(\R^n)$ is measurable.

\vspace{2mm}

\noindent\emph{Definition of $\bar X^{m+1,t,\bar\xi}$, for every integer $m\geq1$.} We define $\bar X^{m+1,t,\bar\xi}$ recursively, assuming that $\bar X^{m,t,\bar\xi}$ has already been defined. We also assume that the map $\bar X^{m,t,\bar\xi}\colon([t,T]\times\Omega\times\Omega^1,\Bc([t,T])\otimes\Fc\otimes\Fc^1)\rightarrow(\R^n,\Bc(\R^n))$ is measurable and that $(\hat\P_s^{m,t,\pi})_{s\in[t,T]}$ is the $\mathscr P_{\text{\tiny$2$}}(\R^n)$-valued $\F^\mu$-predictable stochastic process given by formula \eqref{RCPD} with $\bar X^{m,t,\bar\xi}$ in place of $\bar X^{t,\bar\xi}$. Moreover, we suppose that the map $\hat\P_\cdot^{m,t,\cdot}\colon\Omega^1\times[t,T]\times\mathscr P_{\text{\tiny$2$}}(\R^n)\rightarrow\mathscr P_{\text{\tiny$2$}}(\R^n)$ is measurable. Notice that $(\hat\P_s^{m,t,\pi})_{s\in[t,T]}$ is a version of $(\P_{\text{\tiny$\bar X_s^{m,t,\bar\xi}$}}^{\text{\tiny$\bar\Fc_s^\mu$}})_{s\in[t,T]}$.

Then, we define $\bar X^{m+1,t,\bar\xi}$ as follows:
\[
\bar X_s^{m+1,t,\bar\xi} \ = \ \bar\xi + \int_t^s b\big(r,\bar X_r^{m,t,\bar\xi},\hat\P_r^{m,t,\pi},\bar I_r\big) dr + \int_t^s \sigma\big(r,\bar X_r^{m,t,\bar\xi},\hat\P_r^{m,t,\pi},\bar I_r\big) d\bar B_r,
\]
for all $s\in[t,T]$. Notice that, by construction, the map $\bar X^{m+1,t,\bar\xi}\colon([t,T]\times\Omega\times\Omega^1,\Bc([t,T])\otimes\bar\Fc)\rightarrow(\R^n,\Bc(\R^n))$ is measurable. Therefore, up to indistinguishability, we can suppose that $\bar X^{m+1,t,\bar\xi}\colon([t,T]\times\Omega\times\Omega^1,\Bc([t,T])\otimes\Fc\otimes\Fc^1)\rightarrow(\R^n,\Bc(\R^n))$ is measurable. Then, by Fubini's theorem, we can define the $\mathscr P_{\text{\tiny$2$}}(\R^n)$-valued $\F^\mu$-predictable stochastic process $(\hat\P_s^{m+1,t,\pi})_{s\in[t,T]}$ by formula \eqref{RCPD} with $\bar X^{m+1,t,\bar\xi}$ in place of $\bar X^{t,\bar\xi}$, namely
\[
\hat\P_s^{m+1,t,\pi}(\omega^1)[\varphi] \ = \ \E\big[\varphi\big(\bar X_s^{m+1,t,\bar\xi}(\cdot,\omega^1)\big)\big],
\]
for every $\omega^1\in\Omega^1$, $\varphi\in\mathscr B_{\text{\tiny$2$}}(\R^n)$, $s\in[t,T]$. In particular, we have
\begin{align*}
\hat\P_s^{m+1,t,\pi}(\omega^1)[\varphi] \ &= \ \E\bigg[\varphi\bigg(\bar\xi + \int_t^s b\big(r,\bar\xi+\cdots,\hat\P_r^{m,t,\pi}(\omega^1),\bar I_r(\cdot,\omega^1)\big) dr \\
&\quad \ + \int_t^s \sigma\big(r,\bar\xi+\cdots,\hat\P_r^{m,t,\pi}(\omega^1),\bar I_r(\cdot,\omega^1)\big) d\bar B_r\bigg)\bigg] \\
&= \ \int_{\R^n} \Phi_{m+1,\varphi}(\omega^1,s,x,\pi)\,\pi(dx),
\end{align*}
for some measurable $\Phi_{m+1,\varphi}\colon\Omega^1\times[t,T]\times\R^n\times\mathscr P_{\text{\tiny$2$}}(\R^n)\rightarrow\R$, with at most quadratic growth in $(x,\pi)$ uniformly with respect to $(\omega^1,s)$ (the dependence of $\Phi_{m+1,\varphi}$ on $\pi$ is due to the presence of $\hat\P_r^{m,t,\pi}$). Then, we see that the map $\hat\P_\cdot^{m+1,t,\cdot}[\varphi]\colon\Omega^1\times[t,T]\times\mathscr P_{\text{\tiny$2$}}(\R^n)\rightarrow\R$ is measurable, as it can be deduced using a monotone class argument, first taking $\Phi_{m+1,\varphi}$ of the form $\Phi_{m+1,\varphi}(\omega^1,s,x,\pi)=\ell(\omega^1,s,\pi)h(x)$, for some $h\in\mathscr B_{\text{\tiny$2$}}(\R^n)$, and some measurable function $\ell$ with at most quadratic growth in $\pi$ uniformly with respect to $(\omega^1,s)$. Then, by Remark \ref{R:MeasWass}, we see that the map $\hat\P_\cdot^{m+1,t,\cdot}\colon\Omega^1\times[t,T]\times\mathscr P_{\text{\tiny$2$}}(\R^n)\rightarrow\mathscr P_{\text{\tiny$2$}}(\R^n)$ is measurable.

\vspace{2mm}

\noindent\textbf{End of the proof of Lemma \ref{L:P}.} Now that we have constructed the sequence $(\bar X^{m,t,\bar\xi})_m$, we notice that it can be proved (proceeding for instance along the same lines as in the proof of Theorem IX.2.1 in \cite{RevuzYor}) that
\begin{equation}\label{ConvProb_supX}
\sup_{s\in[t,T]}\big|\bar X_s^{m,t,\bar\xi}-\bar X_s^{t,\bar\xi}\big|\ \underset{m\rightarrow\infty}{\overset{\bar\P}{\longrightarrow}} \ 0,
\end{equation}
where the convergence holds in probability. Fix $s\in[t,T]$ and let us prove that \eqref{ConvProb_supX} implies the following convergence in probability:
\begin{equation}\label{ConvProbW2}
\Wc_{\text{\tiny$2$}}\big(\hat\P_s^{m,t,\pi},\hat\P_s^{t,\pi}\big) \ \underset{m \rightarrow\infty}{\overset{\P^1}{\longrightarrow}} 0.
\end{equation}
In order to prove \eqref{ConvProbW2}, it is enough to show that every subsequence $(\hat\P_s^{m_\ell,t,\pi})_\ell$ admits a subsubsequence $(\hat\P_s^{m_{\ell_h},t,\pi})_h$ for which \eqref{ConvProbW2} holds. Let us fix a subsequence $(\hat\P_s^{m_\ell,t,\pi})_\ell$. We begin noting that, by \eqref{ConvProb_supX}, we have, for every $\varphi\in C_{\text{\tiny$2$}}(\R^n)$,
\[
\hat\P_s^{m_\ell,t,\pi}[\varphi] \ \underset{\ell\rightarrow\infty}{\overset{\P^1}{\longrightarrow}} \ \hat\P_s^{t,\pi}[\varphi].
\]
Let $(\varphi_k)_k\subset C_{\text{\tiny$2$}}(\R^n)$ be a countable convergence determining class for the $\Wc_{\text{\tiny$2$}}$-convergence, whose existence follows from Lemma \ref{L:CountableDeterm}. Then, there exists a unique $\P^1$-null set $N^1\in\Fc^1$ and a subsubsequence $(\hat\P_s^{m_{\ell_h},t,\pi})_h$ such that, for all $\omega^1\in\Omega^1\backslash N^1$,
\[
\hat\P_s^{m_{\ell_h},t,\pi}(\omega^1)[\varphi_k] \ \overset{h\rightarrow\infty}{\longrightarrow} \ \hat\P_s^{t,\pi}(\omega^1)[\varphi_k], \qquad \text{for every }k.
\]
By Theorem 6.9 in \cite{Villani} it follows that, for all $\omega^1\in\Omega^1\backslash N^1$,
\[
\Wc_{\text{\tiny$2$}}\big(\hat\P_s^{m_{\ell_h},t,\pi}(\omega^1),\hat\P_s^{m,t,\pi}(\omega^1)\big) \ \overset{h\rightarrow\infty}{\longrightarrow} 0.
\]
In particular, the above convergence holds in probability. This concludes the proof of \eqref{ConvProbW2}.

Notice that convergence \eqref{ConvProbW2} holds for every $s\in[t,T]$ and $\pi\in\mathscr P_{\text{\tiny$2$}}(\R^n)$. Moreover, for every $m\in\N$, $\hat\P_\cdot^{m,t,\cdot}$ is jointly measurable with respect to $(\omega^1,s,\pi)$. Then, we deduce (proceeding for instance as in the first item of Exercise IV.5.17 in \cite{RevuzYor} or as in Proposition 1 of \cite{stricker_yor}) that there exists a measurable map $\P_\cdot^{t,\cdot}\colon\Omega^1\times[t,T]\times\mathscr P_{\text{\tiny$2$}}(\R^n)\rightarrow\mathscr P_{\text{\tiny$2$}}(\R^n)$ such that
\[
\Wc_{\text{\tiny$2$}}\big(\hat\P_s^{m,t,\pi},\P_s^{t,\pi}\big) \ \underset{m \rightarrow\infty}{\overset{\P^1}{\longrightarrow}} 0,
\]
for every $s\in[t,T]$ and $\pi\in\mathscr P_{\text{\tiny$2$}}(\R^n)$. This implies that $\P_s^{t,\pi}$ coincides $\P^1$-a.s. with $\hat\P_s^{t,\pi}$. By Lemma \ref{L:hatP} we conclude that $(\P_s^{t,\pi})_{s\in[t,T]}$ is a version of $(\P_{\text{\tiny$\bar X_s^{t,\bar\xi}$}}^{\text{\tiny$\bar\Fc_s^\mu$}})_{s\in[t,T]}$. 
\ep

\subsection{Stability lemma}
\label{App:Stability}

\setcounter{Theorem}{0}
\setcounter{Definition}{0}
\setcounter{Proposition}{0}
\setcounter{Assumption}{0}
\setcounter{Lemma}{0}
\setcounter{Corollary}{0}
\setcounter{Remark}{0}
\setcounter{Example}{0}
\setcounter{equation}{0}

For the proof of Theorem \ref{Thm:Equivalence}, we need the following stability result.

\begin{Lemma}\label{L:Stability}
Suppose that Assumption {\bf (A1)} holds.
\begin{itemize}
\item Let $(\tilde\Omega,\tilde\Fc,\Q)$ be a probability space, on which a $d$-dimensional Brownian motion $\tilde B=(\tilde B_t)_{t\geq0}$ is defined.
\item For every $\ell\in\N$, let $\tilde\F^\ell=(\tilde\Fc_s^\ell)_{s\geq0}$ be a filtration on $(\tilde\Omega,\tilde\Fc,\Q)$ such that $\tilde B$ is a Brownian motion with respect to $\tilde\F^\ell$.
\item For every $\ell\in\N$, let $\tilde\F^{\mu_\ell}=(\tilde\Fc_s^{\mu_\ell})_{s\geq0}$, with $\tilde\Fc_s^{\mu_\ell}\subset\tilde\Fc_s^\ell$, be a filtration on $(\tilde\Omega,\tilde\Fc,\Q)$ independent of $\tilde B$.
\item Let $(t,x,\tilde\xi)\in[0,T]\times\R^n\times L^2(\tilde\Omega,\tilde\Fc,\Q;\R^n)$, where $\tilde\xi$ is $\tilde\Fc_t^\ell$-measurable for every $\ell\in\N$ and $\pi=\P_{\tilde\xi}$ under $\Q$.
\end{itemize}
For every $\ell\in\N$, consider the system of equations:
\begin{align*}
d\tilde X_s^{t,\tilde\xi,\ell} \ &= \ b\big(s,\tilde X_s^{t,\tilde\xi,\ell},\P_{\text{\tiny$\tilde X_s^{t,\tilde\xi,\ell}$}}^{\text{\tiny$\tilde\Fc_s^{\mu_\ell}$}},\tilde I_s^\ell\big)\,ds + \sigma\big(s,\tilde X_s^{t,\tilde\xi,\ell},\P_{\text{\tiny$\tilde X_s^{t,\tilde\xi,\ell}$}}^{\text{\tiny$\tilde\Fc_s^{\mu_\ell}$}},\tilde I_s^\ell\big)\,d\tilde B_s, \qquad &\tilde X_t^{t,\tilde\xi,\ell} \ = \ \tilde\xi, \\
d\tilde X_s^{t,x,\pi,\ell} \ &= \ b\big(s,\tilde X_s^{t,x,\pi,\ell},\P_{\text{\tiny$\tilde X_s^{t,\tilde\xi,\ell}$}}^{\text{\tiny$\tilde\Fc_s^{\mu_\ell}$}},\tilde I_s^\ell\big)\,ds + \sigma\big(s,\tilde X_s^{t,x,\pi,\ell},\P_{\text{\tiny$\tilde X_s^{t,\tilde\xi,\ell}$}}^{\text{\tiny$\tilde\Fc_s^{\mu_\ell}$}},\tilde I_s^\ell\big)\,d\tilde B_s, &\tilde X_t^{t,x,\pi,\ell} \ = \ x,
\end{align*}
for all $s\in[t,T]$, where $(\tilde I_s^\ell)_{s\in[t,T]}$ is an $A$-valued $\tilde\F^\ell$-progressive process. Then
\begin{align*}
&\E^\Q\bigg[\int_t^T f\big(s,\tilde X_s^{t,x,\pi,\ell},\P_{\text{\tiny$\tilde X_s^{t,\tilde\xi,\ell}$}}^{\text{\tiny$\tilde\Fc_s^{\mu_\ell}$}},\tilde I_s^\ell\big)\,ds + g\big(\tilde X_T^{t,x,\pi,\ell},\P_{\text{\tiny$\tilde X_T^{t,\tilde\xi,\ell}$}}^{\text{\tiny$\tilde\Fc_T^{\mu_\ell}$}}\big)\bigg] \\
&\overset{\ell\rightarrow\infty}{\longrightarrow} \ \E^\Q\bigg[\int_t^T f\big(s,\tilde X_s^{t,x,\pi,0},\P_{\text{\tiny$\tilde X_s^{t,\tilde\xi,0}$}}^{\text{\tiny$\tilde\Fc_s^{\mu_0}$}},\tilde I_s^0\big)\,ds + g\big(\tilde X_T^{t,x,\pi,0},\P_{\text{\tiny$\tilde X_T^{t,\tilde\xi,0}$}}^{\text{\tiny$\tilde\Fc_T^{\mu_0}$}}\big)\bigg].
\end{align*}
whenever $\tilde\rho^\Q(\tilde I^\ell,\tilde I^0):=\E^\Q[\int_0^T\rho(\tilde I_s^\ell,\tilde I_s^0)\,ds]\rightarrow0$ as $\ell\rightarrow\infty$.
\end{Lemma}
\textbf{Proof.}
We begin noting that, by standard arguments (based on the Burkholder-Davis-Gundy and Gronwall inequalities), we have
\begin{equation}\label{EstimateRand_m}
\sup_{\ell\in\N}\E^\Q\Big[\sup_{s\in[t,T]}\big(\big|\tilde X_s^{t,\tilde\xi,\ell}\big|^2 + \big|\tilde X_s^{t,x,\pi,\ell}\big|^q\big)\Big] \ < \ \infty,
\end{equation}
for all $q\geq1$. We also have
\begin{align}\label{Xm-X}
\E^\Q\Big[\sup_{s\in[t,T]}\big|\tilde X_s^{t,\tilde\xi,\ell} - \tilde X_s^{t,\tilde\xi,0}\big|^2\Big] \ &\leq \ C\,\E^\Q\bigg[\int_t^T \big(\big|b\big(s,\tilde X_s^{t,\tilde\xi,0},\P_{\text{\tiny$\tilde X_s^{t,\tilde\xi,0}$}}^{\text{\tiny$\tilde\Fc_s^{\mu_0}$}},\tilde I_s^\ell\big) - b\big(s,\tilde X_s^{t,\tilde\xi,0},\P_{\text{\tiny$\tilde X_s^{t,\tilde\xi,0}$}}^{\text{\tiny$\tilde\Fc_s^{\mu_0}$}},\tilde I_s^0\big)\big|^2 \notag \\
& + \big|\sigma\big(s,\tilde X_s^{t,\tilde\xi,0},\P_{\text{\tiny$\tilde X_s^{t,\tilde\xi,0}$}}^{\text{\tiny$\tilde\Fc_s^{\mu_0}$}},\tilde I_s^\ell\big) - \sigma\big(s,\tilde X_s^{t,\tilde\xi,0},\P_{\text{\tiny$\tilde X_s^{t,\tilde\xi,0}$}}^{\text{\tiny$\tilde\Fc_s^{\mu_0}$}},\tilde I_s^0\big)\big|^2\big)\,ds\bigg],
\end{align}
for some positive constant $C$, independent of $\ell$. Now, we notice that $\tilde\rho^\Q(\tilde I^\ell,\tilde I^0)\rightarrow0$ implies $\tilde I^\ell\rightarrow\tilde I^0$ in $d\Q\,ds$-measure, which in turn implies the convergence to zero in $d\Q\,ds$-measure of the integrand in the right-hand side of \eqref{Xm-X}. By uniform integrability (which follows from \eqref{EstimateRand_m} and Assumption {\bf (A1)}(ii)), we deduce
\[
\Wc_{\text{\tiny$2$}}\big(\P_{\text{\tiny$\tilde X_s^{t,\tilde\xi,\ell}$}}^{\text{\tiny$\tilde\Fc_s^{\mu_\ell}$}},\P_{\text{\tiny$\tilde X_s^{t,\tilde\xi,0}$}}^{\text{\tiny$\tilde\Fc_s^{\mu_0}$}}\big)^2 \ \leq \ \E^\Q\Big[\big|\tilde X_s^{t,\tilde\xi,\ell} - \tilde X_s^{t,\tilde\xi,0}\big|^2\Big|\bigvee_{\ell\in\N}\tilde\Fc_\infty^{\mu_\ell}\Big] \ \overset{\ell\rightarrow\infty}{\longrightarrow} \ 0,
\]
$\Q$-a.s., for all $s\in[t,T]$. Moreover
\begin{equation}\label{Xm-->X}
\sup_{s\in[t,T]}\Wc_{\text{\tiny$2$}}\big(\P_{\text{\tiny$\tilde X_s^{t,\tilde\xi,\ell}$}}^{\text{\tiny$\tilde\Fc_s^{\mu_\ell}$}},\P_{\text{\tiny$\tilde X_s^{t,\tilde\xi,0}$}}^{\text{\tiny$\tilde\Fc_s^{\mu_0}$}}\big)^2 \ \leq \ \E^\Q\Big[\sup_{s\in[t,T]}\big|\tilde X_s^{t,\tilde\xi,\ell} - \tilde X_s^{t,\tilde\xi,0}\big|^2\Big|\bigvee_{\ell\in\N}\tilde\Fc_\infty^{\mu_\ell}\Big] \ \overset{\ell\rightarrow\infty}{\longrightarrow} \ 0.
\end{equation}
Similarly, we have
\begin{align*}
&\E^\Q\Big[\sup_{s\in[t,T]}\big|\tilde X_s^{t,x,\pi,\ell} - \tilde X_s^{t,x,\pi,0}\big|^2\Big] \ \leq \ C\,\E^\Q\bigg[\int_t^T \big(\big|b\big(s,\tilde X_s^{t,x,\pi,0},\P_{\text{\tiny$\tilde X_s^{t,\tilde\xi,\ell}$}}^{\text{\tiny$\tilde\Fc_s^{\mu_\ell}$}},\tilde I_s^\ell\big) \\
&- |b\big(s,\tilde X_s^{t,x,\pi,0},\P_{\text{\tiny$\tilde X_s^{t,\tilde\xi,0}$}}^{\text{\tiny$\tilde\Fc_s^{\mu_0}$}},\tilde I_s^0\big)\big|^2 + \big|\sigma\big(s,\tilde X_s^{t,x,\pi,0},\P_{\text{\tiny$\tilde X_s^{t,\tilde\xi,\ell}$}}^{\text{\tiny$\tilde\Fc_s^{\mu_\ell}$}},\tilde I_s^\ell\big) - \sigma\big(s,\tilde X_s^{t,x,\pi,0},\P_{\text{\tiny$\tilde X_s^{t,\tilde\xi,0}$}}^{\text{\tiny$\tilde\Fc_s^{\mu_0}$}},\tilde I_s^0\big)\big|^2\big)\,ds\bigg].
\end{align*}
Then, by \eqref{Xm-->X}, the convergence $\tilde I^\ell\rightarrow\tilde I^0$ in $d\Q\,ds$-measure, estimate \eqref{EstimateRand_m}, and Assumption {\bf (A1)}(ii), we obtain
\begin{equation}\label{Xm-->X_2}
\E^\Q\Big[\sup_{s\in[t,T]}\big|\tilde X_s^{t,x,\pi,\ell} - \tilde X_s^{t,x,\pi,0}\big|^2\Big] \ \overset{\ell\rightarrow\infty}{\longrightarrow} \ 0.
\end{equation}
Then, by \eqref{Xm-->X} and \eqref{Xm-->X_2}, we see that $f(s,\tilde X_s^{t,x,\pi,\ell},\P_{\text{\tiny$\tilde X_s^{t,\tilde\xi,\ell}$}}^{\text{\tiny$\tilde\Fc_s^{\mu_\ell}$}},\tilde I_s^\ell)\rightarrow f(s,\tilde X_s^{t,x,\pi,0},\P_{\text{\tiny$\tilde X_s^{t,\tilde\xi,0}$}}^{\text{\tiny$\tilde\Fc_s^{\mu_0}$}},\tilde I_s^0)$ as $\ell\rightarrow\infty$ in $d\Q\,ds$-measure. Therefore, by uniform integrability (which follows from estimate \eqref{EstimateRand_m} and Assumption {\bf (A1)}(ii)), we deduce
\[
\E^\Q\bigg[\int_t^T f\big(s,\tilde X_s^{t,x,\pi,\ell},\P_{\text{\tiny$\tilde X_s^{t,\tilde\xi,\ell}$}}^{\text{\tiny$\tilde\Fc_s^{\mu_\ell}$}},\tilde I_s^\ell\big)\,ds\bigg] \ \overset{\ell\rightarrow\infty}{\longrightarrow} \ \E^\Q\bigg[\int_t^T f\big(s,\tilde X_s^{t,x,\pi,0},\P_{\text{\tiny$\tilde X_s^{t,\tilde\xi,0}$}}^{\text{\tiny$\tilde\Fc_s^{\mu_0}$}},\tilde I_s^0\big)\,ds\bigg].
\]
Using again \eqref{Xm-->X} and \eqref{Xm-->X_2}, we obtain the $\Q$-a.s. pointwise convergence $g(\tilde X_T^{t,x,\pi,\ell},\P_{\text{\tiny$\tilde X_T^{t,\tilde\xi,\ell}$}}^{\text{\tiny$\tilde\Fc_T^{\mu_\ell}$}})\rightarrow g(\tilde X_T^{t,x,\pi,0},\P_{\text{\tiny$\tilde X_T^{t,\tilde\xi,0}$}}^{\text{\tiny$\tilde\Fc_T^{\mu_0}$}})$ as $\ell\rightarrow\infty$. By estimate \eqref{Estimate} together with the polynomial growth condition of $g$ in Assumption {\bf (A1)}(ii), we can apply Lebesgue's dominated convergence theorem and obtain
\[
\E^\Q\Big[g\big(\tilde X_T^{t,x,\pi,\ell},\P_{\text{\tiny$\tilde X_T^{t,\tilde\xi,\ell}$}}^{\text{\tiny$\tilde\Fc_T^{\mu_\ell}$}}\big)\Big] \ \overset{\ell\rightarrow\infty}{\longrightarrow} \ \E^\Q\Big[g\big(\tilde X_T^{t,x,\pi,0},\P_{\text{\tiny$\tilde X_s^{t,\tilde\xi,0}$}}^{\text{\tiny$\tilde\Fc_s^{\mu_0}$}}\big)\Big],
\]
which concludes the proof.
\ep

\subsection{On a different randomization of the control}
\label{App:Randomization}

\setcounter{Theorem}{0}
\setcounter{Definition}{0}
\setcounter{Proposition}{0}
\setcounter{Assumption}{0}
\setcounter{Lemma}{0}
\setcounter{Corollary}{0}
\setcounter{Remark}{0}
\setcounter{Example}{0}
\setcounter{equation}{0}

In the present appendix we introduce, following \cite{KP15}, a different kind of randomization, which in our paper turns out to be useful in the proof of Theorem \ref{Thm:Feynman-Kac}. More precisely, for every $t\in[0,T]$, $a_0\in A$, consider the $A$-valued piecewise constant process $\bar I^{t,a_0}=(\bar I_s^{t,a_0})_{s\geq t}$ on $(\bar\Omega,\bar\Fc,\bar\P)$ given by:
\begin{equation}
\label{Ihat_a}
\bar I_s^{t,a_0}(\omega,\omega^1) = \sum_{\substack{n\geq 0 \\ t<T_{n+1}(\omega^1)}} \big(a_0 1_{\{T_n(\omega^1)<t\}} + (\Ac_n(\omega^1))_{s\wedge T}(\omega) 1_{\{t\leq T_n(\omega^1)\}}\big)\,1_{[T_n(\omega^1),T_{n+1}(\omega^1))}(s),
\end{equation}
for all $s\geq t$, where we recall that $T_0=0$ and $\Ac_0=\bar\alpha$. The process $\bar I=(\bar I_s)_{s\geq0}$ defined in \eqref{barI} corresponds to $\bar I^{0,a_0}=(\bar I_s^{0,a_0})_{s\geq0}$, for any $a_0\in A$ (when $t=0$, $a_0$ plays no role in \eqref{Ihat_a}).

Let $\bar\F^{B,t}=(\bar\Fc_s^{B,t})_{s\geq t}$ (resp. $\bar\F^{\mu,t}=(\bar\Fc_s^{\mu,t})_{s\geq t}$) be the $\P$-completion of the filtration generated by $(\bar B_s-\bar B_t)_{s\geq t}$ (resp. $\bar\mu\,1_{(t,\infty)\times\Ac}$), and let $\bar\F^{B,\mu,t}=(\bar\Fc_s^{B,\mu,t})_{s\geq t}$ denote the $\P$-completion of the filtration generated by $(\bar B_s-\bar B_t)_{s\geq t}$ and $\bar\mu\,1_{(t,\infty)\times\Ac}$. If we randomize the control in \eqref{StateEq1}-\eqref{StateEq2} by means of the process $\bar I^{t,a_0}$, we obtain, for every $(x,\bar\xi)\in\R^n\times L^2(\bar\Omega,\bar\Gc,\bar\P;\R^n)$, with $\pi=\P_{\text{\tiny$\xi$}}$ under $\bar\P$:
\begin{align}
d\bar X_s^{t,\bar\xi,a_0} \ &= \ b\big(s,\bar X_s^{t,\bar\xi,a_0},\P_{\text{\tiny$\bar X_s^{t,\bar\xi,a_0}$}}^{\text{\tiny$\bar\Fc_s^{\mu,t}$}},\bar I_s^{t,a_0}\big)\,ds + \sigma\big(s,\bar X_s^{t,\bar\xi,a_0},\P_{\text{\tiny$\bar X_s^{t,\bar\xi,a_0}$}}^{\text{\tiny$\bar\Fc_s^{\mu,t}$}},\bar I_s^{t,a_0}\big)\,d\bar B_s, \label{StateEq1Rand_a} \\
d\bar X_s^{t,x,\pi,a_0} \ &= \ b\big(s,\bar X_s^{t,x,\pi,a_0},\P_{\text{\tiny$\bar X_s^{t,\bar\xi,a_0}$}}^{\text{\tiny$\bar\Fc_s^{\mu,t}$}},\bar I_s^{t,a_0}\big)\,ds + \sigma\big(s,\bar X_s^{t,x,\pi,a_0},\P_{\text{\tiny$\bar X_s^{t,\bar\xi,a_0}$}}^{\text{\tiny$\bar\Fc_s^{\mu,t}$}},\bar I_s^{t,a_0}\big)\,d\bar B_s, \label{StateEq2Rand_a}
\end{align}
for all $s\in[t,T]$, with $\bar X_t^{t,\bar\xi,a_0}=\bar\xi$ and $\bar X_t^{t,x,\pi,a_0}=x$. Under Assumption {\bf (A1)}, there exists a unique (up to indistinguishability) pair $(\bar X_s^{t,\bar\xi,a_0},\bar X_s^{t,x,\pi,a_0})_{s\in[t,T]}$ of continuous $(\bar\Fc_s^{B,\mu,t}\vee\bar\Gc)_s$-adapted processes solution to equations \eqref{StateEq1Rand_a}-\eqref{StateEq2Rand_a}, satisfying
\[
\bar\E\Big[\sup_{s\in[t,T]}\big(\big|\bar X_s^{t,\bar\xi,a_0}\big|^2 + \big|\bar X_s^{t,x,\pi,a_0}\big|^q\big)\Big] \ < \ \infty,
\]
for all $q\geq1$.

Let $\F^{\mu,t}=(\Fc_s^{\mu,t})_{s\geq t}$ be the $\P^1$-completion of the filtration generated by $\mu\,1_{(t,\infty)\times\Ac_{\textup{\tiny{step}}}}$, and denote by $\Pc(\F^{\mu,t})$ the predictable $\sigma$-algebra on $\Omega^1\times[t,\infty)$ corresponding to $\F^{\mu,t}$. Then, we define $\Vc_t$ as the set of $\Pc(\F^{\mu,t})\otimes\Bc(\Ac)$-measurable maps $\nu\colon\Omega^1\times[t,\infty)\times\Ac\rightarrow(0,\infty)$, with $0<\inf_{\Omega^1\times[t,\infty)\times\Ac}\nu\leq\sup_{\Omega^1\times[t,\infty)\times\Ac}\nu<\infty$. Given $\nu\in\Vc_t$, we define $\nu^*\in\Vc$ as $\nu^*=1_{\Omega^1\times[0,t)\times\Ac} + \nu\,1_{\Omega^1\times[t,\infty)\times\Ac}$. We denote $\P^\nu$ (resp. $\bar\P^\nu$) the probability $\P^{\nu^*}$ (resp. $\bar\P^{\bar\nu^*}$), and $\E^\nu$ (resp. $\bar\E^\nu$) the expectation $\E^{\nu^*}$ (resp. $\bar\E^{\bar\nu^*}$). Then, for every $\nu\in\Vc_t$, we define the gain functional (notice that $J^\Rc(t,x,\pi,a_0,\nu)$ does not depend on the value of $\nu^*$ on $\Omega^1\times[0,t)\times\Ac$)
\[
J^\Rc(t,x,\pi,a_0,\nu) \ = \ \bar\E^\nu\bigg[\int_t^T f\big(s,\bar X_s^{t,x,\pi,a_0},\P_{\text{\tiny$\bar X_s^{t,\bar\xi,a_0}$}}^{\text{\tiny$\bar\Fc_s^{\mu,t}$}},\bar I_s^{t,a_0}\big)\,ds + g\big(\bar X_T^{t,x,\pi,a_0},\P_{\text{\tiny$\bar X_T^{t,\bar\xi,a_0}$}}^{\text{\tiny$\bar\Fc_T^{\mu,t}$}}\big)\bigg]
\]
and the value function
\[
V^\Rc(t,x,\pi,a_0) \ = \ \sup_{\nu\in\Vc_t} J^\Rc(t,x,\pi,a_0,\nu).
\]
Finally, let $\F^{B,t}=(\Fc_s^{B,t})_{s\geq t}$ be the $\P$-completion of the filtration generated by $(B_s-B_t)_{s\geq t}$, and let $\Ac_t$ denote the set of $\F^{B,t}$-progressive processes $\alpha\colon\Omega\times[t,T]\rightarrow A$. Given $\alpha\in\Ac_t$, we define $\alpha^*\in\Ac$ as $\alpha^*=\bar a\,1_{\Omega\times[0,t)}+\alpha\,1_{\Omega\times[t,T]}$, for some deterministic and fixed point $\bar a\in A$. Then, we denote $J(t,x,\pi,\alpha^*)$ simply by $J(t,x,\pi,\alpha)$ (notice that $J(t,x,\pi,\alpha^*)$ does not depend on the value of $\alpha^*$ on $\Omega\times[0,t)$, namely on $\bar a$).

\begin{Theorem}\label{Thm:Equivalence2}
Under Assumption {\bf (A1)}, we have the following identities:
\begin{align}\label{Equivalence_a}
V(t,x,\pi) \ := \ \sup_{\alpha\in\Ac} J(t,x,\pi,\alpha) \ &= \ \sup_{\alpha\in\Ac_t} J(t,x,\pi,\alpha) \ = \ \sup_{\nu\in\Vc_t} J^\Rc(t,x,\pi,a_0,\nu) \ =: \ V^\Rc(t,x,\pi,a_0) \notag \\
&= \ \sup_{\nu\in\Vc} J^\Rc(t,x,\pi,\nu) \ =: \ V^\Rc(t,x,\pi),
\end{align}
for all $(t,x,\pi,a_0)\in[0,T]\times\R^n\times\mathscr P_{\text{\tiny$2$}}(\R^n)\times A$.
\end{Theorem}

\begin{Remark}
{\rm
From Theorem \ref{Thm:Equivalence} we conclude that the function $V^\Rc(t,x,\pi,a_0)$ does not depend on $a_0\in A$ and coincides with the function $V^\Rc(t,x,\pi)$ defined in \eqref{ValueRand}.
\ep
}
\end{Remark}

\noindent\textbf{Proof.}
When $t=0$, we see that, for every $a_0\in A$, we have $\bar I^{0,a_0}=\bar I$, $\Ac_0=\Ac$, and $\Vc_0=\Vc$. Therefore, $V^\Rc(0,x,\pi,a_0)$ coincides with $V^\Rc(0,x,\pi)$, so the result follows from Theorem \ref{Thm:Equivalence}. When $t>0$, we proceed along the same lines as in the proof of Theorem \ref{Thm:Equivalence} for the case $t=0$, with $(\bar B_s)_{s\geq0}$, $\bar\F^B=(\bar\Fc_s^B)_{s\geq0}$, $\Ac$, $\bar\mu$, $\bar\F^{B,\mu}=(\bar\Fc_s^{B,\mu})_{s\geq0}$, $\Vc$ replaced respectively by $(\bar B_s-\bar B_t)_{s\geq t}$, $\bar\F^{B,t}=(\bar\Fc_s^{B,t})_{s\geq t}$, $\Ac_t$, $\bar\mu\,1_{(t,\infty)\times\Ac}$, $\bar\F^{B,\mu,t}=(\bar\Fc_s^{B,\mu,t})_{s\geq t}$, $\Vc_t$. Then, we obtain
\[
\sup_{\alpha\in\Ac_t} J(t,x,\pi,\alpha) \ = \ \sup_{\nu\in\Vc_t} J^\Rc(t,x,\pi,a_0,\nu).
\]
This implies that $V^\Rc(t,x,\pi,a_0)$ does not depend on $a_0\in A$, since the left-hand side of the above inequality does not depend on it.

By Theorem \ref{Thm:Equivalence}, equivalence \eqref{Equivalence_a} follows if we prove the following inequalities
\begin{align}\label{TwoIneq}
V(t,x,\pi) \ \geq \ \sup_{\alpha\in\Ac_t} J(t,x,\pi,\alpha), \qquad\qquad \sup_{\nu\in\Vc_t} J^\Rc(t,x,\pi,a_0,\nu) \ &\geq \ V^\Rc(t,x,\pi).
\end{align}
Since for every $\alpha\in\Ac_t$ we have, by definition, $J(t,x,\pi,\alpha)=J(t,x,\pi,\alpha^*)$, where $\alpha^*=\bar a\,1_{\Omega\times[0,t)}+\alpha\,1_{\Omega\times[t,T]}$, we see that $\sup_{\alpha\in\Ac_t} J(t,x,\pi,\alpha)\leq\sup_{\alpha\in\Ac} J(t,x,\pi,\alpha)=V(t,x,\pi)$. Therefore, the first inequality in \eqref{TwoIneq} is proved.

In order to establish the second inequality in \eqref{TwoIneq}, we fix $(t,x,\bar\xi,\pi)\in[0,T]\times\R^n\times L^2(\bar\Omega,\bar\Gc,\bar\P;\R^n)\times\mathscr P_{\text{\tiny$2$}}(\R^n)$, with $\pi=\P_{\text{\tiny$\xi$}}$ under $\P$, and we take a particular probabilistic setting for the randomized McKean-Vlasov control problem. More precisely, we first consider another probabilistic framework for randomized problem, where the objects $(\Omega,\Fc,\P)$, $(\Omega^1,\Fc^1,\P^1)$, $(\bar\Omega,\bar\Fc,\bar\P)$, $\bar B$, $\bar\mu$, $(T_n,\Ac_n)$, $\bar I$ are replaced respectively by $(\Omega^0,\Fc^0,\P^0)$, $(\check\Omega^1,\check\Fc^1,\check\P^1)$, $(\check\Omega,\check\Fc,\check\P)$, $\check B$, $\check\mu$, $(\check T_n,\check\Ac_n)$, $\check I$.

Let $\hat\Omega=\check\Omega\times\bar\Omega$, $\hat\Fc$ the $\check\P\otimes\bar\P$-completion of $\check\Fc\otimes\bar\Fc$, $\hat\P$ the extension of $\check\P\otimes\bar\P$ to $\hat\Fc$, and $\hat\E$ the $\hat\P$-expected value. Let also $\hat\Gc$ be the canonical extension of $\bar\Gc$ to $\hat\Omega$. Define $\hat\xi(\check\omega,\bar\omega):=\bar\xi(\bar\omega)$ and
\begin{align*}
\hat B_s(\check\omega,\bar\omega) \ &:= \ \check B_s(\check\omega)\,1_{\{s\leq t\}} + (\bar B_s(\bar\omega) - \bar B_t(\bar\omega) + \check B_t(\check\omega))\,1_{\{s>t\}}, \\
\hat\mu(\check\omega,\bar\omega;ds\,d\alpha) \ &:= \ \check\mu(\check\omega;ds\,d\alpha)\,1_{\{s\leq t\}} + \bar\mu(\bar\omega;ds\,d\alpha)\,1_{\{s>t\}}.
\end{align*}
Notice that $\pi=\P_{\text{\tiny$\hat\xi$}}$ under $\hat\P$, $\hat B=(\hat B_s)_{s\geq 0}$ is a Brownian motion on $(\hat\Omega,\hat\Fc,\hat\P)$, and $\hat\mu$ is a Poisson random measure with compensator $\lambda(d\alpha)\,ds$ under $\hat\P$, with respect to its natural filtration. We also define as in \eqref{barI} the $A$-valued piecewise constant process $\hat I=(\hat I_s)_{s\geq0}$ associated to $\hat\mu$, which in the present case takes the following form:
\begin{align*}
\hat I_s(\check\omega,\bar\omega) \ &= \ \check I_s(\check\omega)\,1_{\{s\leq t\}} \\
&\quad \ + \!\!\!\!\!\! \sum_{\substack{n\geq 0 \\ t<T_{n+1}(\omega^1)}} \!\!\! \big(\check I_s(\check\omega) 1_{\{T_n(\omega^1)<t\}} + (\Ac_n(\omega^1))_{s\wedge T}(\omega) 1_{\{t\leq T_n(\omega^1)\}}\big)\,1_{[T_n(\omega^1),T_{n+1}(\omega^1))}(s)\,1_{\{s>t\}}.
\end{align*}
In particular, $\hat I_t=\check I_t$. We define $\hat\F^{B,\mu}=(\hat\Fc_s^{B,\mu})_{s\geq0}$ (resp. $\hat\F^\mu=(\hat\Fc_s^\mu)_{s\geq0}$) as the $\hat\P$-completion of the filtration generated by $\hat B$ and $\hat\mu$ (resp. $\hat\mu$). We denote $(\hat X_s^{t,\hat\xi},\hat X_s^{t,x,\pi})_{s\in[t,T]}$ the unique (up to indistinguishability) continuous $(\hat\Fc_s^{B,\mu}\vee\hat\Gc)_s$-adapted solution to equations \eqref{StateEq1Rand}-\eqref{StateEq2Rand} on $(\hat\Omega,\hat\Fc,\hat\P)$ with $\bar \xi$, $\bar B$, $\bar I$, $\bar\Fc_\cdot^\mu$ replaced respectively by $\hat\xi$, $\hat B$, $\hat I$, $\hat\Fc_\cdot^\mu$. For later use, we also consider, for every $\check\omega\in\check\Omega$, the unique (up to indistinguishability) continuous $(\bar\Fc_s^{B,\mu,t}\vee\bar\Gc)_s$-adapted solution $(\bar X_s^{t,\bar\xi,\check I_t(\check\omega)},\bar X_s^{t,x,\pi,\check I_t(\check\omega)})_{s\in[t,T]}$ to equations \eqref{StateEq1Rand_a}-\eqref{StateEq2Rand_a} with $a_0$ replaced by $\check I_t(\check\omega)$. Then, we see that, for $\check\P$-a.e. $\check\omega\in\check\Omega$, $(\hat X_s^{t,\hat\xi}(\check\omega,\cdot),\hat X_s^{t,x,\pi}(\check\omega,\cdot))_{s\in[t,T]}$ and $(\bar X_s^{t,\bar\xi,\check I_t(\check\omega)},\bar X_s^{t,x,\pi,\check I_t(\check\omega)})_{s\in[t,T]}$ solve the same system of equations. Therefore, by pathwise uniqueness, for $\check\P$-a.e. $\check\omega\in\check\Omega$, we have $\hat X_s^{t,\hat\xi}(\check\omega,\bar\omega)=\bar X_s^{t,\bar\xi,\check I_t(\check\omega)}(\bar\omega)$ and $\hat X_s^{t,x,\pi}(\check\omega,\bar\omega)=\bar X_s^{t,x,\pi,\check I_t(\check\omega)}(\bar\omega)$, for all $s\in[t,T]$, $\bar\P(d\bar\omega)$-almost surely.

Let $\Pc(\hat\F^\mu)$ be the predictable $\sigma$-algebra on $\hat\Omega\times\R_+$ corresponding to $\hat\F^\mu$. In order to define the randomized McKean-Vlasov control problem on $(\hat\Omega,\hat\Fc,\hat\P)$, we introduce the set $\hat\Vc$ of all $\Pc(\hat\F^\mu)\otimes\Bc(\Ac)$-measurable maps $\hat\nu\colon\hat\Omega\times\R_+\times\Ac\rightarrow(0,\infty)$, satisfying $0<\inf_{\hat\Omega\times\R_+\times\Ac}\hat\nu\leq\sup_{\hat\Omega\times\R_+\times\Ac}\hat\nu<\infty$. Then, we define in an obvious way $\kappa^{\hat\nu}$, $\hat\P^{\hat\nu}$, $\hat\E^{\hat\nu}$, $\hat J^\Rc(t,x,\pi,\hat\nu)$, and the corresponding value function $\hat V^\Rc(t,x,\pi)$. We recall from step I of the proof of Theorem \ref{Thm:Equivalence} that $\hat V^\Rc(t,x,\pi)=V^\Rc(t,x,\pi)$.

We can now prove the second inequality in \eqref{TwoIneq}, namely
\begin{equation}\label{TwoIneq_Proof}
V^\Rc(t,x,\pi) \ = \ \hat V^\Rc(t,x,\pi) \ := \ \sup_{\hat\nu\in\hat\Vc}\hat J^\Rc(t,x,\pi,\hat\nu) \ \leq \ \sup_{\nu\in\Vc_t} J^\Rc(t,x,\pi,a_0,\nu).
\end{equation}
Fix $\hat\nu\in\hat\Vc$. We begin noting that, since $\hat\nu$ is $\Pc(\hat\F^\mu)\otimes\Bc(\Ac)$-measurable, up to a $\hat\P$-null set, $\hat\nu$ depends only $(\check\omega^1,\omega^1)$. Now, by a monotone class argument, we see that there exists a $\check\P^1$-null set $\check N^1\in\check\Fc^1$ such that $\nu^{\check\omega^1}=\nu_s^{\check\omega^1}(\omega^1,\alpha)\colon\Omega^1\times[t,\infty)\times\Ac\rightarrow(0,\infty)$, given by
\[
\nu_s^{\check\omega^1}(\omega^1,\alpha) \ := \ \hat\nu_s(\check\omega^1,\omega^1,\alpha), \qquad \text{for all }(\check\omega^1,\omega^1,s,\alpha)\in\check\Omega^1\times\Omega^1\times[t,\infty)\times\Ac,
\]
is an element of $\Vc_t$, for every $\check\omega^1\notin\check N^1$. In other words, for every $\check\omega^1\notin\check N^1$, $\nu^{\check\omega^1}$ is a $\Pc(\F^{\mu,t})\otimes\Bc(\Ac)$-measurable map satisfying $0<\inf_{\Omega^1\times[t,\infty)\times\Ac}\nu^{\check\omega^1}\leq\sup_{\Omega^1\times[t,\infty)\times\Ac}\nu^{\check\omega^1}<\infty$. Therefore, by Fubini's theorem,
\begin{align*}
&\hat J^\Rc(t,x,\pi,\hat\nu) \ = \ \hat\E\bigg[\kappa_T^{\hat\nu}\bigg(\int_t^T f\big(s,\hat X_s^{t,x,\pi},\P_{\text{\tiny$\hat X_s^{t,\hat\xi}$}}^{\text{\tiny$\hat\Fc_s^\mu$}},\hat I_s\big)\,ds + g\big(\hat X_T^{t,x,\pi},\P_{\text{\tiny$\hat X_T^{t,\hat\xi}$}}^{\text{\tiny$\hat\Fc_T^\mu$}}\big)\bigg)\bigg] \\
&= \ \int_{\check\Omega} \bar\E\bigg[\kappa_T^{\nu^{\check\omega^1}}\bigg(\int_t^T f\big(s,\bar X_s^{t,x,\pi,\check I_t(\check\omega)},\P_{\text{\tiny$\bar X_s^{t,\bar\xi,\check I_t(\check\omega)}$}}^{\text{\tiny$\bar\Fc_s^\mu$}},\bar I_s^{t,\check I_t(\check\omega)}\big)\,ds + g\big(\bar X_T^{t,x,\pi,\check I_t(\check\omega)},\P_{\text{\tiny$\bar X_T^{t,\bar\xi,\check I_t(\check\omega)}$}}^{\text{\tiny$\bar\Fc_T^\mu$}}\big)\bigg)\bigg]\,\check\P(d\check\omega) \\
&= \ \int_{\check\Omega} J^\Rc(t,x,\pi,\check I_t(\check\omega),\nu^{\check\omega^1})\,\check\P(d\check\omega) \ \leq \ \sup_{\nu\in\Vc_t} J^\Rc(t,x,\pi,a_0,\nu),
\end{align*}
for any $a_0\in A$ (recall that $\sup_{\nu\in\Vc_t} J^\Rc(t,x,\pi,a_0,\nu)$ does not depend on $a_0\in A$). From the arbitrariness of $\hat\nu\in\hat\Vc$, we deduce that $\sup_{\hat\nu\in\hat\Vc}\hat J^\Rc(t,x,\pi,\hat\nu)\leq\sup_{\nu\in\Vc_t} J^\Rc(t,x,\pi,a_0,\nu)$, hence establishing \eqref{TwoIneq_Proof}, and consequently the second inequality in \eqref{TwoIneq}.
\ep

\begin{Corollary}\label{C:EquivConditional}
Under Assumption {\bf (A1)}, we have
\begin{equation}\label{EquivConditional}
V(t,x,\pi) \ = \ \esssup_{\nu\in\Vc} \E^\nu\bigg[\int_t^T \E\big[f\big(s,\bar X_s^{t,x,\pi},\P_{\text{\tiny$\bar X_s^{t,\bar\xi}$}}^{\text{\tiny$\bar\Fc_s^\mu$}},\bar I_s\big)\big]\,ds + \E\big[g\big(\bar X_T^{t,x,\pi},\P_{\text{\tiny$\bar X_T^{t,\bar\xi}$}}^{\text{\tiny$\bar\Fc_T^\mu$}}\big)\big]\bigg|\Fc_t^\mu\bigg],
\end{equation}
$\P^1$-a.s., for all $(t,x,\bar\xi)\in[0,T]\times\R^n\times L^2(\bar\Omega,\bar\Gc,\bar\P;\R^n)$, with $\pi=\P_{\text{\tiny$\bar\xi$}}$ under $\bar\P$.
\end{Corollary}
\textbf{Proof.}
Fix $(t,x,\bar\xi)\in[0,T]\times\R^n\times L^2(\bar\Omega,\bar\Gc,\bar\P;\R^n)$, with $\pi=\P_{\text{\tiny$\xi$}}$ under $\bar\P$. We have
\begin{align*}
&\E^1\bigg[\esssup_{\nu\in\Vc} \E^\nu\bigg[\int_t^T \E\big[f\big(s,\bar X_s^{t,x,\pi},\P_{\text{\tiny$\bar X_s^{t,\bar\xi}$}}^{\text{\tiny$\bar\Fc_s^\mu$}},\bar I_s\big)\big]\,ds + \E\big[g\big(\bar X_T^{t,x,\pi},\P_{\text{\tiny$\bar X_T^{t,\bar\xi}$}}^{\text{\tiny$\bar\Fc_T^\mu$}}\big)\big]\bigg|\Fc_t^\mu\bigg]\bigg] \\
&\geq \ \E^1\bigg[\esssup_{\nu\in\Vc_{1,t}} \E^\nu\bigg[\int_t^T \E\big[f\big(s,\bar X_s^{t,x,\pi},\P_{\text{\tiny$\bar X_s^{t,\bar\xi}$}}^{\text{\tiny$\bar\Fc_s^\mu$}},\bar I_s\big)\big]\,ds + \E\big[g\big(\bar X_T^{t,x,\pi},\P_{\text{\tiny$\bar X_T^{t,\bar\xi}$}}^{\text{\tiny$\bar\Fc_T^\mu$}}\big)\big]\bigg|\Fc_t^\mu\bigg]\bigg] \\
&\geq \ \sup_{\nu\in\Vc_{1,t}} \E^1\bigg[\E^\nu\bigg[\int_t^T \E\big[f\big(s,\bar X_s^{t,x,\pi},\P_{\text{\tiny$\bar X_s^{t,\bar\xi}$}}^{\text{\tiny$\bar\Fc_s^\mu$}},\bar I_s\big)\big]\,ds + \E\big[g\big(\bar X_T^{t,x,\pi},\P_{\text{\tiny$\bar X_T^{t,\bar\xi}$}}^{\text{\tiny$\bar\Fc_T^\mu$}}\big)\big]\bigg|\Fc_t^\mu\bigg]\bigg].
\end{align*}
By the Bayes formula, and recalling that $\kappa_t^\nu=1$ whenever $\nu\in\Vc_{1,t}$, we obtain
\begin{align*}
&\sup_{\nu\in\Vc_{1,t}} \E^1\bigg[\E^\nu\bigg[\int_t^T \E\big[f\big(s,\bar X_s^{t,x,\pi},\P_{\text{\tiny$\bar X_s^{t,\bar\xi}$}}^{\text{\tiny$\bar\Fc_s^\mu$}},\bar I_s\big)\big]\,ds + \E\big[g\big(\bar X_T^{t,x,\pi},\P_{\text{\tiny$\bar X_T^{t,\bar\xi}$}}^{\text{\tiny$\bar\Fc_T^\mu$}}\big)\big]\bigg|\Fc_t^\mu\bigg]\bigg] \\
&= \ \sup_{\nu\in\Vc_{1,t}} \E^1\bigg[\E^1\bigg[\kappa_T^\nu\bigg(\int_t^T \E\big[f\big(s,\bar X_s^{t,x,\pi},\P_{\text{\tiny$\bar X_s^{t,\bar\xi}$}}^{\text{\tiny$\bar\Fc_s^\mu$}},\bar I_s\big)\big]\,ds + \E\big[g\big(\bar X_T^{t,x,\pi},\P_{\text{\tiny$\bar X_T^{t,\bar\xi}$}}^{\text{\tiny$\bar\Fc_T^\mu$}}\big)\big]\bigg)\bigg|\Fc_t^\mu\bigg]\bigg] \\
&= \ \sup_{\nu\in\Vc_{1,t}} \E^\nu\bigg[\int_t^T \E\big[f\big(s,\bar X_s^{t,x,\pi},\P_{\text{\tiny$\bar X_s^{t,\bar\xi}$}}^{\text{\tiny$\bar\Fc_s^\mu$}},\bar I_s\big)\big]\,ds + \E\big[g\big(\bar X_T^{t,x,\pi},\P_{\text{\tiny$\bar X_T^{t,\bar\xi}$}}^{\text{\tiny$\bar\Fc_T^\mu$}}\big)\big]\bigg] \ = \ V(t,x,\pi),
\end{align*}
where the last equality follows from Remark \ref{R:Vc1,t}. Then, we conclude that
\begin{equation}\label{IneqConditional1}
\E^1\bigg[\esssup_{\nu\in\Vc} \E^\nu\bigg[\int_t^T \E\big[f\big(s,\bar X_s^{t,x,\pi},\P_{\text{\tiny$\bar X_s^{t,\bar\xi}$}}^{\text{\tiny$\bar\Fc_s^\mu$}},\bar I_s\big)\big]\,ds + \E\big[g\big(\bar X_T^{t,x,\pi},\P_{\text{\tiny$\bar X_T^{t,\bar\xi}$}}^{\text{\tiny$\bar\Fc_T^\mu$}}\big)\big]\bigg|\Fc_t^\mu\bigg]\bigg] \ \geq \ V(t,x,\pi).
\end{equation}
Let us now prove the following inequality: for every $\nu\in\Vc$, $\P^1$-a.s.,
\begin{equation}\label{IneqConditional2}
\E^\nu\bigg[\int_t^T \E\big[f\big(s,\bar X_s^{t,x,\pi},\P_{\text{\tiny$\bar X_s^{t,\bar\xi}$}}^{\text{\tiny$\bar\Fc_s^\mu$}},\bar I_s\big)\big]\,ds + \E\big[g\big(\bar X_T^{t,x,\pi},\P_{\text{\tiny$\bar X_T^{t,\bar\xi}$}}^{\text{\tiny$\bar\Fc_T^\mu$}}\big)\big]\bigg|\Fc_t^\mu\bigg] \ \leq \ V(t,x,\pi).
\end{equation}
Suppose we have already proved \eqref{IneqConditional2}. Hence, $\P^1$-a.s.,
\[
\esssup_{\nu\in\Vc} \E^\nu\bigg[\int_t^T \E\big[f\big(s,\bar X_s^{t,x,\pi},\P_{\text{\tiny$\bar X_s^{t,\bar\xi}$}}^{\text{\tiny$\bar\Fc_s^\mu$}},\bar I_s\big)\big]\,ds + \E\big[g\big(\bar X_T^{t,x,\pi},\P_{\text{\tiny$\bar X_T^{t,\bar\xi}$}}^{\text{\tiny$\bar\Fc_T^\mu$}}\big)\big]\bigg|\Fc_t^\mu\bigg] \ \leq \ V(t,x,\pi).
\]
From the above inequality and \eqref{IneqConditional1}, it is then easy to see that equality \eqref{EquivConditional} holds. It remains to prove \eqref{IneqConditional2}. To this end, we notice that \eqref{IneqConditional2} holds if and only if the following inequality holds: for every $\nu\in\Vc$, $\bar\P$-a.s.,
\begin{equation}\label{IneqConditional2_bis}
\bar\E^\nu\bigg[\int_t^T f\big(s,\bar X_s^{t,x,\pi},\P_{\text{\tiny$\bar X_s^{t,\bar\xi}$}}^{\text{\tiny$\bar\Fc_s^\mu$}},\bar I_s\big)\,ds + g\big(\bar X_T^{t,x,\pi},\P_{\text{\tiny$\bar X_T^{t,\bar\xi}$}}^{\text{\tiny$\bar\Fc_T^\mu$}}\big)\bigg|\bar\Fc_t^\mu\bigg] \ \leq \ V(t,x,\pi).
\end{equation}
Now, consider the same probabilistic setting introduced in the proof of Theorem \ref{Thm:Equivalence2}: $(\hat\Omega,\hat\Fc,\hat\P)$, $\hat\Gc$, $\hat B$, $\hat\mu$, $\hat\F^{B,\mu}=(\hat\Fc_s^{B,\mu})_{s\geq0}$, $\hat\F^\mu=(\hat\Fc_s^\mu)_{s\geq0}$, $\hat I$, $\hat X^{t,\hat\xi}$, $\hat X^{t,x,\pi}$, $\hat\Vc$, $\hat\Vc_{1,t}$, $\hat\P^{\hat\nu}$, $\hat\E^{\hat\nu}$, $\hat J^\Rc(t,x,\pi,\hat\nu)$, $\hat V^\Rc(t,x,\pi)$. Observe that \eqref{IneqConditional2_bis} holds if and only if the following inequality holds: for every $\hat\nu\in\hat\Vc$, $\hat\P$-a.s.,
\begin{equation}\label{EquivConditional_hat}
\hat\E^{\hat\nu}\bigg[\int_t^T f\big(s,\hat X_s^{t,x,\pi},\P_{\text{\tiny$\hat X_s^{t,\hat\xi}$}}^{\text{\tiny$\hat\Fc_s^\mu$}},\hat I_s\big)\,ds + g\big(\hat X_T^{t,x,\pi},\P_{\text{\tiny$\hat X_T^{t,\hat\xi}$}}^{\text{\tiny$\hat\Fc_T^\mu$}}\big)\bigg|\hat\Fc_t^\mu\bigg] \ \leq \ V(t,x,\pi).
\end{equation}
Indeed, let us prove that if \eqref{EquivConditional_hat} holds then \eqref{IneqConditional2_bis} holds as well (the other implication has a similar proof). Fix $\nu\in\Vc$. Then, proceeding as in step I of the proof of Theorem \ref{Thm:Equivalence}, we see that there exists $\hat\nu\in\hat\Vc$ such that
\[
\frac{\bar\kappa_T^\nu}{\bar\kappa_t^\nu}\bigg(\int_t^T f\big(s,\bar X_s^{t,x,\pi},\P_{\text{\tiny$\bar X_s^{t,\bar\xi}$}}^{\text{\tiny$\bar\Fc_s^\mu$}},\bar I_s\big)\,ds + g\big(\bar X_T^{t,x,\pi},\P_{\text{\tiny$\bar X_T^{t,\bar\xi}$}}^{\text{\tiny$\bar\Fc_T^\mu$}}\big)\bigg), \; \bar B, \; \bar \mu
\]
and
\[
\frac{\kappa_T^{\hat\nu}}{\kappa_t^{\hat\nu}}\bigg(\int_t^T f\big(s,\hat X_s^{t,x,\pi},\P_{\text{\tiny$\hat X_s^{t,\hat\xi}$}}^{\text{\tiny$\hat\Fc_s^\mu$}},\hat I_s\big)\,ds + g\big(\hat X_T^{t,x,\pi},\P_{\text{\tiny$\hat X_T^{t,\hat\xi}$}}^{\text{\tiny$\hat\Fc_T^\mu$}}\big)\bigg), \; \hat B, \; \hat \mu
\]
have the same joint law. As a consequence,
\[
\bar\E^\nu\bigg[\int_t^T f\big(s,\bar X_s^{t,x,\pi},\P_{\text{\tiny$\bar X_s^{t,\bar\xi}$}}^{\text{\tiny$\bar\Fc_s^\mu$}},\bar I_s\big)\,ds + g\big(\bar X_T^{t,x,\pi},\P_{\text{\tiny$\bar X_T^{t,\bar\xi}$}}^{\text{\tiny$\bar\Fc_T^\mu$}}\big)\bigg|\bar\Fc_t^\mu\bigg]
\]
and
\[
\hat\E^{\hat\nu}\bigg[\int_t^T f\big(s,\hat X_s^{t,x,\pi},\P_{\text{\tiny$\hat X_s^{t,\hat\xi}$}}^{\text{\tiny$\hat\Fc_s^\mu$}},\hat I_s\big)\,ds + g\big(\hat X_T^{t,x,\pi},\P_{\text{\tiny$\hat X_T^{t,\hat\xi}$}}^{\text{\tiny$\hat\Fc_T^\mu$}}\big)\bigg|\hat\Fc_t^\mu\bigg]
\]
have the same law. In particular, we have
\begin{align*}
&\bar\P\bigg(\bar\E^\nu\bigg[\int_t^T f\big(s,\bar X_s^{t,x,\pi},\P_{\text{\tiny$\bar X_s^{t,\bar\xi}$}}^{\text{\tiny$\bar\Fc_s^\mu$}},\bar I_s\big)\,ds + g\big(\bar X_T^{t,x,\pi},\P_{\text{\tiny$\bar X_T^{t,\bar\xi}$}}^{\text{\tiny$\bar\Fc_T^\mu$}}\big)\bigg|\bar\Fc_t^\mu\bigg]\leq V(t,x,\pi)\bigg) \\
&= \ \hat\P\bigg(\hat\E^{\hat\nu}\bigg[\int_t^T f\big(s,\hat X_s^{t,x,\pi},\P_{\text{\tiny$\hat X_s^{t,\hat\xi}$}}^{\text{\tiny$\hat\Fc_s^\mu$}},\hat I_s\big)\,ds + g\big(\hat X_T^{t,x,\pi},\P_{\text{\tiny$\hat X_T^{t,\hat\xi}$}}^{\text{\tiny$\hat\Fc_T^\mu$}}\big)\bigg|\hat\Fc_t^\mu\bigg]\leq V(t,x,\pi)\bigg) \ = \ 1,
\end{align*}
where the last equality follows from the assumption that \eqref{EquivConditional_hat} holds. This implies that \eqref{IneqConditional2_bis} also holds for $\nu$. Since $\nu$ was arbitrary, the claim follows.

Let us now prove that \eqref{EquivConditional_hat} holds. For every $\hat\nu\in\hat\Vc$, by the Bayes formula, and proceeding as in the proof of Theorem \ref{Thm:Equivalence2}, we find
\begin{align*}
&\hat\E^{\hat\nu}\bigg[\int_t^T f\big(s,\hat X_s^{t,x,\pi},\P_{\text{\tiny$\hat X_s^{t,\hat\xi}$}}^{\text{\tiny$\hat\Fc_s^\mu$}},\hat I_s\big)\,ds + g\big(\hat X_T^{t,x,\pi},\P_{\text{\tiny$\hat X_T^{t,\hat\xi}$}}^{\text{\tiny$\hat\Fc_T^\mu$}}\big)\bigg|\hat\Fc_t^\mu\bigg] \\
&= \ \hat\E\bigg[\frac{\kappa_T^{\hat\nu}}{\kappa_t^{\hat\nu}}\bigg(\int_t^T f\big(s,\hat X_s^{t,x,\pi},\P_{\text{\tiny$\hat X_s^{t,\hat\xi}$}}^{\text{\tiny$\hat\Fc_s^\mu$}},\hat I_s\big)\,ds + g\big(\hat X_T^{t,x,\pi},\P_{\text{\tiny$\hat X_T^{t,\hat\xi}$}}^{\text{\tiny$\hat\Fc_T^\mu$}}\big)\bigg)\bigg|\hat\Fc_t^\mu\bigg] \\
&= \ \hat\E\bigg[\frac{\kappa_T^{\nu^{\cdot}}}{\kappa_t^{\nu^{\cdot}}}\bigg(\int_t^T f\big(s,\bar X_s^{t,x,\pi,\check I_t},\P_{\text{\tiny$\bar X_s^{t,\bar\xi,\check I_t}$}}^{\text{\tiny$\bar\Fc_s^\mu$}},\bar I_s^{t,\check I_t}\big)\,ds + g\big(\bar X_T^{t,x,\pi,\check I_t},\P_{\text{\tiny$\bar X_T^{t,\bar\xi,\check I_t}$}}^{\text{\tiny$\bar\Fc_T^\mu$}}\big)\bigg)\bigg|\hat\Fc_t^\mu\bigg].
\end{align*}
Then, by the freezing lemma (see for instance Proposition 10.1.2 in \cite{Zabczyk}), we obtain
\begin{align*}
&\hat\E\bigg[\frac{\kappa_T^{\nu^{\cdot}}}{\kappa_t^{\nu^{\cdot}}}\bigg(\int_t^T f\big(s,\bar X_s^{t,x,\pi,\check I_t},\P_{\text{\tiny$\bar X_s^{t,\bar\xi,\check I_t}$}}^{\text{\tiny$\bar\Fc_s^\mu$}},\bar I_s^{t,\check I_t}\big)\,ds + g\big(\bar X_T^{t,x,\pi,\check I_t},\P_{\text{\tiny$\bar X_T^{t,\bar\xi,\check I_t}$}}^{\text{\tiny$\bar\Fc_T^\mu$}}\big)\bigg)\bigg|\hat\Fc_t^\mu\bigg] \\
&= \ \E\bigg[\kappa_T^{\nu^{\check\omega^1}}\bigg(\int_t^T f\big(s,\bar X_s^{t,x,\pi,\check I_t(\check\omega)},\P_{\text{\tiny$\bar X_s^{t,\bar\xi,\check I_t(\check\omega)}$}}^{\text{\tiny$\bar\Fc_s^\mu$}},\bar I_s^{t,\check I_t(\check\omega)}\big)\,ds + g\big(\bar X_T^{t,x,\pi,\check I_t(\check\omega)},\P_{\text{\tiny$\bar X_T^{t,\bar\xi,\check I_t(\check\omega)}$}}^{\text{\tiny$\bar\Fc_T^\mu$}}\big)\bigg)\bigg] \\
&= \ J^\Rc(t,x,\pi,\check I_t(\check\omega),\nu^{\check\omega^1}) \ \leq \ \sup_{\nu\in\Vc_t} J^\Rc(t,x,\pi,a_0,\nu),
\end{align*}
$\hat\P$-a.s., for any $a_0\in A$ (recall from Theorem \ref{Thm:Equivalence2} that $\sup_{\nu\in\Vc_t} J^\Rc(t,x,\pi,a_0,\nu)$ does not depend on $a_0\in A$). Then, since by Theorem \ref{Thm:Equivalence2} we have that $\sup_{\nu\in\Vc_t} J^\Rc(t,x,\pi,a_0,\nu)=V(t,x,\pi)$, we deduce that \eqref{EquivConditional_hat} holds, which concludes the proof.
\ep

\vspace{9mm}

\small
\bibliographystyle{plain}
\bibliography{references}

\end{document}